\documentclass[final,onefignum,onetabnum]{siamart220329}

\usepackage{amsmath}
\usepackage{amsfonts}
\usepackage{amssymb}
\usepackage{bm}
\usepackage{bbm}   
\usepackage{latexsym}
\usepackage{stmaryrd}
\usepackage{epsfig}
\usepackage{url}
\usepackage{fullpage}
\usepackage{subfigure}
\usepackage{appendix}
\usepackage{color}
\usepackage{xcolor}
\usepackage{multirow}
\usepackage{hhline}
\usepackage[ntheorem]{empheq}
\usepackage[top=2.50cm,left=2.50cm,right=2.50cm,bottom=3.15cm]{geometry}
\usepackage{hyperref}

\newcommand{\R}{\mathbb R}

\def\P{{\mathbb P}}

\def\1{{\mathbbm 1}}

\newcommand{\dd}{\mathrm{d}}

\newtheorem{remark}{Remark}

\newcommand{\Min}{\displaystyle\min}
\newcommand{\Max}{\displaystyle\max}
\newcommand{\Sup}{\displaystyle\sup}

\newcommand{\xory}{{x\backslash y}}
\newcommand{\Int}{\displaystyle\int}
\newcommand{\Sum}{\displaystyle\sum}
\newcommand{\Lim}[2]{\displaystyle\lim\limits_{#1\,\rightarrow\, #2}}
\newcommand{\Frac}{\displaystyle\frac}
\newcommand{\Bigcup}{\displaystyle\bigcup}
\newcommand{\ov}{\overline}

\newcommand{\wt}{\widetilde}

\newcommand{\wh}[1]{\widehat{#1}}

\newcommand{\tra}{^{\text{t}}}
\newcommand{\thf}{^{\text{th}}}

\newcommand{\pds}{\centerdot}
\newcommand{\demi}{\frac{1}{2}}

\newcommand{\Dt}{\Delta t}

\newcommand{\gradx}[1]{\nabla_x{#1}}

\newcommand{\divx}{\nabla_x\pds}

\newcommand{\dt}[1]{\Frac{\partial {\,#1}}{\partial t}}

\newcommand{\ddt}[1]{\Frac{\mathrm{d} {\,#1}}{\mathrm{d} t}}

\newcommand{\vdt}{\partial_t}
\newcommand{\vdx}{\partial_x}

\newcommand{\vdy}{\partial_y}

\newcommand{\mc}[1]{\mathcal{#1}}

\newcommand{\inv}[1]{\frac{1}{#1}}
\newcommand{\Inv}[1]{\Frac{1}{#1}}

\newcommand{\bk}[1]{\llbracket{#1}\rrbracket}
\newcommand{\bs}[1]{\boldsymbol{#1}}

\newcommand{\wb}[1]{\widebreve{#1}}

\def\veps{\varepsilon}

\def\<{\left<}
\def\>{\right>}
\def\({\left(}
\def\){\right)}

\makeatletter
\def\widebreve{\mathpalette\wide@breve}
\def\wide@breve#1#2{\sbox\z@{$#1#2$}%
     \mathop{\vbox{\m@th\ialign{##\crcr
\kern0.08em\brevefill#1{0.8\wd\z@}\crcr\noalign{\nointerlineskip}%
                    $\hss#1#2\hss$\crcr}}}\limits}
\def\brevefill#1#2{$\m@th\sbox\tw@{$#1($}%
  \hss\resizebox{#2}{\wd\tw@}{\rotatebox[origin=c]{90}{\upshape(}}\hss$}
\makeatletter

\newcommand{\ext}{pdf}
\graphicspath{{FIGURES/\ext/}}

\newcommand{\ncell}{576}

\makeatletter
\def\input@path{{SOURCES/}}
\makeatother

\title{\textit{A Posteriori} local subcell correction of high-order discontinuous Galerkin scheme for conservation laws on two-dimensional unstructured grids}

\author{Fran\c{c}ois Vilar \thanks{IMAG, Univ Montpellier, CNRS, Montpellier, France (\email{francois.vilar@umontpellier.fr}).}
\and R\'emi Abgrall\thanks{Instit\"ut f\"ur Mathematik, Universit\"at Z\"urich, CH-8057 Z\"urich, Switzerland (\email{remi.abgrall@uzh.ch}).}}

\begin{document} 

\maketitle
\begin{abstract}
  In this paper, we present the two-dimensional unstructured grids extension of the \textit{a posteriori} local subcell correction (APLSC) of discontinuous Galerkin (DG) schemes introduced in \cite{vilar_aplsc_1D}. The technique is based on the reformulation of DG scheme as a finite volume (FV) like method through the definition of some specific numerical fluxes referred to as reconstructed fluxes. High-order DG numerical scheme combined with this new local subcell correction technique ensures positivity preservation of the solution, along with a low oscillatory and sharp shocks representation.\\

  The main idea of this correction procedure is to retain as much as possible the high accuracy and the very precise subcell resolution of DG schemes, while ensuring the robustness and stability of the numerical method, as well as preserving physical admissibility of the solution. Consequently, an \textit{a posteriori} correction will only be applied locally at the subcell scale where it is needed, but still ensuring the scheme conservation. Practically, at each time step, we compute a DG candidate solution and check if this solution is admissible (for instance positive, non-oscillating, \dots). If it is the case, we go further in time. Otherwise, we return to the previous time step and correct locally, at the subcell scale, the numerical solution. To this end, each cell is subdivided into subcells. Then, if the solution is locally detected as bad, we substitute the DG reconstructed flux on the subcell boundaries by a robust first-order numerical flux. And for subcell detected as admissible, we keep the high-order DG reconstructed flux which allows us to retain the very high accurate resolution and conservation of the DG scheme. As a consequence, only the solution inside troubled subcells and its first neighbors will have to be recomputed, elsewhere the solution remains unchanged. Another technique blending in a convex combination fashion DG reconstructed fluxes and first-order FV fluxes for admissible subcells in the vicinity of troubled areas will also be presented and prove to improve results in comparison to the original algorithm introduced in \cite{vilar_aplsc_1D}. Numerical results on various type problems and test cases will be presented to assess the very good performance of the design correction algorithm.
\end{abstract}

\begin{keyword}
  a posteriori correction, subcell correction,  arbitrary high-order, DG subcell FV formulation, positivity-preserving scheme,  hyperbolic conservation laws, subcell conservative scheme
\end{keyword}

%

\section{Introduction}
\label{intro}

This paper is devoted to the extension of the \textit{a posteriori} local subcell correction (APLSC) method introduced in \cite{vilar_aplsc_1D} to the two-dimensional unstructured case. It is well known that the Cauchy problem for any non linear hyperbolic problem leads to handling, in the generic situation, solutions with discontinuities. In numerical simulations, this aspect has to be dealt with, and this issue has been one of the main questions to address for designing numerical methods for hyperbolic problems. Another core problem is the one of accuracy. The last one is that the problem is well set only if the solution belongs to an admissibility set called invariant domain. For example, for the gas dynamics Euler equations, the density and the internal energy must stay positive. This is particularly difficult to address, simultaneously, these three questions together, because those constraints are antagonistic. These questions has been at the center of algorithmic developments since decades, one may mention \cite{FCT,vanL,Harten,jiang_weno,shulimit} and the reference herein.\\

The Discontinuous Galerkin (DG) method is one of the most widely used numerical scheme, especially in the context of computational fluid dynamics.  Cockburn, Shu \textit{et al} in a series  papers \cite{Cockburn_lcs2,Cockburn_lcs3,Cockburn_lcs4,Cockburn_lcs5} have paved the road to efficient methods for fluid dynamics.  DG methods allow to reach any arbitrary order of accuracy, while keeping the stencil compact, along with other good properties such as built-in entropy stability  and $hp$-adaptivity. However, though recent progress have been done \cite{shu}, designing method that are oscillation free and compliant with the invariant domain is not a trivial matter. Controlling spurious oscillations has been studied in many papers, among which \cite{bisw,burb01,kriv07,MYang,MLP_unstruct,Kuzmin13,Li_vertex_based_lim}. Staying in the invariant domain has also  been made possible, see for example \cite{shulimit,siamreview, guermond}. However, this is often achieved to the cost of enlarging the width of the discontinuous patterns, and some time a lost of accuracy.\\
 
Other methods have recently gained in popularity, the so-called subcell techniques. Here the idea is to subdivide the bad cells, and to adopt a special procedure with the hope of curing the negative aspects of the original scheme. Examples of this strategy can be found in \cite{peraire_2012,peraire_2013,munz_subFV_2014,munz_subFV_2018}. For example,  \cite{peraire_2012}, the authors use a convex combination between high-order DG schemes and first-order Finite-Volumes (FV) on a subgrid, allowing them to retain the very high accurate resolution of DG in smooth areas and ensuring the scheme robustness in the presence of shocks. Similarly, in \cite{munz_subFV_2014,munz_subFV_2018}, after having detected the troubled zones, cells are then subdivided into subcells, and a robust first-order finite volume scheme is performed on the subgrid in troubled cells. Alternatively, some robust high-order scheme as MUSCL or WENO could either be used to avoid too much accuracy discrepancy. Note that in general this methods are tuned for Cartesian meshes.\\
 
Another approach which is worth to be mentioned is the so-called MOOD technique, \cite{mood1,mood2,mood3}. Through this procedure, the order of approximation of the numerical scheme is locally reduced in an \textit{a posteriori} sequence until the solution becomes admissible, $i.e$ oscillation free and the solution lives in the invariant domain. In \cite{dumbser_subFV_lim,dumbser_subFV_lim_tri}, a subcell FV technique similar to the one presented in \cite{munz_subFV_2014} has been applied to the \textit{a posteriori} paradigm. In practice, if the numerical solution in a cell is detected as bad, the cell is then subdivided into subcells and a first-order finite volume, or alternatively other robust scheme (second-order TVD FV scheme, WENO scheme, \dots), is applied on each subcell. Then, through these new subcell mean values, a high-order polynomial is reconstructed on the primal cell. This correction procedure has the benefit to be very simple and robust, and is able to preserve the high accuracy of DG schemes in smooth areas.\\
 
In all these aforementioned limitation techniques, \textit{a priori} and \textit{a posteriori}, the high-order DG polynomial is either globally modified in the cell, or even discard as it is in the (H)WENO limiter or any \textit{a posteriori} correction technique in the troubled cells. Since one of the main advantage of high-order scheme is to be able to use coarse grids while still being very precise, one can see that there is a waste of information here, as well as unnecessary computational effort made. This is particularly the case in the vicinity of discontinuities since the polynomials are globally modified. This problem was addressed in the one-dimensional case in \cite{vilar_aplsc_1D}. This new technique relies on the reformulation of DG schemes as a FV-like scheme defined on a subgrid. First, as the number of subcells matches the dimension of the polynomial solution space, the numerical solution inside a cell can be uniquely defined by either a high degree polynomial, or as a piecewise constant solution through its different subcell mean values; the connection between the two being done via a projector. Second, by means of particular basis functions introduced in \cite{vilar_aplsc_1D}, which are nothing but the $L_2$ projection over the polynomial space of the subcells indicator function, the DG volume and boundary terms can be rewritten as flux differences. Such theoretical analysis is relatively simple in one dimension in space, but much more challenging in two dimensions.\\

The format of the paper is as follow. Extending the theoretical analysis introduced in \cite{vilar_aplsc_1D} and using ideas from \cite{remark_RD_abgrall}, we first reinterpret  unstructured grid DG scheme as a subgrid FV-like scheme, through the definition of particular fluxes that we referred to as reconstructed fluxes. These DG reconstructed fluxes are analytically computed, and the analysis shows how those fluxes are connected to the interior polynomial flux and the jump at the cell interface between the interior flux and the DG numerical flux. Let us emphasize that reformulation of DG scheme as a subcell finite-volume method can be performed regardless the form of the element: this is thus not limited to quadrilateral nor triangular cells, but can be done on general polygonal elements. This theoretical part is done in section \ref{sect_DG_as_FV}, where a discussion on the type of subcells is also provided. Using this equivalent formulation, we can proceed by means of a MOOD paradigm as follow: at each time step, we compute a DG candidate solution and check if this solution is admissible. If it is the case, we go further in time. Otherwise, we return to the previous time step and correct locally, at the subcell scale, the numerical solution. In the subcells where the solution was detected as bad, we substitute the DG reconstructed flux on the subcell boundaries by a robust first-order numerical flux. And for subcells detected as admissible, we keep the high-order reconstructed flux which allows us to retain the very high accurate resolution and conservation of DG schemes. Consequently, only the solution inside troubled subcells and their first neighbors will have to be recomputed. Elsewhere, the solution remains unchanged. This correction procedure is then extremely local. Another technique blending in a convex combination fashion DG reconstructed fluxes and first-order FV fluxes for admissible subcells in the vicinity of troubled areas will also be presented and prove to improve results in comparison to the original algorithm. This procedure and all the technical details are  described  in section  \ref{sect_DG_lim}. Section \ref{sect_results} provides numerical results, and a conclusion follows.

\newpage
\section{Discontinuous Galerkin method as a subcell Finite Volume scheme}
\label{sect_DG_as_FV}

This section is devoted to the demonstration of the equivalency between DG schemes and a FV-like method on a subgrid provided the definition of particular fluxes. To remain as simple as possible, two-dimensional Scalar Conservation Laws (SCL) will be considered. Let $u=u(\bs{x},t)$, for $\bs{x} \in \omega\subset \mathbb{R}^2$ and $t \in [0,T]$, be the solution of the following system,

\begin{subequations}
\label{lcs_eq2D}
\begin{empheq}[left = \empheqlbrace\,]{align}
&\vdt{\,u}(\bs{x},t)+\divx{\bs{F}\(u(\bs{x},t)\)}=0, && (\bs{x},t)\in\,\omega\times[0,T], \label{lcs1}\\[3mm]
&u(\bs{x},0)=u_0(\bs{x}), &&\bs{x}\in\,\omega, \label{lcs2}
\end{empheq}
\end{subequations}\vspace*{5mm}

where $u_0$ is the initial data and $\bs{F}(u)$ the 2D flux function. For the subsequent discretization, let us introduce the following notation. Let $\{\omega_c\}_c$ be a generic partition of the domain $\omega$ into non-overlapping cells, with $|\omega_c|$ being the size of $\omega_c$. We also partition the time domain in intermediate times $(t^n)_n$ with $\Delta t^n=t^{n+1}-t^n$ the $n^{th}$ time step. In the DG frame, the numerical solution is considered piecewise polynomial over the domain, and hence developed on each cell onto $\P^{\,k} (\omega_c)$, the set of polynomials of degree up to $k$ defined on cell $\omega_c$. This space approximation theoretically leads to a $(k+1)^{th}$ space order accurate scheme. Let $u_h^c$ be the restriction of $u_h$, the piecewise polynomial approximation of the solution $u$, over the cell $\omega_c$
\begin{align}
  u_h^c(\bs{x},t)=\Sum_{m=1}^{N_k}u_m^c(t)\,\sigma_m^c(\bs{x}),
\label{u_h}
\end{align}
where the $u_m^c$ are the $N_k$ successive components of $u_h^c$ over the polynomial basis $\{\sigma_m^c\}_m$. We recall that in the two-dimensional case  $N_k=\frac{(k+1) (k+2)} {2}$. The coefficients $u_m^c$ present in \eqref{u_h} are the solution moments to be computed through a local variational formulation on $\omega_c$. To this end, one has to multiply equation \eqref{lcs1} by $\psi \in\,\P^{\,k}(\omega_c)$, a polynomial test function, and integrate it on $\omega_c$. By means of an integration by parts and substituting the solution $u$ by its approximated polynomial counterpart $u_h^c$, one gets\\[-2mm]
\begin{align}
  \label{DG_2D_1}
  \Int_{\omega_c} \dt{u_h^c}\,\psi\,\dd V=\Int_{\omega_c} \bs{F}(u_h^c)\pds\gradx{\psi}\,\dd V-\Int_{\partial \omega_c}\hspace*{-2mm}\psi\;\mc{F}_n\,\dd S, \qquad \forall \,\psi\in\P^k(\omega_c).
\end{align}
The numerical solution $u_h^c$ is then the unique polynomial function defined in $\,\P^{\,k}(\omega_c)$ satisfying equation \eqref{DG_2D_1} for all function $\psi\in\P^k(\omega_c)$. In \eqref{DG_2D_1}, the numerical flux function $\mc{F}_n$, in addition to ensure the scheme conservation, is the cornerstone of any finite volume or DG scheme regarding fundamental considerations as stability, positivity and entropy among others. In this context, this numerical flux is defined as a function of the two states on the left and right of each interface\\[-3mm]
\begin{align}
  \label{num_flux_0}
  \mc{F}_n=\mc{F}\(u^-,u^+,\bs{n}\),
\end{align}
with $u^-=\Lim{\epsilon}{0^+}u_h^c(\bs{x}_i-\epsilon\,\bs{n},\,t)$ and $u^+=\Lim{\epsilon}{0^+}u_h^v(\bs{x}_i+\epsilon\,\bs{n},\,t)$, where $\omega_v$ is a face neighboring cell of $\omega_c$, while $\bs{x}_i$ and $\bs{n}$ respectively stand for a point and the unit outward normal of the separating interface. From now on, we refer $\mc{V}_c$ to as the set containing the face neighboring cells of $\omega_c$. This function is generally obtained through the resolution of an exact or approximated Riemann problem. In the remainder of this paper, for sake of simplicity, we make use of the very well-known global Lax-Friedrichs numerical flux which reads
\begin{align}
  \label{num_flux}
  \mc{F}(u,v,\bs{n})=\Frac{\big(\bs{F}(u)+\bs{F}(v)\big)}{2}\pds\bs{n}-\Frac{\gamma}{2}\,(v-u),
\end{align}
where $\gamma=\Sup_w\(\|d_w \, \bs{F}(w)\|_2\)$.\\

Now, taking in \eqref{DG_2D_1} the test function $\psi$ among the polynomial basis functions leads to the following linear system allowing the calculation of the solution moments $u_m^c$
\begin{align}
  \label{DG_2D_2}
  \Sum_{m=1}^{N_k}\ddt{u_m^c}\,\Int_{\omega_c} 
  \sigma^c_m\,\sigma^c_p\,\dd V=\Int_{\omega_c} \bs{F}(u_h^c)\pds\gradx{\sigma^c_p}\,\dd V-\Int_{\partial \omega_c}
  \sigma^c_p\;\mc{F}_n\,\dd S, \qquad \forall \,p\in\,[1, N_k].
\end{align}
The terms $\int_{\omega_c} \bs{F}(u_h^c)\pds\gradx{\sigma^c_p}\,\dd V$ and $\int_{\partial \omega_c}\sigma^c_p\;\mc{F}_n\,\dd S$ are respectively referred to as volume and surface integrals. In \eqref{DG_2D_2}, we identify $\int_{\omega_c} \sigma^c_m\,\sigma^c_p\,\dd V=\(M_c\)_{mp}$ as the generic coefficient of the symmetric mass matrix $M_c$. The scheme \eqref{DG_2D_2} can then be reformulated in a compact matrix-vector form as follows
\begin{align}
  \label{DG_2D_3}
  M_c\,\ddt{U_c}=\Phi_c,
\end{align}\vspace*{1mm}
with $(U_c)_m=u_m^c$ the solution vector filled with the polynomial moments, and where the so-called DG residuals $\Phi_c$ write
\begin{align}
  \label{DG_residual}
(\Phi_c)_m=\Int_{\omega_c} \bs{F}(u_h^c)\pds\gradx{\sigma^c_m}\,\dd V-\Int_{\partial \omega_c}\hspace*{-2mm}\sigma^c_m\;\mc{F}_n\,\dd S.
\end{align}

Similarly to what has been done in the one-dimensional case in \cite{vilar_aplsc_1D}, let us now demonstrate the equivalency between discontinuous Galerkin schemes and a finite volume like method on a subgrid, and exhibit the corresponding subcell numerical fluxes that will be referred to as high-order reconstructed fluxes. To do so, we first need to subdivide the mesh cells into subcells. Let us emphasize that to obtain a relation of equivalency, one would need the same number of Degrees of Freedom (DoF) as number of subcells. Indeed, it will then be possible to get a one-to-one relation between subcell mean values obtained through a FV-like method provided with reconstructed fluxes and the solution polynomial moments deriving from the DG schemes previously introduced. But the procedure presented in the following is not limited to this case, and this is absolutely possible to exhibit the corresponding reconstructed fluxes even in the case where the number of subcells does not fit the dimension of the polynomial space $\P^{\,k}$. In the case where more subcells than DoF are used, one can either apply a least square procedure to the set of submean values obtained through the FV-like method to return to a $(k+1)\thf$ order DG scheme, or directly return to an higher-order DG scheme in the particular case where the number of subcells corresponds to the dimension of higher degree polynomial space. If one uses less subcells than DoF, then with the same techniques than in the previous case, one can return to a lower DG scheme from the subcell mean values set. In the remainder, for sake of positivity-preserving purpose, we will consider $N_k$ subcells.\\

\begin{figure}[!ht]
  \begin{center}
    \subfigure[]{\includegraphics[height=4.5cm]{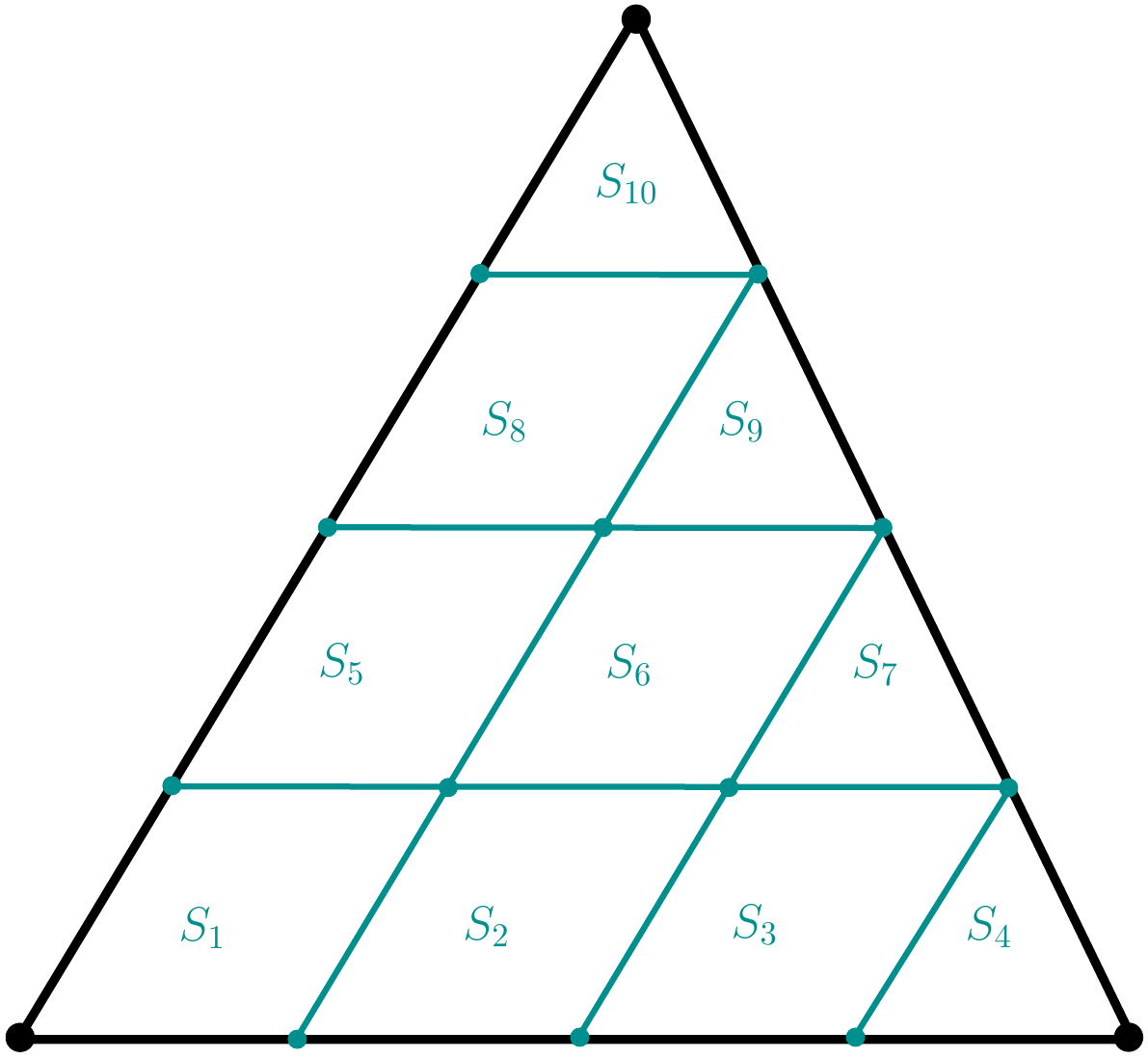}\label{fig_tri1}}\hspace*{4mm}
    \subfigure[]{\includegraphics[height=4.5cm]{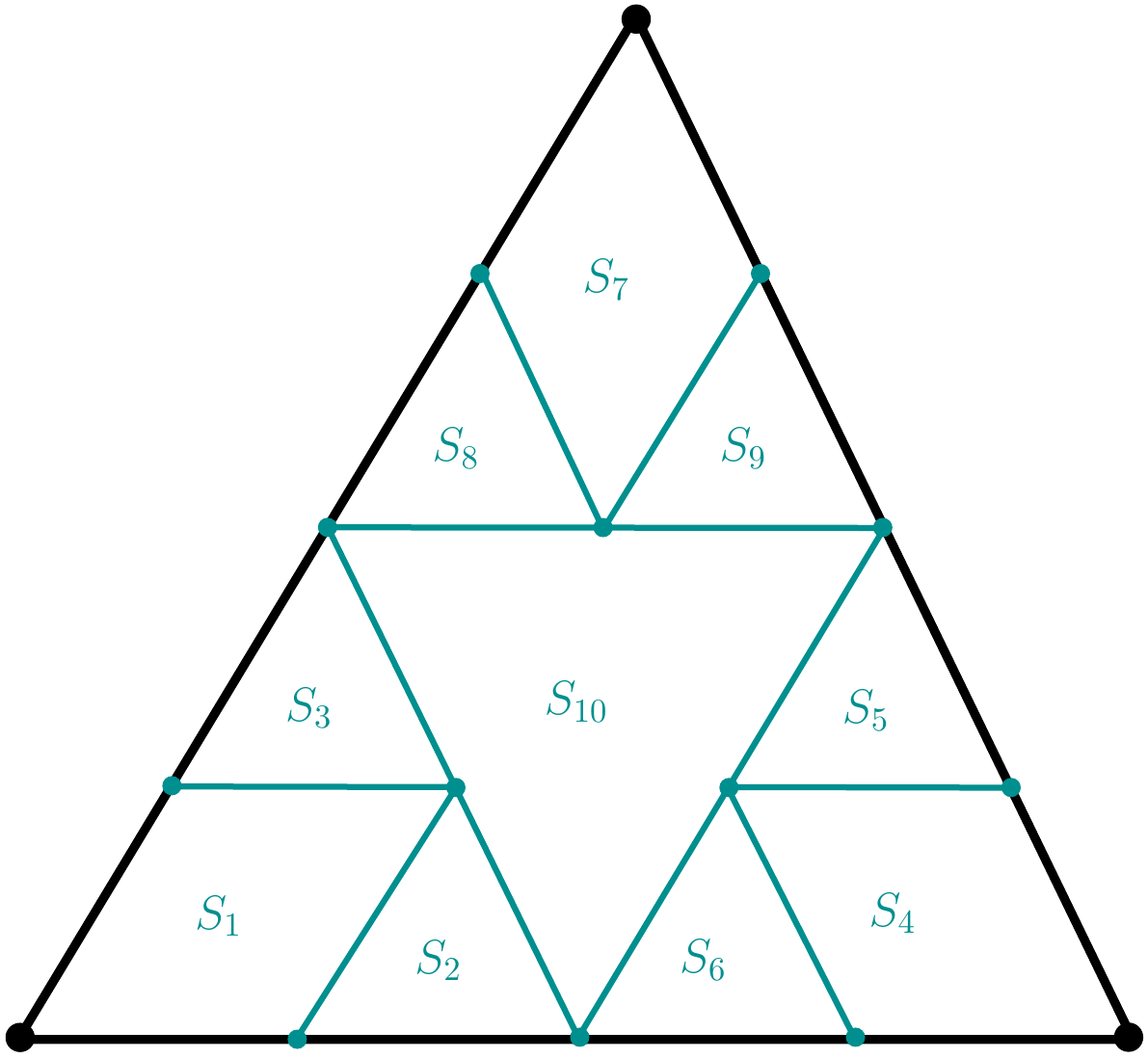}\label{fig_tri2}}\hspace*{4mm}
    \subfigure[]{\includegraphics[height=4.5cm]{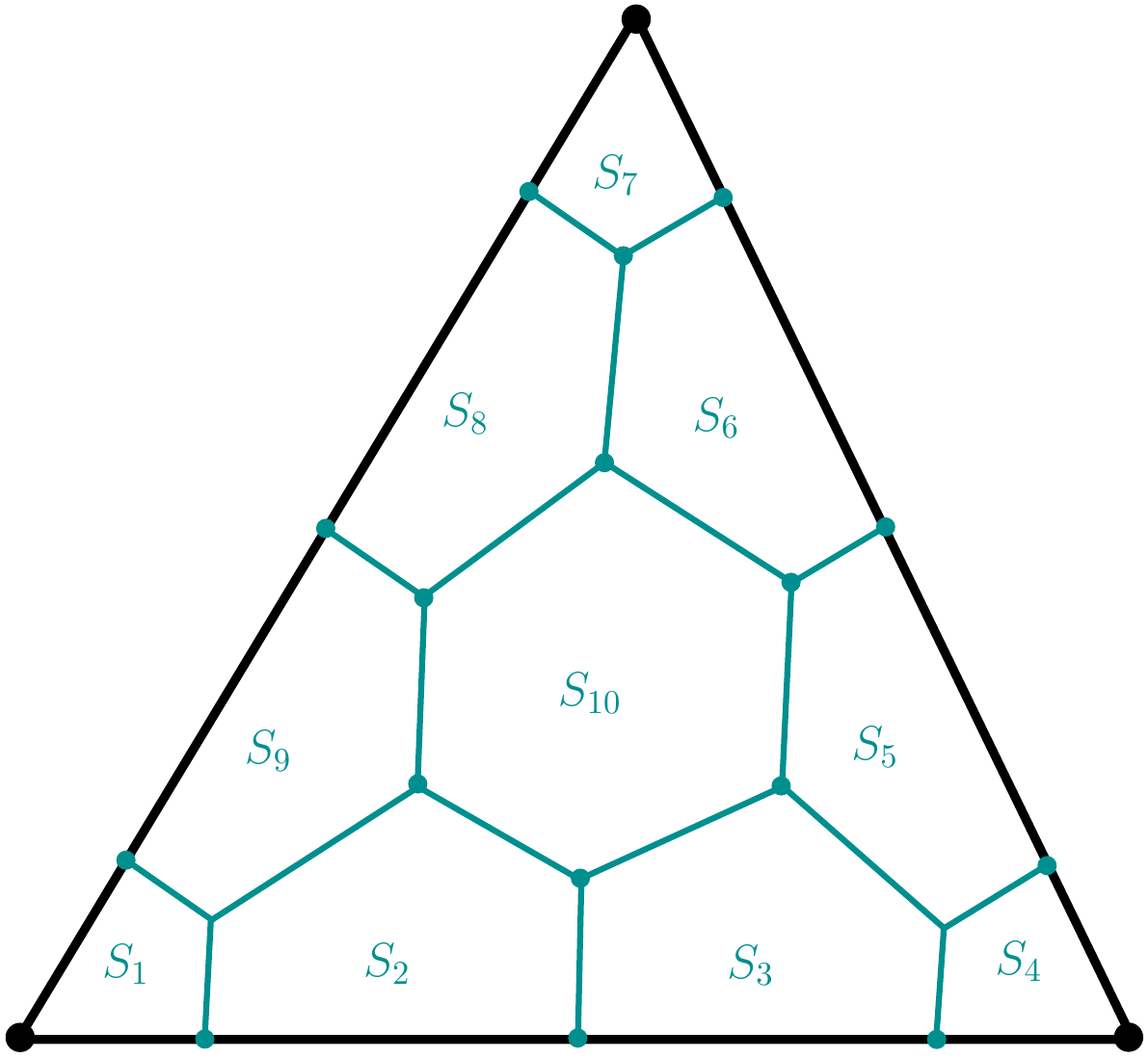}\label{fig_tri3}}
    \caption{Examples of subdivision for a triangular cell and a $\P^3$ DG scheme.}
  \label{fig_subdiv_tri}
  \end{center}
\end{figure}

Even if the choice of the cells subdivision may have an effect of the DG correction technique, for the following theoretical part it has no influence what so ever. The only constrain is that the projection matrix $P_c$ further introduced in \eqref{proj_mat} has to been non-singular. Many subdivisions can be found in the literature for other methods relying of a subgrid, as spectral volume methods for instance \cite{Wang_SV_2002_II,Wang_SV_2009} or subcell shock capturing technique as \cite{peraire_2012}. In Figure~\ref{fig_subdiv_tri}, three types of subdivision are displayed for a triangular cell, in the case of a $4\thf$-order DG scheme. Let us emphasize that the one depicted in Figure~\ref{fig_tri3} is the most widely used in subgrid techniques, and has the advantages to be invariant by rotation and can be generated for any order of accuracy. For the numerical applications, this subdivision will be compared to a way simpler one, Figure~\ref{fig_tri1}, which has the benefits to be extremely simple to implement and where the subcells normals are nothing but the primal triangular cell ones.\\

\begin{figure}[!ht]
  \begin{center}
    \subfigure[$\P^1$.]{\includegraphics[height=4.5cm]{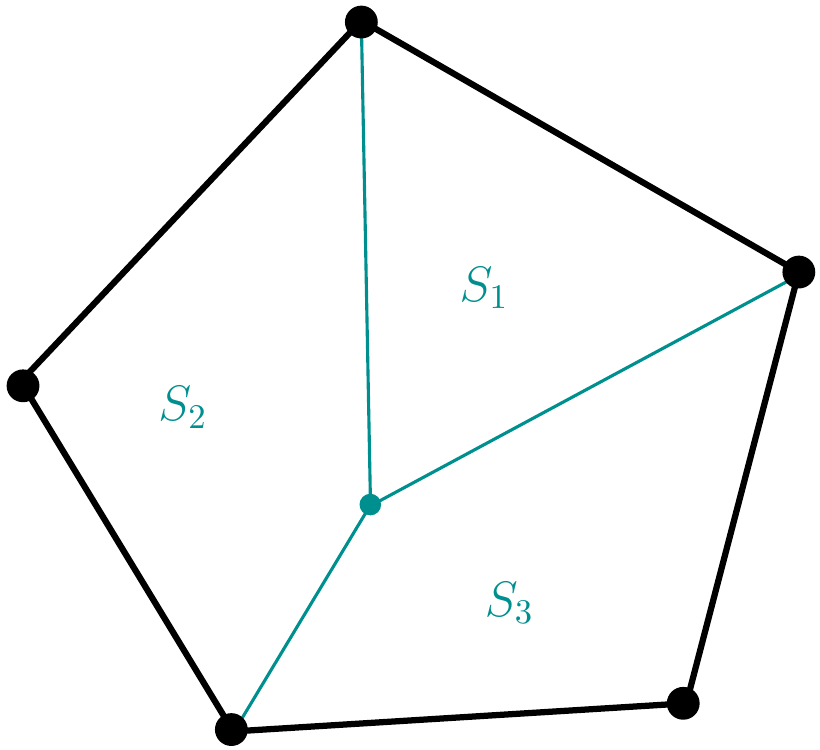}\label{fig_poly1}}\hspace*{4mm}
    \subfigure[$\P^2$.]{\includegraphics[height=4.5cm]{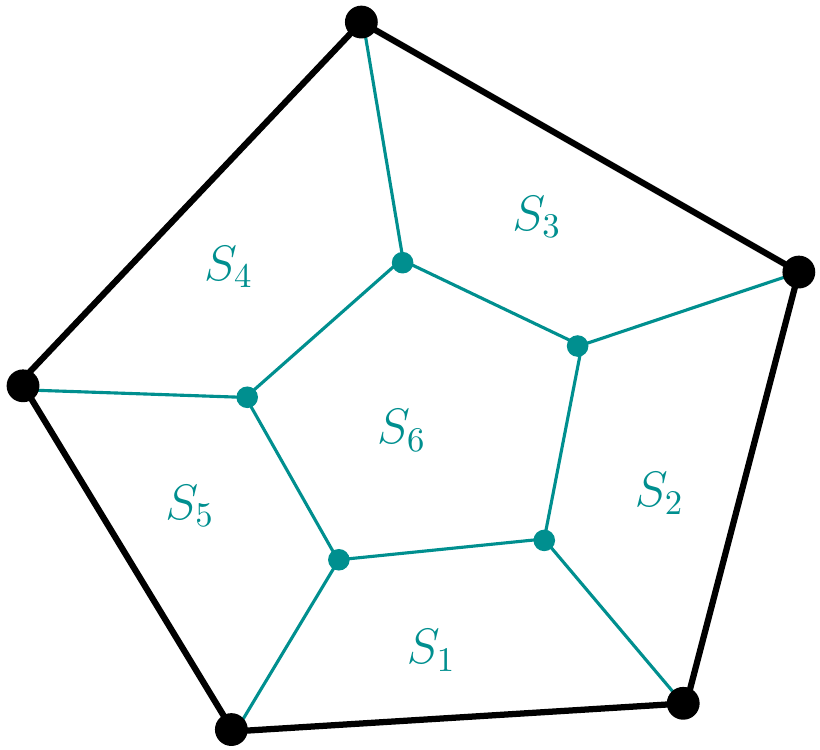}\label{fig_poly2}}\hspace*{4mm}
    \subfigure[$\P^3$.]{\includegraphics[height=4.5cm]{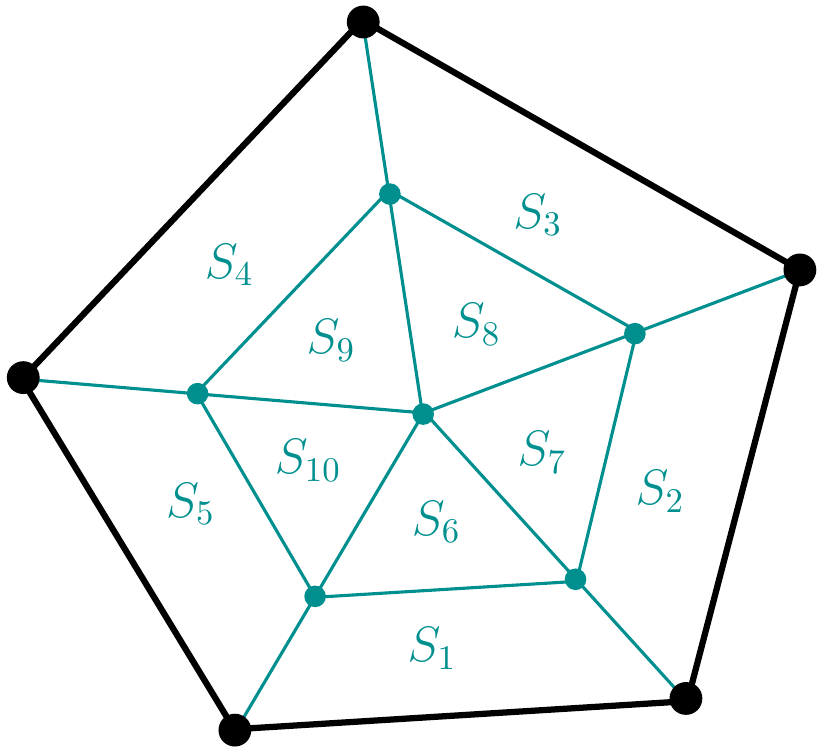}\label{fig_poly3}}
    \caption{Examples of subdivision for a polygonal cell and DG schemes from $\P^1$ to $\P^3$.}
  \label{fig_subdiv_polyg}
  \end{center}
\end{figure}

Even if only triangular grids are considered for numerical applications, let us emphasize that the following demonstration as well as the subcell correction technique presented in this paper are not limited to this case. Any grid made of generic polygonal cells can be considered. Examples of possible subdivisions for DG schemes up to $\P^3$ are displayed in Figure~\ref{fig_subdiv_polyg}. Curvilinear meshes with curvilinear cell subdivision could also be used, and will be the topic of a near future paper.\\

That being said, let us consider a cell $\omega_c$ and its subdivision into $N_k$ subcells $S_m^c$, for $m\in[1,N_k]$. For any function $\psi \in L^2(\omega_c)$, we define the corresponding subcell mean values, also referred to as submean values

\begin{align}
  \label{submean}
  \ov{\psi}_m^{\,c}=\Inv{|S_m^c|}\Int_{S_m^c}\psi\,\dd V.
\end{align}\vspace*{2mm}

Let us now reformulated DG as a FV-like scheme provided the definition of the so-called reconstructed fluxes.

\subsection{Reconstructed flux through residuals}
\label{FR_residual}

Applying definition \eqref{submean} to $\vdt{u_h^c}$ and by means of \eqref{u_h}, one gets
\begin{align}
  \label{submean_rel1}
  \ddt{\ov{u}_m^c}=\Inv{|S_m^c|}\,\Sum_{q=1}^{N_k}\ddt{u_q^c}\,\Int_{S_m^c}\sigma^c_q\,\dd V,
\end{align}
which can be put into into a matrix-vector form as follows
\begin{align}
  \label{submean_rel2}
  \ddt{\ov{U}_c}=P_c\,\ddt{U_c}.
\end{align}
In \eqref{submean_rel2}, the projection matrix $P_c$, which has to be invertible for the subdivision to be admissible, is defined
\begin{align}
  \label{proj_mat}
  (P_c)_{mp}=\Inv{|S_m^c|}\,\Int_{S_m^c} \sigma^c_p\,\dd V.
\end{align}

\begin{remark}
  \label{reconstruction_matrix}
  Let us note that if the number of subcells does not coincide with the polynomial space dimension $N_k$, then the projection matrix $P_c$ is not square. For the subdivision under consideration to be admissible, this projection matrix has to present a left inverse, the reconstruction operator $R_c$, such that $R_c\,P_c=I_{N_k}$. In the case of $N_k$ subcells, $R_c$ reduces to $P_c^{-1}$.
\end{remark}\vspace*{3mm}

In relation \eqref{submean_rel2}, $\ov{U}_c$ is nothing but the vector containing the cell submean values, $i.e.\, (\ov{U}_c)_m=\ov{u}_m^c$. Now, by means of DG scheme definition \eqref{DG_2D_3}, it follows that
\begin{align}
  \label{submean_rel3}
  \ddt{\ov{U}_c}=P_c\,M_c^{-1}\,\Phi_c.
\end{align}
To express \eqref{submean_rel3} as a FV-like scheme, we now introduce the DG reconstructed flux $\wh{F_n}$ such that
\begin{align}
  \label{recons_flux}
  \ddt{\ov{u}_m^c}=-\Inv{|S_m^c|}\,\Int_{\partial S_m^c}\wh{F_n}\;\dd S.
\end{align}
By introducing $\mc{V}_m^{\,c}$, the set of face neighboring subcells of $S_m^c$, this last expression rewrites
\begin{align}
  \label{recons_flux1}
  \ddt{\ov{u}_m^c}=-\Inv{|S_m^c|}\,\Sum_{S_{p}^v\,\in \,\mc{V}_m^{\,c}}\Int_{f_{mp}^c} \wh{F_n}\;\dd S,
\end{align}
where $f_{mp}^c$ denotes the face between subcells $S_m^c$ and $S_p^v$. Let us emphasize that $S_{p}^v\,\in \,\mc{V}_m^{\,c}$ can either be inside cell $\omega_c$ or in one of its neighbors $\omega_v$ with $v\neq c$. This situation is displayed in Figure~\ref{fig_flux_lim} where $S_m^c$ would be colored red, while its faces neighboring subcells would be colored green.\\

Now, similarly to what has been done in the 1D case, we impose that on the boundary of cell $\omega_c$ the reconstructed flux coincide with the DG numerical flux
\begin{align}
  \label{recons_flux_boundary}
  \wh{F_n}_{|_{\partial\omega_c}}=\mc{F}_n,
\end{align}

where $\mc{F}_n$ has been previously introduced in \eqref{num_flux_0}. This particular choice, while not unique, permits to recover in the most natural way the time evolution governing equation of the cell average value $\ov{u}_c$
\begin{align*}
  \ddt{\ov{u}_c} := \Inv{|\omega_c|}\,\Int_{\omega_c} \dt{u_h^c}\,\dd V = -\Inv{|\omega_c|}\, \Int_{\partial \omega_c}\hspace*{-2mm}\mc{F}_n\,\dd S,
\end{align*}

obtained through DG scheme definition \eqref{DG_2D_1} with $\psi=1$. By means of \eqref{recons_flux_boundary}, definition \eqref{recons_flux1} then rewrites
\begin{align}
  \label{recons_flux2}
  \ddt{\ov{u}_m^c}=-\Inv{|S_m^c|}\,\(\Sum_{S_{p}^c\,\in \,\wb{\mc{V}_m^{\,c}}}\Int_{f_{mp}^c} \wh{F_n}\;\dd S\;+\Int_{\partial S_m^c\cap\partial \omega_c}\hspace*{-0.75cm}\mc{F}_n\;\dd S\),
\end{align}
where  $\wb{\mc{V}_m^{\,c}}$ stands for the set containing only the face neighboring subcells of $S_m^c$ inside $\omega_c$. For now on, an orientation will be assigned to each face. Then, taking two subcells $S_m^c$ and $S_p^v$, we introduce the sign function $\veps^c_{mp}$  depending on the orientation of face $f_{mp}^c$
\begin{align}
  \label{sign_fct}
  \veps_{mp}^c=\left\{\begin{array}{ll} 
           1\quad &\text{if face }\; f_{mp}^c \; \text{is direct or if }f_{mp}^c\subset\partial\omega_c,\\[3mm] 
           -1\quad &\text{if face }\; f_{mp}^c \; \text{is indirect},\\[3mm]
           0\quad &\text{if }\; S_p^v\notin \mc{V}_m^c.
           \end{array}\right.
\end{align}

Obviously, $\forall S_p^c\in\wb{\mc{V}_m^c}$, we have $\veps_{pm}^c=-\veps_{mp}^c$. Now, unlike the 1D case where a pointwise definition of the reconstructed flux was given inside the cell, we make use here of a face integrated version of the high-order DG reconstructed flux. Indeed, for a face $f_{mp}^c$, let $\wh{F_{mp}}$ be defined as follows
\begin{align}
  \label{recons_flux_int}
  \Int_{f_{mp}^c} \wh{F_n}\;\dd S=\veps^c_{mp}\,\wh{F_{mp}}.
\end{align}
Since the face orientation has been carried through $\veps^c_{mp}$, the face integrated quantity $\wh{F_{mp}}$ is then continuous, $i.e. \;\forall S_p^c\in\wb{\mc{V}_m^c},\quad \wh{F_{pm}}=\wh{F_{mp}}$. Now, denoting by $N_f^c$ the number of subcells faces inside $\omega_c$, meaning not belonging to $\partial \omega_c$, let us introduce $\wh{F_c}\in\R^{N_f^c}$ the vector containing all the interior faces reconstructed fluxes. The subcell mean values governing equations \eqref{recons_flux2} then yield the following system
\begin{align}
  \label{recons_flux3}
  -A_c\,\wh{F_c}=D_c\,\ddt{\ov{U}_c}+B_c,
\end{align}
where $A_c\in\mc{M}_{N_k\times N_f^c}$, defined as $(A_c)_{mp}=\veps^c_{mp}$, stands for the adjacency matrix, and where $D_c=\text{diag}\(|S_1^c|,\dots,|S_{N_k}^c|\)$ is the subcells volume matrix and $(B_c)_m=\Int_{\partial S_m^c\cap\partial \omega_c}\hspace*{-0.75cm}\mc{F}_n\;\dd S$ carries the cell boundary contribution. Finally, introducing the matrix $Q_c=D_c\,P_c$ and substituting equation \eqref{submean_rel3} into \eqref{recons_flux3}, one gets
\begin{align}
  \label{recons_flux4}
  -A_c\,\wh{F_c}=Q_c\,M_c^{-1}\,\Phi_c+B_c.
\end{align}
To solve such system, and obtain an explicit expression of the reconstructed fluxes $\wh{F_c}$, we make use of the technique employed in \cite{remark_RD_abgrall} in a similar context. Since the graph is connected, the rank and the kernel of the Laplacian matrix of the interior subgrid graph $L_c=A_c\,A_c\tra$ respectively are $N_k-1$ and $\text{span}\{\bs{1}\}$. Then the inverse of $L_c$ on the orthogonal of its kernel is defined, and is henceforth denoted by $\mc{L}_c^{-1}$. For any $\lambda \neq 0$, this inverse writes
\begin{align}
  \label{graph_laplacian}
  \mc{L}_c^{-1}=(L_c+\lambda\,\Pi)^{-1}-\Inv{\lambda}\,\Pi,
\end{align}
with $\Pi=\inv{N_k}\,\(\bs{1}\otimes\bs{1}\)\in\,\mc{M}_{N_k}$. Finally, we are able to exhibit the following definition of the reconstructed flux
\begin{align}
  \label{recons_flux_definition}
  \wh{F_c}=-A_c\tra\,\mc{L}_c^{-1}\(Q_c\,M_c^{-1}\,\Phi_c+B_c\).
\end{align}
\begin{remark}
\label{remark_recons_flux_res}
In the definition of the reconstructed flux through DG residual, \eqref{recons_flux_definition}, the only terms depending on time are $\Phi_c$, the residual which is directly available in any DG code, and the boundary contribution $B_c$. All the other terms can be evaluated initially, once and for all.\\
Furthermore, let us emphasize that if all the mesh cells have the same structure, as triangles for example, then by means of a mapping to a referential triangle, one can prove that the projection matrix $P_c$, the adjacency matrix $A_c$ and the generalized inverse of the graph Laplacian matrix $\mc{L}_c^{-1}$ do not depend on the cell under consideration, but only on the order of approximation and on the choice of the subdivision.
\end{remark}\vspace*{2mm}

\begin{remark}
\label{remark_proj_mat}
As said previously, for the subdivision to be admissible, the projection matrix has to be invertible, which is almost always true. However, the type of cell subdivision may have an influence on the condition number of the projection matrix, and then may introduce more numerical errors during the simulation. In Section~\ref{sect_results}, the projection matrix condition numbers are gathered for the four types of cell subdivisions used in practical applications, for different orders of accuracy, see Table~\ref{table_cond_number}.  
\end{remark}\vspace*{2mm}

Once we have computed the reconstructed flux $\wh{F_c}$, we can simply recover the polynomial solution governing equation as follows
\begin{align}
  \label{recons_flux_polynomial}
  \ddt{U_c}=-Q_c^{-1}\(A_c\,\wh{F_c}+B_c\).
\end{align}
Let us gather all the previous results into the following theorem.
\begin{theorem}
  \label{main_thm_residual}
  Discontinuous Galerkin schemes, expressed in the following equation in cell $\omega_c$ in a compact matrix-vector form through residuals
  \begin{align}
    \label{DG_2D_thm_res}
    \ddt{U_c}=M_c^{-1}\,\Phi_c,
  \end{align}
  where the DG residual $\Phi_c$ has been defined in \eqref{DG_residual}, can be recast into $N_k$ subcell finite volume like schemes as
  \begin{align}
    \label{sub_FV_thm_res}
    \ddt{\ov{U}_c}=-D_c^{-1}\(A_c\,\wh{F_c}+B_c\),
  \end{align}
  where the finite volume fluxes $\wh{F_c}$, referred to as reconstructed fluxes, is defined through equation \eqref{recons_flux_definition}. All the other quantities involved in \eqref{DG_2D_thm_res} and \eqref{sub_FV_thm_res} have been previously introduced.
\end{theorem}

Now, for a deeper understanding of the reconstructed flux $\wh{F_c}$, let us seek its relation with the interior flux $\bs{F}(u_h^c)$ and the DG numerical flux $\mc{F}_n$.
\subsection{Reconstructed flux through fluxes}
\label{FR_flux}
In this section, we aim at expressing the reconstructing fluxes only through the interior flux $\bs{F}(u_h^c)$ and a correction term taking into account the jump at the cell boundaries, similarly to what is done in SBP operator with SAT boundary treatment \cite{SBP_carpenter,DG_SBP_gassner} or CPR schemes \cite{Huynh_FR,Wang_lifting}. To do so, we will not make use of the definition of DG scheme through residual \eqref{DG_2D_thm_res} but through fluxes \eqref{DG_2D_1}. The first step is to substitute in \eqref{DG_2D_1} the interior flux $\bs{F}(u_h^c)$ by $\bs{F}_h^c$ its $L_2$ projection onto $\(\P^{k+1}(\omega_c)\)^2$. If one uses nodal DG or any collocation of the interior flux, this step can obviously be skipped. Performing a second integration by parts leads to the so-called strong form of DG scheme
\begin{align}
  \label{strong_DG}
  \Int_{\omega_c} \dt{u_h^c}\,\psi\,\dd V=-\hspace*{-0.5mm}\Int_{\omega_c}\hspace*{-1.2mm} \psi\,\divx{\bs{F}_h^c}\,\dd V+\Int_{\partial \omega_c}\hspace*{-2.7mm}\psi\(\bs{F}_h^c\pds\bs{n}-\mc{F}_n\)\,\dd S,\qquad \forall \,\psi\in\P^k(\omega_c).
\end{align}
Similarly to what has been done in \cite{vilar_aplsc_1D} for the one-dimensional case, let us introduce the $N_k$ sub-resolution basis functions $\{\phi_m\}_m$. These particular basis functions of $\P^k(\omega_c)$, which can be seen as the $L_2$ projection of the subcell indicator functions $\1_{S_m^c}(\bs{x})$ onto $\P^k(\omega_c)$, are defined such that $\forall\,\psi\in\P^k(\omega_c)$
\begin{align}
  \label{subresolution}
  \Int_{\omega_c}\phi_m\,\psi\,\dd V=\Int_{S_m^c}\psi\,\dd V,\qquad \forall\,m=1,\ldots,N_k.
\end{align}
Now, because equation \eqref{strong_DG} holds for any polynomial function $\psi$ of degree $k$, let us substitute $\phi_m$ for $\psi$ in DG schemes
\begin{align}
  \label{recons_flux5}
  \Int_{\omega_c} \dt{u_h^c}\;\phi_m\;\dd V=-\Int_{\omega_c} \phi_m\;\divx{\bs{F}_h^c}\;\dd V+\Int_{\partial \omega_c}\hspace*{-2mm}\phi_m\;\(\bs{F}_h^c\pds\bs{n}-\mc{F}_n\)\,\dd S.
\end{align}
Through the sub-resolution property \eqref{subresolution}, equation \eqref{recons_flux5} can be recast as
\begin{align}
  \label{recons_flux6}
  |S_m^c|\,\ddt{\ov{u}_m^c}=-\Int_{S_m^c} \divx{\bs{F}_h^c}\;\dd V+\Int_{\partial \omega_c}\hspace*{-2mm}\phi_m\;\(\bs{F}_h^c\pds\bs{n}-\mc{F}_n\)\,\dd S.
\end{align}
Finally, we obtain the following governing equation of the submean values time evolution
\begin{align}
  \label{recons_flux7}
  \ddt{\ov{u}_m^c}=-\Inv{|S_m^c|}\(\Int_{\partial S_m^c}\hspace*{-2mm}\bs{F}_h^c\pds\bs{n}\;\,\dd S-\Int_{\partial \omega_c}\hspace*{-2mm}\phi_m\;\(\bs{F}_h^c\pds\bs{n}-\mc{F}_n\)\,\dd S\).
\end{align}
The use of the reconstructed flux definition \eqref{recons_flux} directly leads to the following relation between the different types of fluxes
\begin{align}
  \label{recons_flux8}
  \Int_{\partial S_m^c}\wh{F_n}\;\dd S=\Int_{\partial S_m^c}\hspace*{-2mm}\bs{F}_h^c\pds\bs{n}\;\,\dd S-\Int_{\partial \omega_c}\hspace*{-2mm}\phi_m\;\(\bs{F}_h^c\pds\bs{n}-\mc{F}_n\)\,\dd S.
\end{align}
One can see how the reconstructed flux is connected to the interior polynomial flux and the jump at the cell interface between the interior flux and the DG numerical flux. To further develop relation \eqref{recons_flux8}, we make the same assumption on the cell boundary reconstructed fluxes as before \eqref{recons_flux_boundary}. Expression \eqref{recons_flux8} then rewrites
\begin{align}
  \label{recons_flux9}
  \Int_{\partial S_m^c\setminus\partial\omega_c}\wh{F_n}\;\dd S=\Int_{\partial S_m^c\setminus\partial\omega_c}\hspace*{-2mm}\bs{F}_h^c\pds\bs{n}\,\;\dd S-\Int_{\partial \omega_c}\hspace*{-2mm}\wt{\phi_m}\;\(\bs{F}_h^c\pds\bs{n}-\mc{F}_n\)\,\dd S,
\end{align}
where, for sake of  conciseness, $\wt{\phi_m}$ reads as follows
\begin{align}
  \label{subreso_modif}
  \wt{\phi_m}=\left\{\begin{array}{ll}
          \phi_m \qquad &\text{if }\; \bs{x}\in\,\partial\omega_c\setminus\partial S_m^c,\\[3mm]
          \phi_m -1 &\text{if }\; \bs{x}\in\,\partial\omega_c\cap\partial S_m^c.
          \end{array}\right.
\end{align}
Similarly to \eqref{recons_flux_int}, let $F_{mp}$ be the face integrated integrated value of the polynomial interior flux
\begin{align}
  \label{interior_flux_int}
  \Int_{f_{mp}^c} \bs{F}_h^c\pds\bs{n}\,\;\dd S=\veps^c_{mp}\,F_{mp}.
\end{align}
Then, if $F_c$ is the vector containing all the interior faces fluxes, one gets
\begin{align}
  \label{recons_flux10}
  A_c\,\wh{F_c}=A_c\,F_c-G_c,
\end{align}
where the $G_c$ contains the boundary contribution as
\begin{align}
  \label{recons_flux11}
  (G_c)_m=\Int_{\partial \omega_c}\hspace*{-2mm}\wt{\phi_m}\;\(\bs{F}_h^c\pds\bs{n}-\mc{F}_n\)\,\dd S.
\end{align}

Finally, by means of the same graph Laplacian technique used previously, we are able to express the reconstructed flux through the interior flux and a boundary correction term 
\begin{align}
  \label{recons_flux_definition2}
  \wh{F_c}=F_c-A_c\tra\,\mc{L}_c^{-1}\,G_c.
\end{align}
\begin{remark}
\label{remark_recons_flux_flx}
We can rewrite \eqref{recons_flux_definition2} as $\wh{F_c}=F_c-E\hspace*{-0.75mm}\(\bs{F}_h^c\pds\bs{n}-\mc{F}_n\)$, where $E(.)$ would be a correction function taking into account the jump between the polynomial flux and the numerical flux on the cell boundary. This permits to demonstrate once more that in DG schemes the numerical diffusion deriving from the jump in term of flux across the cell interfaces is distributing elsewhere inside the cell, here at the subcells faces. The sub-resolution basis functions act as weighted functions in the diffusion distribution.\\
Let us note that another choice in the correction term function $E(.)$ would lead to a different scheme. For instance, setting $E(.)=0$ leads to spectral volume scheme of Z.J. Wang \cite{Wang_SV_2002_I}.
\end{remark}\vspace*{2mm}

As previously, let us  gather the previous results in another theorem.
\begin{theorem}
  \label{main_thm_flux}
  Discontinuous Galerkin schemes, expressed in the following equation in cell $\omega_c$ through volume and boundary flux contribution
  \begin{align}
    \label{DG_2D_thm_flux}
    \Int_{\omega_c} \dt{u_h^c}\,\psi\,\dd V=\Int_{\omega_c} \bs{F}(u_h^c)\pds\gradx{\psi}\,\dd V-\Int_{\partial \omega_c}\hspace*{-2mm}\psi\;\mc{F}_n\,\dd S, \qquad \forall \,\psi\in\P^k(\omega_c),
  \end{align}
  can be recast into $N_k$ subcell finite volume like schemes as
  \begin{align}
    \label{sub_FV_thm_flux}
  \ddt{\ov{u}_m^c}=-\Inv{|S_m^c|}\,\Sum_{S_{p}^v\,\in \,\mc{V}_m^{\,c}}\veps^c_{mp}\,\wh{F_{mp}},
  \end{align}
  where the finite volume fluxes $\wh{F_{mp}}$, referred to as reconstructed fluxes, are either defined in equation \eqref{recons_flux_definition2} if $S_p^v\subset\omega_c$ and in \eqref{recons_flux_boundary} otherwise.
\end{theorem}

Both theorems \ref{main_thm_residual} and \ref{main_thm_flux} are perfectly equivalent. While the definition of DG reconstructed flux \eqref{recons_flux_definition2} enables a better comprehension of how those fluxes are related to the interior flux $\bs{F}(u_h^c)$ and how the numerical diffusion deriving from the jump at the interface is distributed inside the cell, the definition of reconstructed fluxes through residual \eqref{recons_flux_definition} is a lot easier to implement and do not require the definition of sub-resolution basis functions.

\section{\textit{A posteriori} local subcell correction}
\label{sect_DG_lim}

The previous reformulations of DG scheme into subcell FV-like scheme through the definition of reconstructed fluxes enable us to construct our \textit{a posteriori} local subcell correction (APLSC). In few words, the reconstructed fluxes $\wh{F}_{mp}$ will be modified in a robust way in subcells where the original DG scheme has failed. Let us emphasize that the popularity and the number of subcell correction techniques have extensively grown these past years, see \cite{peraire_2012,munz_subFV_2014,dumbser_subFV_lim_tri,vilar_aplsc_1D,Kuzmin_subcell_flux_limiting,pazner_2021,Gassner_subcell_shock_ES}. Those shock capturing and property preserving methods generally rely on a low-order scheme combined, at the subcell level, with a high-order one. However, let us note that in most cases, all the subcells contained in a cell will be impacted if something bad happened somewhere in the cell. In \cite{vilar_aplsc_1D}, we have introduce, for the one-dimensional case, a new technique permitting to correct the solution in a subcell without modifying the solution elsewhere. This particular feature allow to retain as much as possible the very precise subcell resolution of high-order DG schemes. The present paper aims at presenting the two-dimensional version of this correction. Let us emphasize that, up to our knowledge, this is the only technique working on totally unstructured grids and permitting the modification of the scheme, locally at the subcell level, without impacting the solution everywhere in the cell under consideration.\\

Let us mention that until now, only the semi-discrete version of schemes and their corresponding analysis were presented. To achieve high-accuracy in time, we make use of SSP Runge-Kutta time integration method \cite{Osher}. But, in the light of the fact that these multistage time integration methods write as convex combinations of first-order forward Euler scheme, the correction DG procedure will be presented for the simple case of this latter time numerical scheme, for sake of simplicity. DG schemes \eqref{DG_2D_1} provided with first-order forward Euler time integration writes, $\forall \,\psi\in\P^k(\omega_c)$
\begin{align}
  \label{DG_euler}
  \Int_{\omega_c} u_h^{c,\,n+1}\,\psi\,\dd V=\Int_{\omega_c} u_h^{c,\,n}\,\psi\,\dd V+\Dt\,\(\Int_{\omega_c} \bs{F}(u_h^{c,\,n})\pds\gradx{\psi}\,\dd V-\Int_{\partial \omega_c}\hspace*{-2mm}\psi\;\mc{F}_n\,\dd S\).
\end{align}

The numerical solution $u_h^{c,\,n}$ on cell $\omega_c$, being assumed to be a $\P^k$ polynomial, is uniquely defined through its mean values on the $N_k$ subcells, and reversely. To go from one representation to another, namely from $\{u_m^{c,\,n}\}_m$ the solution moments to $\{\ov{u}_m^{\, c,\, n}\}_m$ the subcell mean values, one has just to make use of the projection matrix $P_c$ or its inverse
\begin{align}
  \label{proj}
    P_c\,\begin{pmatrix}u_1^{c,\,n}\\\vdots\\u_{N_k}^{c,\,n}\end{pmatrix}=\begin{pmatrix}\ov{u}_1^{\,c,\,n}\\\vdots\\\ov{u}_{N_k}^{\,c,\,n}\end{pmatrix}.
\end{align}

Let us now introduce the correction procedure. First, we assume that at time $t^n$ the numerical solution $u_h^n$ is satisfactory in the sens that, on any cell $\omega_c$, the subcell mean values are admissible regarding some criteria yet to be defined. Then, we compute $u_h^{n+1}$ a candidate solution through the uncorrected DG scheme. The third step is then crucial. Indeed, we then have to check if the new uncorrected solution is admissible. If it is the case, we can go further in time without any special treatment. Otherwise we have to return to time $t^n$ and recompute the solution locally by means of a more robust scheme. This step is crucial in the sens that it will tell us if and where a new computation would be required.

\subsection{Troubled zone detector}
\label{subsect_detection}

Regarding the troubled zone detectors, we simply extend to the two-dimensional unstructured case the ones used in the one-dimensional frame \cite{vilar_aplsc_1D}. In this previous work, two detection criteria were mainly used, namely one ensuring the physical admissibility of the numerical solution (PAD) and another addressing the apparition of spurious oscillations (NAD). Let us then recall these two criteria.

\vspace*{3mm}
\paragraph{Physical admissibility detection (PAD)}
\begin{itemize}
\item Check if the different mean values $\ov{u}_m^{\,c,\,n+1}$ lie in a chosen convex physical admissible set (maximum principle for SCL, positivity of the pressure and density for Euler, \dots). Entropy stability can even be added to this admissible set.
\item Check if there is any $NaN$ values
\end{itemize}\vspace*{3mm}

Those are the minimum requirements if one wants to enforce code robustness. Now, in order to tackle the issue of spurious oscillations, we make use of a local maximum principle. Indeed, through the respect of the CFL, the solution in cell $\omega_c$ at time $t^{n+1}$ has to remain in the bounds of the solution at the previous time step $t^n$ wherein $\omega_c$ and its first face neighbors. This condition is reformulated in the following detection criterion.

\vspace*{3mm}
\paragraph{Numerical admissibility detection (NAD)}
\begin{itemize}
\item Check if the following Discrete Maximum Principle (DMP) on submean values is ensured:
  $$\Min_{v\,\in \,\mc{N}(S_m^c)}\(\ov{u}_v^{\,n}\)-\delta \leq \ov{u}_m^{\,c,n+1} \leq \Max_{v\,\in \,\mc{N}(S_m^c)}\(\ov{u}_v^{\,n}\)+\delta,$$
  where $\mc{N}(S_m^c)$ is some set of $S_m^c$ neighboring subcells, including subcell $S_m^c$. The parameter $\delta$, although not required, enables a small relaxation of the NAD criterion. Similarly to \cite{dumbser_subFV_lim_tri}, we set $\delta$ as follows
  \begin{align}
    \label{delta_NAD}
    \delta=\max\(10^{-4},\; 10^{-3}\,\(\Max_{v\,\in \,\mc{N}(S_m^c)}\(\ov{u}_v^{\,n}\)-\Min_{v\,\in \,\mc{N}(S_m^c)}\(\ov{u}_v^{\,n}\)\)\).
  \end{align}
\end{itemize}
\vspace*{3mm}

\begin{remark}
  The bigger the set $\,\mc{N}(S_m^c)$ is, the softer the NAD criterion will be. Reversely, a smaller set would lead to more subcells considered as problematic. Different sets will be considered for linear and non-linear problems, see Sections~\ref{subsect_linear} and \ref{subsect_nonlinear}.
\end{remark}
\vspace*{2mm}

Let us enlighten that NAD criterion relies on a maximum principle based on subcell mean values. And because this maximum principle is not concerned with the whole polynomial set of values, it is very well-known that one has to relax it to preserve scheme accuracy in the presence of smooth extrema. To do so, we make use of the two-dimensional version of the one introduced in the 1D case, \cite{vilar_aplsc_1D}.

\paragraph{Detection of smooth extrema}

This smooth extrema detection criterion is based on an idea present in different limitations, as the hierarchical slope limiter \cite{Kuzmin09}. In this limitation technique, the numerical solution is supposed to exhibit a smooth extrema if at least the linearized version of the numerical solution spatial derivatives, $i.e.$
\begin{subequations}
\label{grad_uh}
\begin{empheq}[left = \empheqlbrace\,]{align}
&v_x^{\,c}(\bs{x})=\ov{\vdx \,u_h^{c,\,n+1}}+\ov{\gradx\, (\vdx\, u_h^{c,\,n+1})}\pds (\bs{x}-\bs{x}_c), \label{dx_uh}\\[3mm]
&v_y^{\,c}(\bs{x})=\ov{\vdy \,u_h^{c,\,n+1}}+\ov{\gradx\, (\vdy\, u_h^{c,\,n+1})}\pds (\bs{x}-\bs{x}_c), \label{dy_uh}
\end{empheq}\\[-4mm]
\end{subequations}

present a monotonous profile. In \eqref{grad_uh}, $\bs{x}_c$ denotes the centroid of cell $\omega_c$, while $\ov{\partial_\xory \,u_h^{c,\,n+1}}$ and $\ov{\gradx\, (\partial_\xory\, u_h^{c,\,n+1})}$ are nothing but the averaged values on $\omega_c$ of the successive partial derivatives of $u_h^c$. In practice, the NAD relaxation used here works as a vertex-based limiter on $v_\xory^{\,c}(\bs{x})$. Due to their linearity, functions $v_\xory^{\,c}(\bs{x})$ attain their extrema at the vertices $\bs{x}_p\in\mc{P}_c$, where $\mc{P}_c$ stands for the set of vertices of cell $\omega_c$. Then, we consider that the exact weak solution underlying the numerical solution $u_h$ presents a smooth profile in cell $\omega_c$ if, for any vertex $\bs{x}_p\in\mc{P}_c$, the linearized spatial derivative functions ensure the following constraint
\begin{align}
  \label{smooth_constraint}
  v_{x,\,p}^{\min}\leq v_x^c(\bs{x_p})\leq v_{x,\,p}^{\max} \qquad \text{and} \qquad v_{y,\,p}^{\min}\leq v_y^c(\bs{x_p})\leq v_{y,\,p}^{\max},
\end{align}

where $v_{\xory,\,p}^{\min}$ and $v_{\xory,\,p}^{\max}$ are respectively defined as
\begin{align}
  \label{min_max}
  v_{\xory,\,p}^{\min}=\Min_{v\in \mc{V}_p} v_\xory^{\,c}(\bs{x}_p)\qquad \text{and} \qquad v_{\xory,\,p}^{\max}=\Max_{v\in \mc{V}_p} v_\xory^{\,c}(\bs{x}_p).
\end{align}

In \eqref{min_max}, $\mc{V}_p$ represents the set of cells that share $\bs{x}_p$ as a vertex, $i.e. \;\; \omega_v\in\mc{V}_p\; \Longrightarrow\; \bs{x}_p\in \mc{P}_v$. Practically, if for any vertex $\bs{x}_p\in\mc{P}_c$, conditions \eqref{smooth_constraint} are ensured, we then consider that the numerical solution presents a smooth profile on cell $\omega_c$. In this particular case, the NAD criterion is relaxed allowing the preservation of smooth extrema along with the order of accuracy for smooth problems, see Section~\ref{sect_results}.

\begin{remark}
  \label{smooth_detection}
  We have presented here the detection based on the linearized first-order derivatives of the solution. This would work for any higher order derivative. Furthermore, such relaxation procedure can also be applied at the subcell level to also preserve smooth extrema even within a cell at a smaller length scale. This would be useful in the context of coarse grids. Actually, because the subcell smooth extrema relaxation technique works well for both coarse and fine meshes, this will be the procedure used for the numerical applications, Section~\ref{sect_results}.
\end{remark}\vspace*{2mm}

Now that we have detailed the troubled subcell detector, the correction procedure will be presented in the next subsection.

\subsection{Correction}
\label{subsect_correction}

In this subsection, the original correction procedure introduced in the one-dimensional frame in \cite{vilar_aplsc_1D} will firstly be extended to the two-dimensional unstructured grid case. Then, another technique blending in a convex combination fashion DG reconstructed fluxes and first-order FV fluxes for admissible subcells in the vicinity of troubled areas will also be presented and prove to improve results in comparison to the original algorithm.

\subsection{Original correction principle}
\label{subsubsect_original_correction}

\begin{figure}[!ht]
  \begin{center}
    \subfigure[Structured subdivision.]{\includegraphics[height=5.cm]{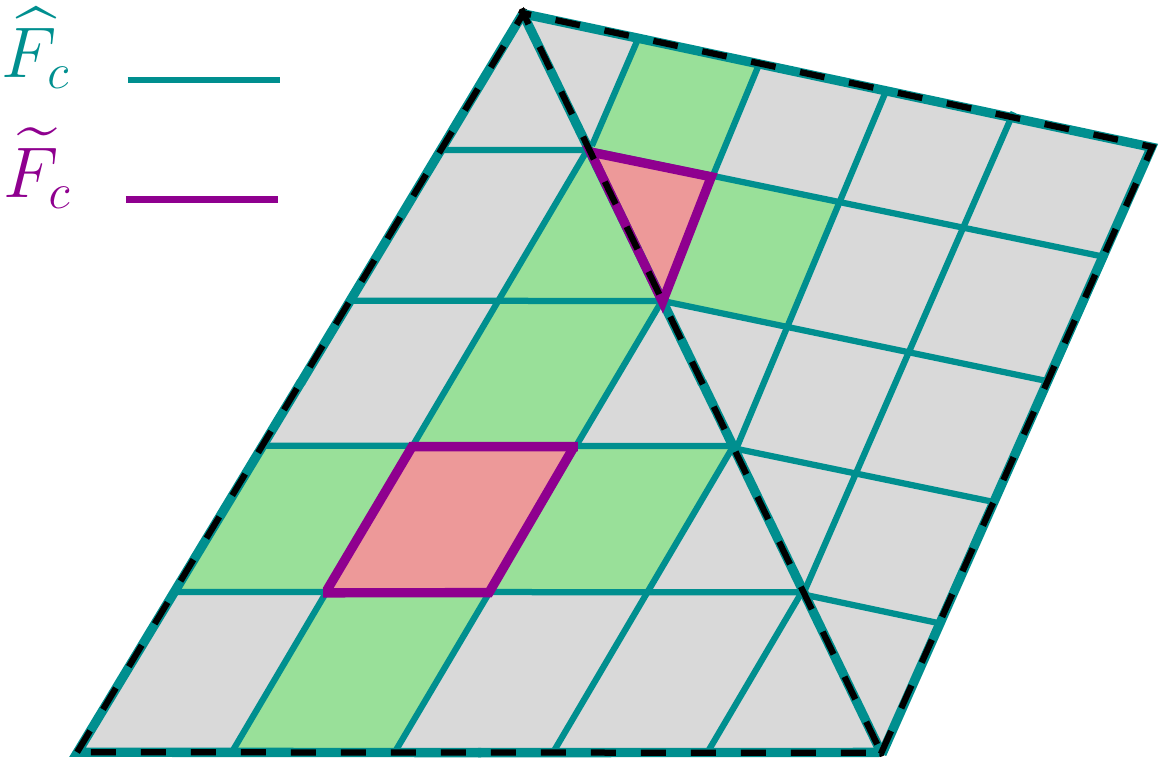}\label{fig_corr_cart}}\hspace*{-4mm}
    \subfigure[Voronoi-type subdivision.]{\includegraphics[height=4.8cm]{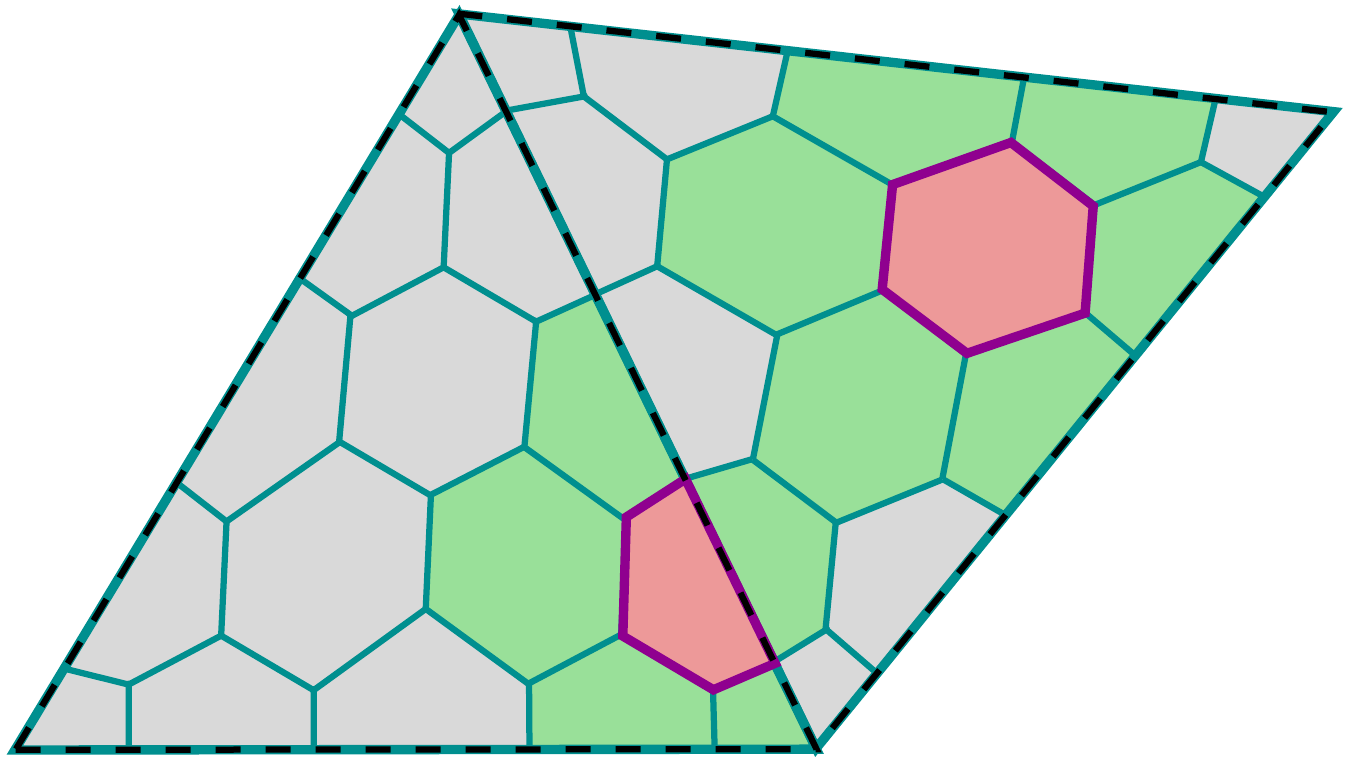}\label{fig_corr}}
    \caption{Original correction of the DG reconstructed flux.}
  \label{fig_flux_lim}
  \end{center}
\end{figure}

The very simple idea that forms the basis of the original correction procedure introduced in \cite{vilar_aplsc_1D} is the following: if the uncorrected DG scheme has produced a numerical solution $u_h^{\,c,\,n+1}$ on cell $\omega_c$, which is not admissible in subcell $S_m^c$ in regards to the detection criteria presented previously, the subcell mean value $\ov{u}_m^{\,c,\,n+1}$ will be recomputed by means of a more robust scheme. To do so, and because uncorrected DG scheme is equivalent to subcell finite volume scheme with the appropriate high-order reconstructed fluxes, see Theorems~\ref{main_thm_residual} and \ref{main_thm_flux}, we substitute on the boundaries of subcell $S_m^c$ the high-order reconstructed fluxes with some first-order finite volume numerical fluxes. The submean value $\ov{u}_m^{\,c,\,n+1}$ will then be recomputed by means of a simple and robust first-order finite volume scheme. This concept is depicted in Figure~\ref{fig_flux_lim}, where the troubled subcells are colored red.

Because it is of fundamental importance to preserve scheme conservation, the first face neighboring subcells, colored green in Figure~\ref{fig_flux_lim}, have to be also recomputed since we have modify the reconstructed fluxes on the boundary of the troubled subcell. The submean values of the neighboring subcells are then computed through a FV-like scheme with first-order numerical flux on one or more faces and high-order reconstructed fluxes on the remaining interfaces. For the remaining subcells, colored gray in Figure~\ref{fig_flux_lim}, because the corresponding reconstructed fluxes have not been modified, there is no need to recompute them. The corresponding submean values are hence the values obtained through the uncorrected DG scheme. This original concept of the \textit{a posteriori} local subcell correction of DG scheme is summarized in the following flowchart:\vspace*{3mm}

\begin{enumerate}
\item Compute the candidate solution $u_h^{\,c,\,n+1}$ by means of the uncorrected DG scheme.\\[0mm]
\item Project $u_h^{\,c,\,n+1}$ through \eqref{proj} to get the $N_k$ submean values $\ov{u}_m^{\,c,\,n+1}$.\\[0mm]
\item Check $\ov{u}_m^{\,c,\,n+1}$ through the troubled zone detection criteria plus relaxation \label{pt_check}.\\[0mm]
\item If all submean values $\ov{u}_m^{\,c,\,n+1}$ are admissible, go further in time. Otherwise, modify the corresponding reconstructed flux value on the troubled subcells faces through a first-order numerical flux as following
  \begin{alignat*}{2}[left = \empheqlbrace\,]
    &\wt{F_{mp}}=\veps^c_{mp}\;l_{mp}^c\; \mc{F}\(\ov{u}_m^{\,c,n},\ov{u}_p^{\,v,n},\bs{n}_{mp}\) \qquad &\text{if } S_{m}^c \text{ or } S_{p}^v\in\mc{V}_m^c \text{ is either marked},\\[2mm]
    &\wt{F_{mp}}=\wh{F_{mp}} & \text{otherwise},
  \end{alignat*}
  where $l_{mp}^c$ is the length of face $f_{mp}^c$ .\\
  
\item Through the corrected reconstructed flux, recompute the submean values for tagged subcells and their first neighboring subcells as
  \begin{align*}
    \ov{u}_m^{\,c,\,n+1}=\ov{u}_m^{\,c,\,n}-\Frac{\Dt}{|S_m^c|}\,\Sum_{S_{p}^v\,\in \,\mc{V}_m^{\,c}}\veps^c_{mp}\,\wt{F_{mp}}.
  \end{align*}
\item By means of $P_c^{-1}$ and equation \eqref{proj}, get the new corrected polynomial solution $u_h^{\,c,\,n+1}$.\\[0mm]
\item Return to point~\ref{pt_check}.
\end{enumerate}\vspace*{4mm}

In light of this correction procedure flowchart, it is clear that the DG solution will only be affected at the subcell scale. Furthermore, the corrected scheme is conservative at the subcell level by construction.\vspace*{2mm}

\begin{remark}
  \label{rmk_correction}
  In the present correction, when needed we substitute the reconstructed flux value $\wh{F_{mp}}$ with a first-order numerical flux $\mc{F}\(\ov{u}_m^{\,c,n},\ov{u}_p^{\,v,n},\bs{n}_{mp}\)$. Obviously, other choices are possible and may even be more appropriate, as for instance a second-order TVD numerical flux or even a WENO numerical flux. In the 1D, results with both first-order and second-order TVD corrections were presented, see \cite{vilar_aplsc_1D}. In the present article, for sake of simplicity and conciseness, only the first-order correction previously detailed will be used.
\end{remark}\vspace*{2mm}

\begin{remark}
  \label{remark_positivity}
  Let us emphasize that since this \textit{a posteriori} correction is based on first-order finite volume scheme, maximum principle or positivity preservation in the case of systems are enforced by construction. However, it is crucial to note that to start from a discrete representation of the initial datum that respects its bounds (or ensures positivity), the initialization has to be carried out through the integration of the initial solution on the subcells to obtain the different submean values, and then through relation \eqref{proj} compute the corresponding polynomial moments of the solution. Traditionally, in DG schemes the initialization is handled by either assigning at the solution points the initial datum value to the numerical solution, or by a $L_2$ projection onto $\P^k$. In either of these procedures, nothing ensures that the submean values would respect the bounds of the initial datum.
\end{remark}\vspace*{2mm}

Let us recall that in the 1D case \cite{vilar_aplsc_1D}, for non-linear problems using very high-order schemes on coarse grids, the numerical solution proved to remain slightly oscillatory at the subcell level. To overcome this issue, we were artificially enlarging the stencil to correct by also marking, for $k\geq 3$, the first face neighboring subcells of a troubled subcell. Now, in this 2D frame, a modified algorithm making use of convex combination of DG reconstructed fluxes and first-order FV fluxes for admissible subcells in the vicinity of troubled areas will be presented and prove to produce better results than the original correction, in this context of very high-orders and coarse meshes.

\subsection{New correction principle}
\label{subsubsect_new_correction}

\begin{figure}[!ht]
  \begin{center}
    \subfigure[Structured subdivision.]{\includegraphics[height=4.4cm]{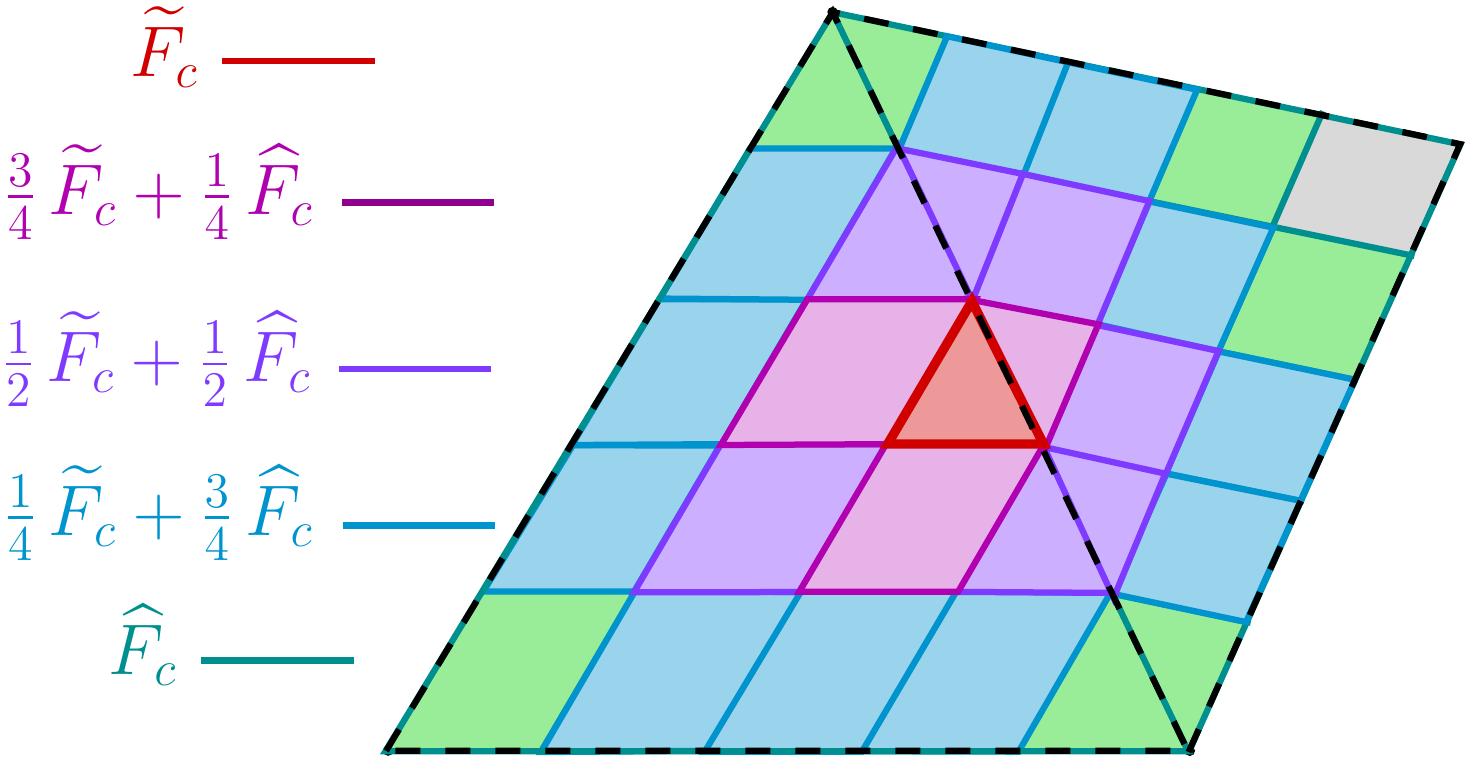}\label{fig_corr_cart_new}}\hspace*{-4mm}
    \subfigure[Voronoi-type subdivision.]{\includegraphics[height=4.3cm]{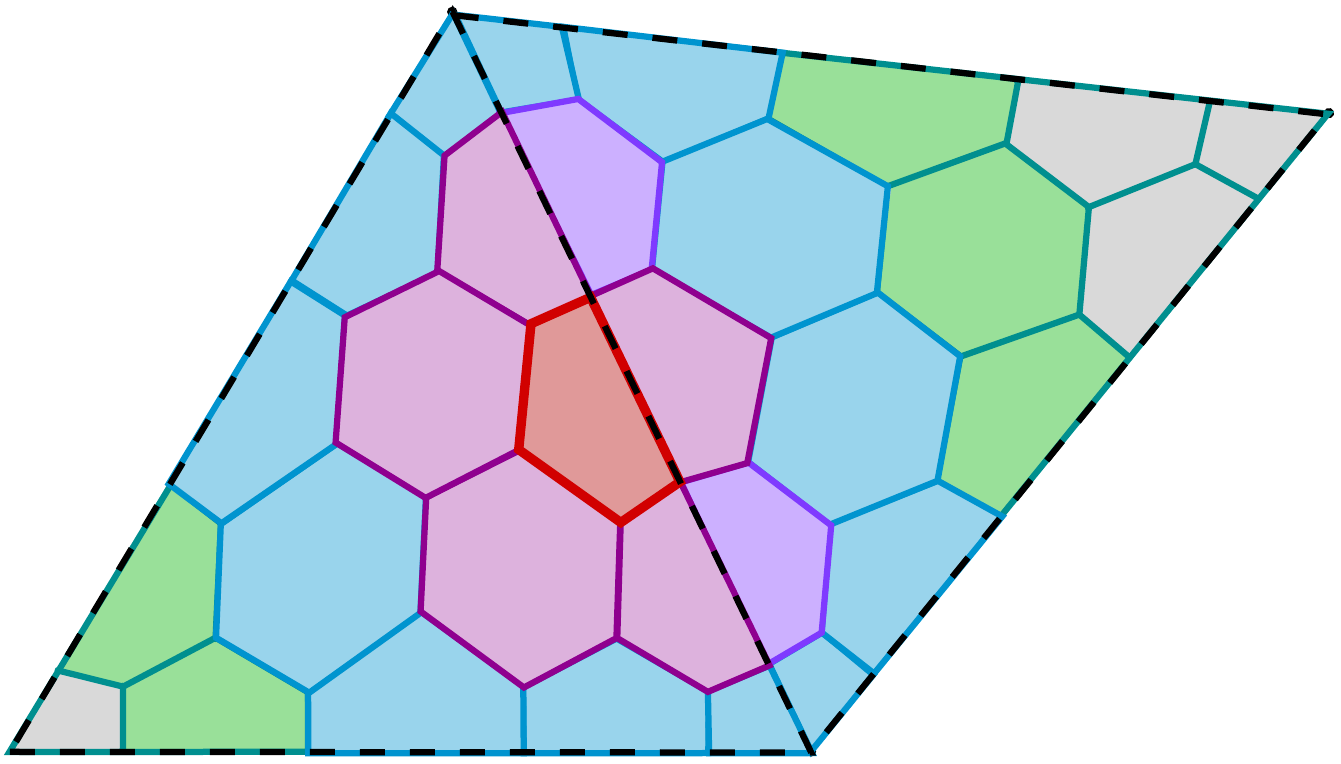}\label{fig_corr_new}}
    \caption{New correction of the DG reconstructed flux.}
  \label{fig_flux_lim_new}
  \end{center}
\end{figure}

As said previously, for very high-order DG methods on coarse grids, the corrected scheme may locally jump, in two subcells scale, from a very precise approximation to a robust but very low accurate first-order representation. In the particular situation, numerical experiments have proved that some subcell oscillations may remain due to this phenomenon. To overcome this issue, we will now introduce a new correction algorithm in which the stiffness of the order transition will be smoothened by means of convex combination between DG reconstructed fluxes and first-order FV fluxes as follows

\begin{align}
  \label{convex_combo}
  \wt{F_{mp}}=\theta_{mp}\;\veps^c_{mp}\,l_{mp}^c\, \mc{F}\(\ov{u}_m^{\,c,n},\ov{u}_p^{\,v,n},\bs{n}_{mp}\)+(1-\theta_{mp})\;\wh{F_{mp}},
\end{align}
where $\theta_{mp}$ is a function of the distance to a non-admissible subcell. Obviously, the farther from the troubled subcell, we are the less of first-order FV we want. Many smoothness indicators exist in the literature, see for instance \cite{vanL,Harten_ENO,burb01,kriv04,art_vis_02,weno_lim_04}, and could be used in this context to determine a relevant bending coefficient. Another way around could be to adopt flux limiting or Flux-Corrected transport (FCT) approach, see \cite{Boris_FCT,Zalesak_FCT,Kuzmin_FCT_book_1,Kuzmin_FCT_book_2,Kuzmin_FCT_limiting_2017,Kuzmin_FCT_limiting_2021,Kuzmin_FCT_limiting_2022}, to find the proper blending ensuring a discrete maximum principle. To remain as simple as possible, and because this paper is concerned with the introduction of the DG-FV equivalency and the basic correction principle on 2D unstructured grids, we make use here of a very naive procedure. The principle is the following, for marked subcells detected through the troubled subcell indicator, first-order FV numerical flux is used on its boundaries, $i.e.\; \theta_{mp}=1$ in \eqref{convex_combo}. Then, for its face neighboring subcells $S_p^v\in\mc{V}_m^c$, convex combination \eqref{convex_combo} with $\theta_{mp}=\frac{3}{4}$ would be chosen for their remaining boundaries. Now, introducing $\wt{\mc{V}_m^c}$ the set containing both the face and node neighboring subcells, fluxes on faces of $S_p^v\in\wt{\mc{V}_m^c}\setminus\mc{V}_m^c$ which yet have not been corrected will be defined through at the blending coefficient $\theta_{mp}=\frac{1}{2}$. Finally, the remaining face neighbors of $S_p^v\in\wt{\mc{V}_m^c}$ will see their remaining boundaries associated with a numerical flux calculated through \eqref{convex_combo} with $\theta_{mp}=\frac{1}{4}$. This naive technique is displayed in Figure~\ref{fig_flux_lim_new}.\\

\begin{figure}[!ht]
  \begin{center}
    \subfigure[$t=0$.]{\includegraphics[height=6.7cm]{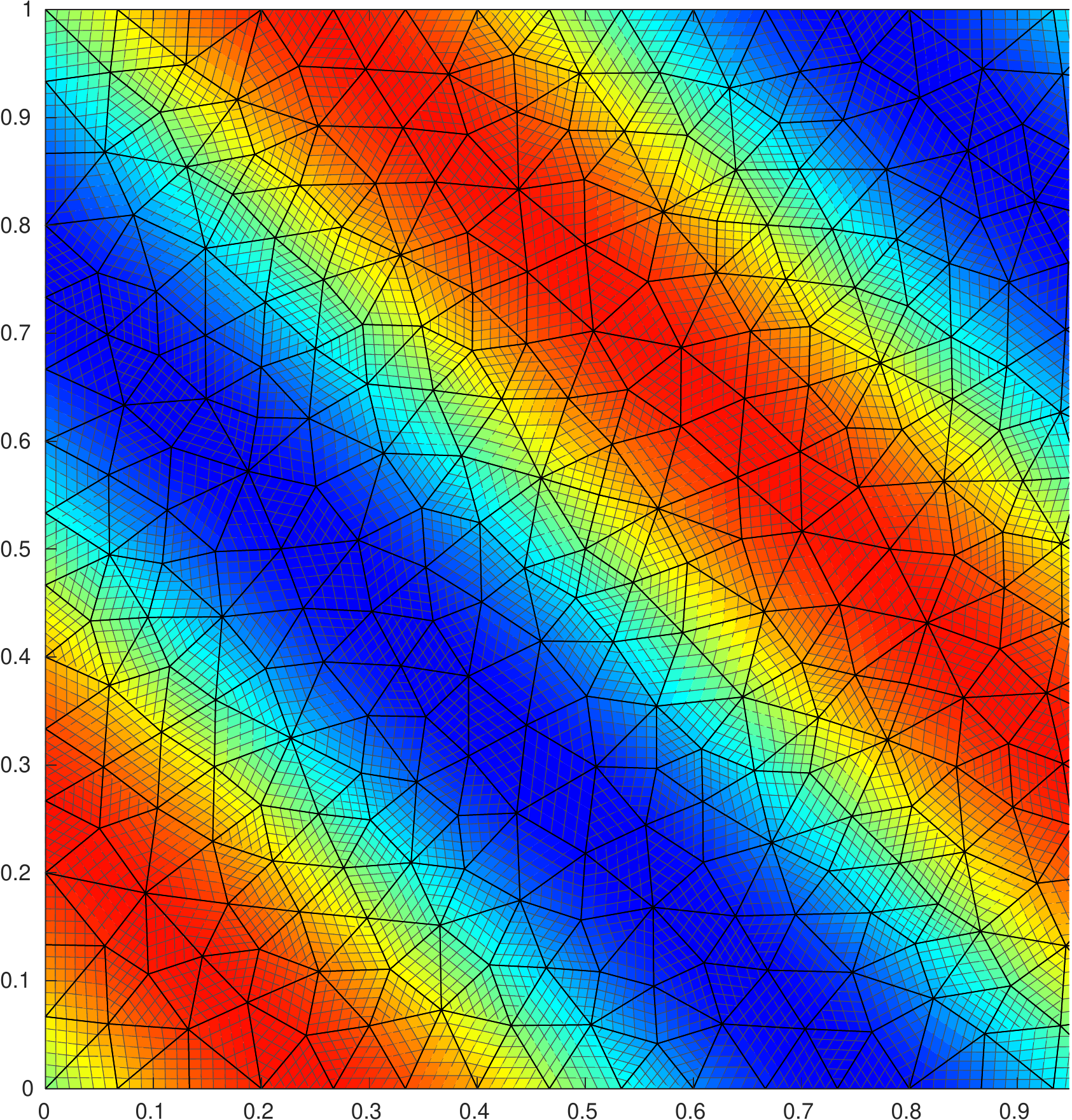}}\hspace*{5.mm}
    \subfigure[$t=0.5$.]{\includegraphics[height=6.7cm]{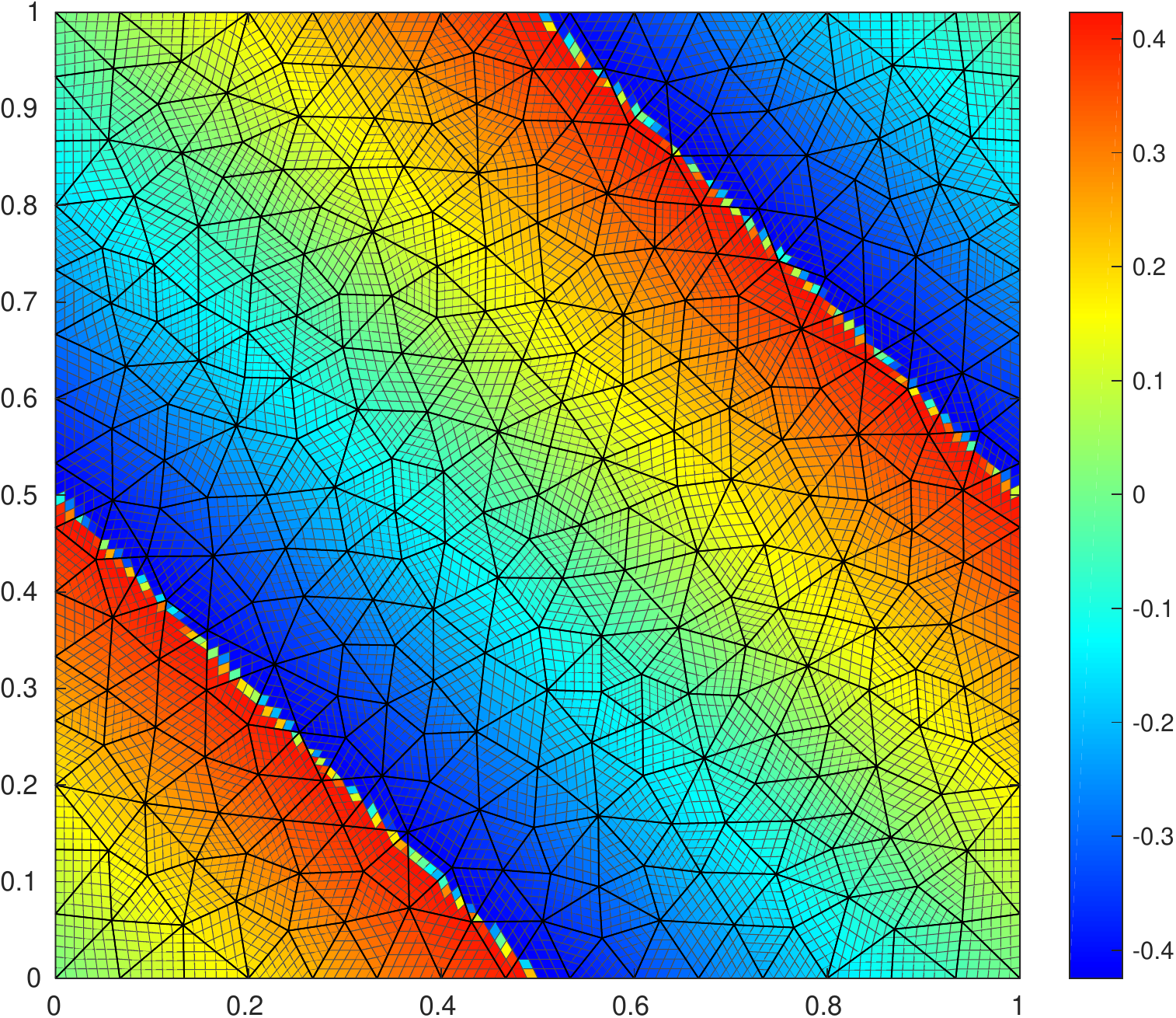}}
    \caption{Subcell mean values of the entropic weak solution of the Burgers equation with $u_0(\bs{x})=\sin(2 \pi \,(x+y))$.}
    \label{fig_burgers_exact}
  \end{center}
\end{figure}

\begin{figure}[!ht]
  \begin{center}
    \subfigure[Original correction.]{\includegraphics[height=7.cm]{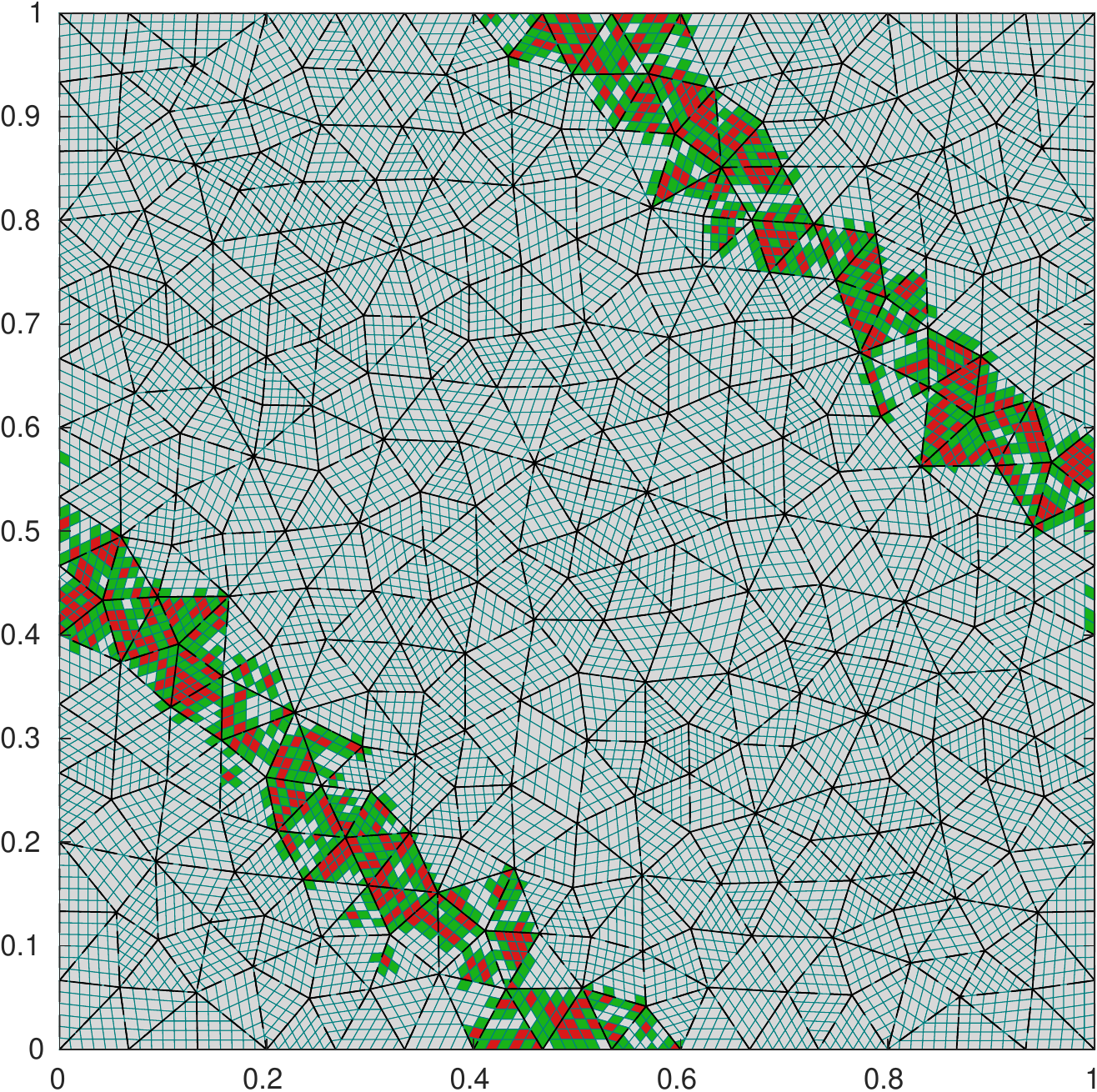}\label{corr1_detect}}\hspace*{12.mm}
    \subfigure[New correction.]{\includegraphics[height=7.cm]{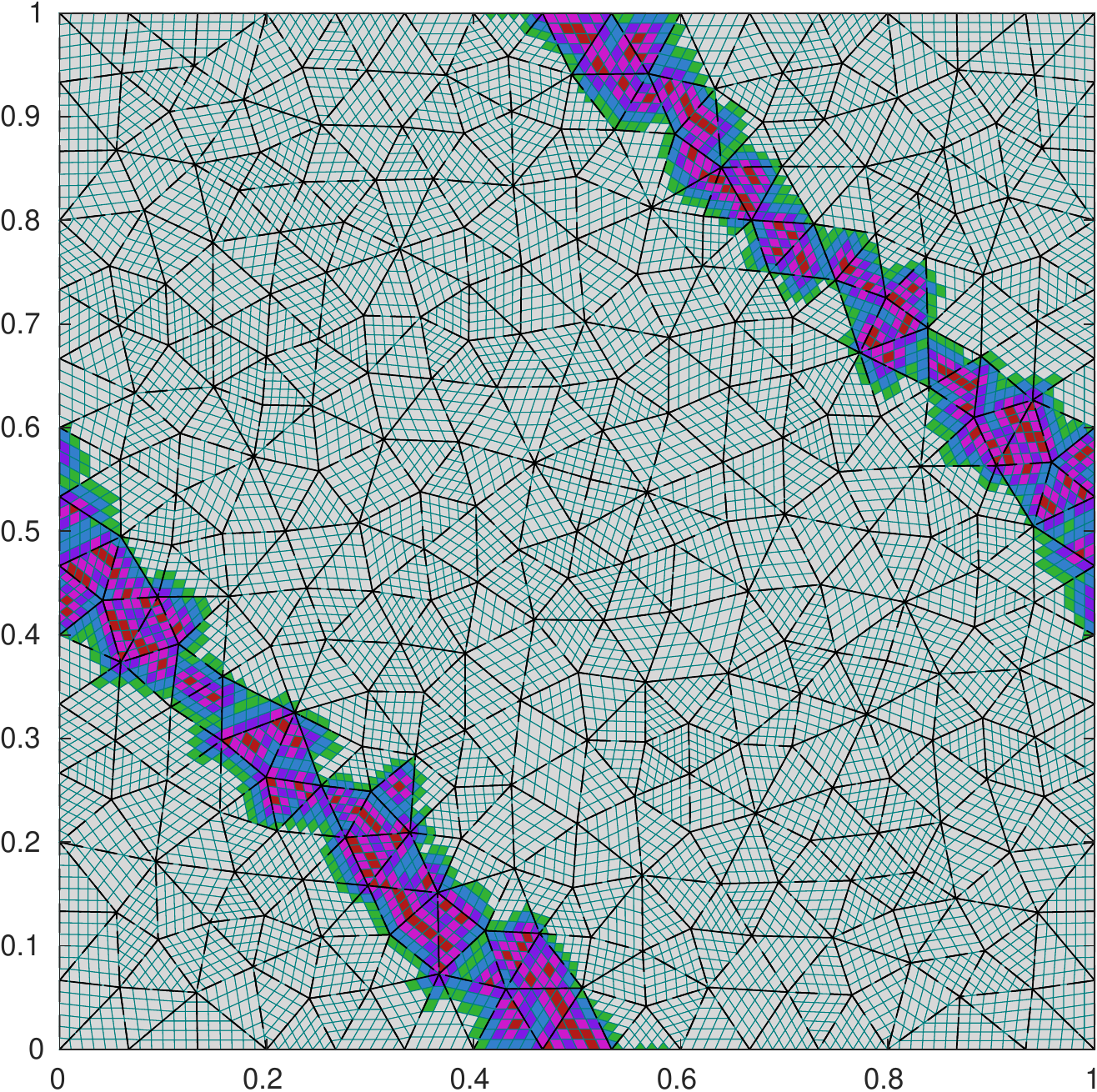}\label{corr2_detect}}
    \caption{Comparison between original and new correction procedure: corrected subcells.}
    \label{fig_burgers_compar_detect}
  \end{center}
\end{figure}

One could think that doing so, we have substantially enlarged the stencil of subcells to be corrected, and thus reduce the simulation code efficiency. However, in practice it is generally not the case, as one can see in Figure~\ref{fig_burgers_compar_detect}. Moreover, the computational time will also be generally slightly reduced. Emphasizing that the correction is done in an \textit{a posteriori} fashion and thus has to be potentially iterated multiple times at a time step to reach an admissible solution (see point 7. of the previous flowchart), if the stiff original correction introduces small oscillations at the subcell scale, the stencil will be automatically enlarged along with an increase of the computational time, through the correction iteration. With the new correction principle, only one iteration is generally needed. So even if at first, the set containing marked subcells is smaller in the original correction, it will end up with approximately the same size than with this new approach, and required more iterations and thus more computational effort. To assess this matter and compare both corrections, we make use of the Burgers equation, defined through \eqref{lcs1} and flux function $\bs{F}(u)=\demi\(u^2,\,u^2\)\tra$, with the smooth initial solution $u_0(\bs{x})=\sin(2 \pi \,(x+y))$. The domain is chosen as the unit square $[0,\,1]^2$ with periodic boundary condition. Through time, the exact solution will exhibit two stationary shocks along the lines defined by $\(\bs{x}\in[0,\,1]^2,\;x+y=0.5\)$ and $\(\bs{x}\in[0,\,1]^2,\;x+y=1.5\)$, see Figure~\ref{fig_burgers_exact}.\\

\begin{figure}[!ht]
  \begin{center}
    \subfigure[Original correction.]{\includegraphics[height=7.cm]{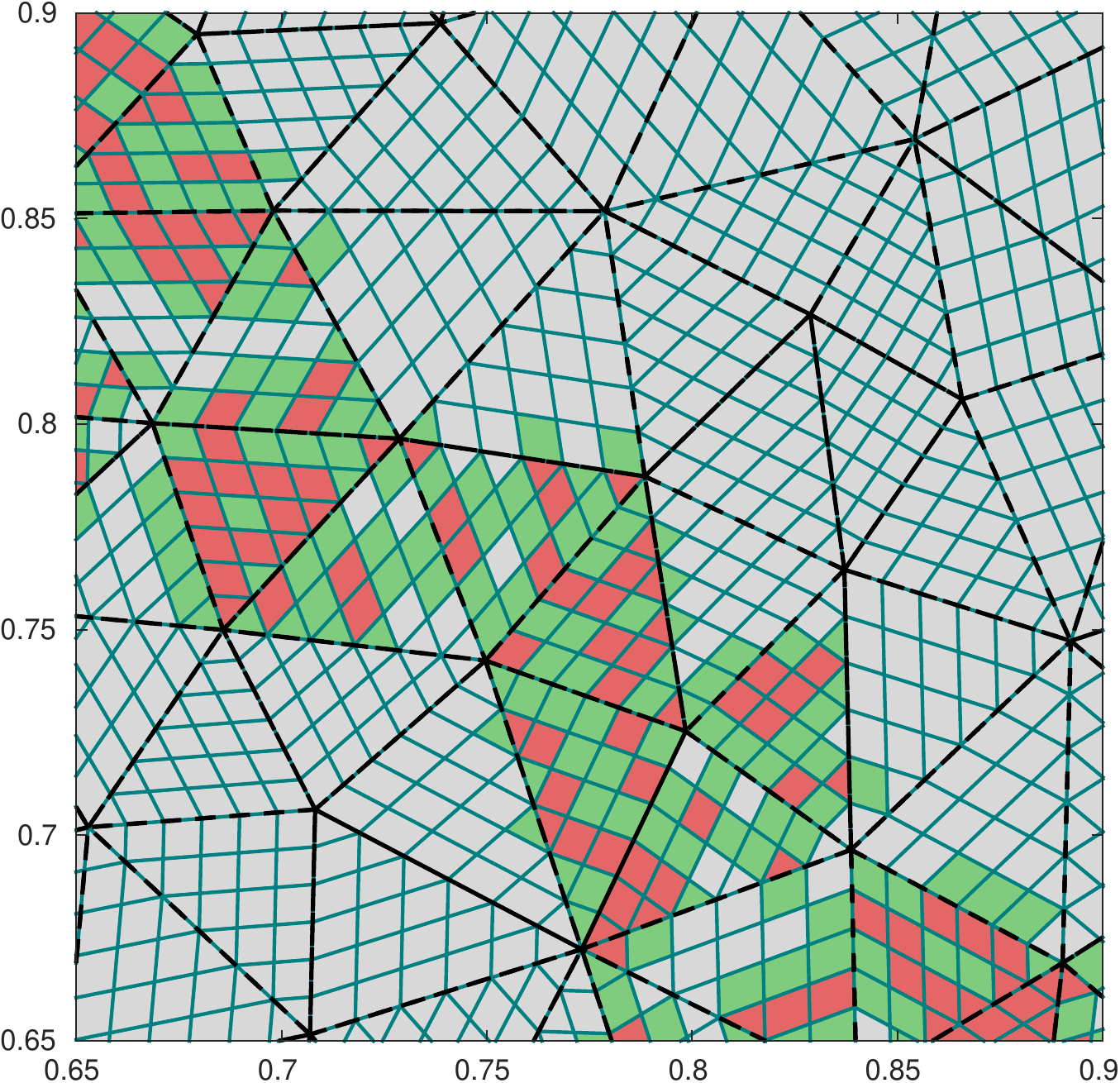}\label{corr1_detect_zoom}}\hspace*{10.mm}
    \subfigure[New correction.]{\includegraphics[height=7.cm]{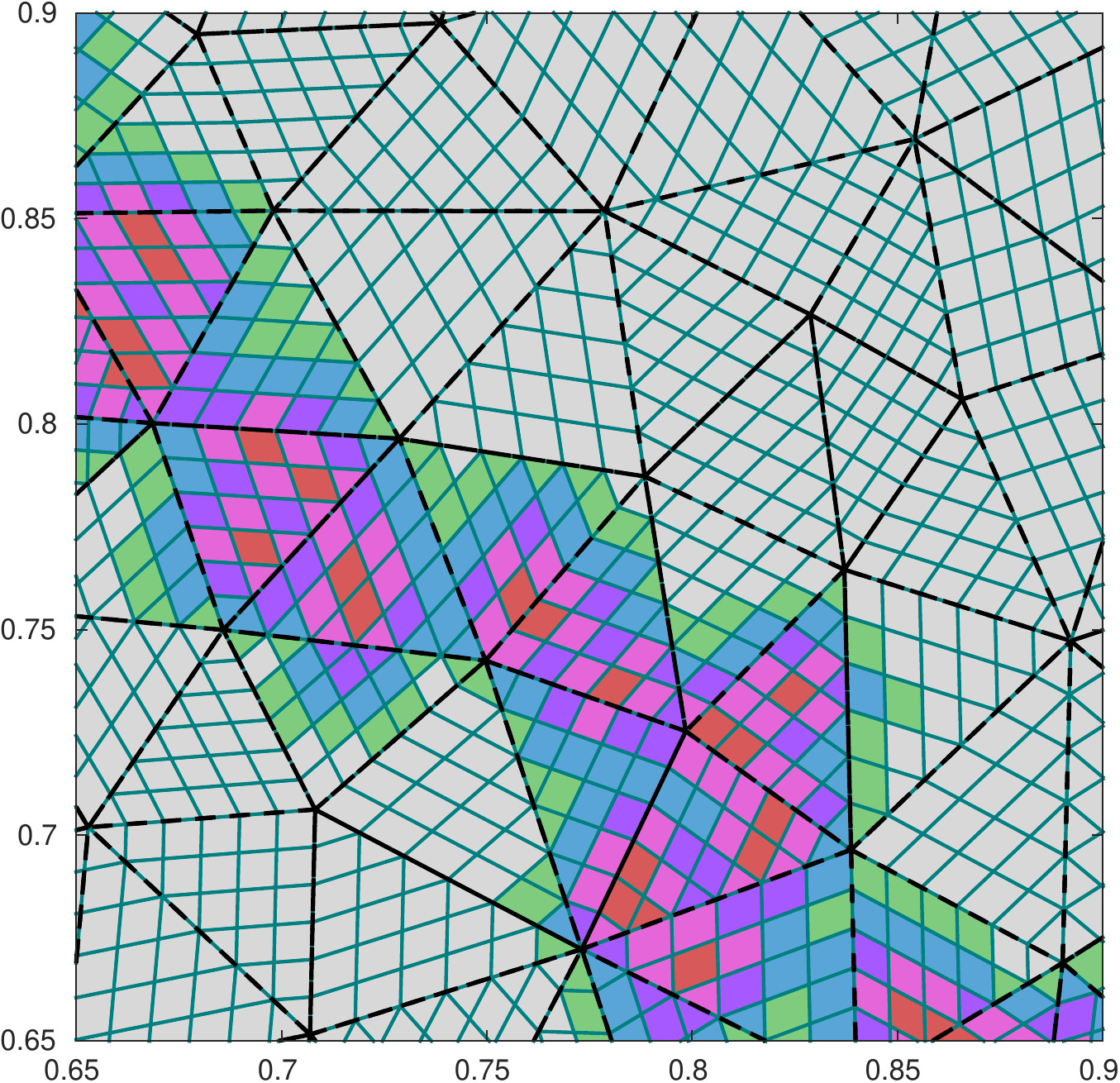}\label{corr2_detect_zoom}}
    \caption{Comparison between original and new correction procedure: corrected subcells - zoom on $[0.65, \,0.9]^2$.}
    \label{fig_burgers_compar_detect_zoom}
  \end{center}
\end{figure}

We run this test case with a sixth order DG scheme, on a quite coarse unstructured triangular grid made of 576 cells, with both \textit{a posteriori} local subcell corrections. We make use here of the simple structured subdivision introduced in Figure~\ref{fig_tri1}. To compare the two approaches, let us first display the marked subcells to be corrected. In Figures~\ref{fig_burgers_compar_detect} and \ref{fig_burgers_compar_detect_zoom}, we color all the subcells that have been corrected in the different Runge-Kutta steps during the last time step.

Firstly, one can note that both corrections, through the NAD criterion, have accurately capture the discontinuities, as the subcells to be corrected remain in a small vicinity of the shocks. Secondly, we can also observe how these corrections operate locally at the subcell scale. Finally, as depicted in Figure~\ref{fig_burgers_compar_detect}, the number of subcells to be corrected remains approximately the same with both approaches. Actually, with the original correction on average 10\% of the total number of subcells have to be corrected through the computation, while 14\% of the subcells with the new approach. But, even if indeed the number of subcells to be corrected has slightly increases through this new procedure, if we compare computational efficiency, it took 2 minutes and 22 seconds to the code with the original correction, and 2 minutes 15 seconds with the new one. As said previously, it comes from the fact that more correction iterations are required with the original procedure. For this calculation, on average 2.86 iterations are needed with the original approach when the correction has been triggered, for a maximum of 6 iterations, while only 1.46 iterations with the new approach, with a maximum of 3 iterations during the whole calculation.

For a better understanding of how work those two different approaches, we zoom on the area $[0.65, \,0.9]^2$. Let us first note that because the dimension of $\P^{\,5}$ is 20, each cell has been accordingly subdivided into 20 subcells. In the light of Figure~\ref{fig_burgers_compar_detect_zoom}, one can observe that a lot less subcells are corrected through a purely first-order finite volume with the new correction than with the original one, which enables even more the preservation of the high accurate subcell resolution of high-order DG schemes. In Figure~\ref{corr1_detect_zoom}, one can see the two types of subcells to be corrected, meaning the troubled subcells colored red and their face neighbors colored green which have to be recalculated to preserve the scheme conservation. In Figure~\ref{corr2_detect_zoom}, there are five types of subcells to be corrected, from the troubled subcells to be recomputed through a first-order FV scheme to the magenta, purple and blue ones which are recalculated through a convex combination between first-order FV numerical flux and the high-order DG reconstructed flux, with different weights. The green ones are the face neighbors of the previous subcells to be recomputed to preserve scheme conservation.

\begin{figure}[!ht]
  \begin{center}
    \subfigure[Original correction.]{\includegraphics[height=6.7cm]{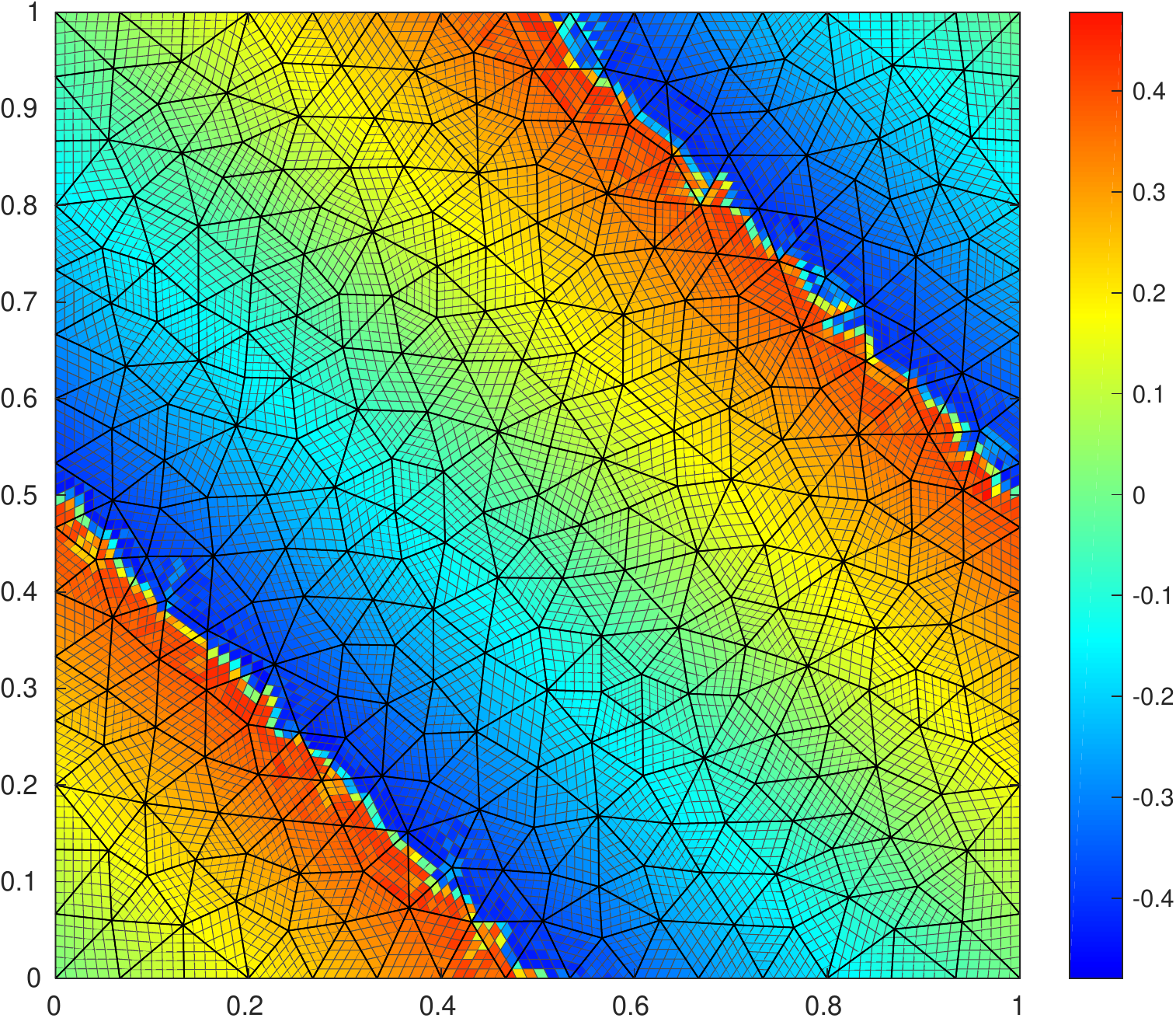}\label{corr1}}\hspace*{5.mm}
    \subfigure[New correction.]{\includegraphics[height=6.7cm]{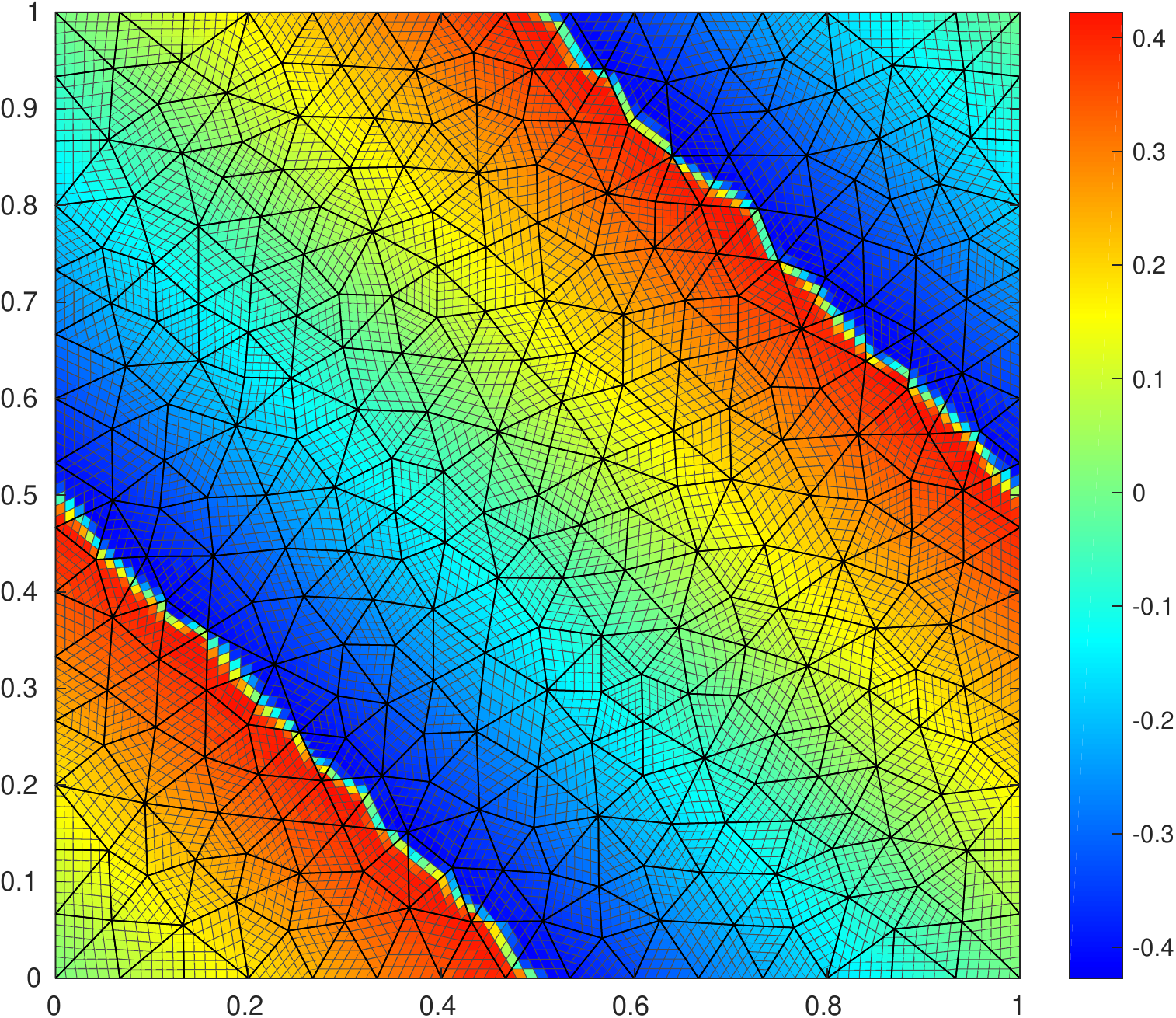}\label{corr2}}
    \caption{Comparison between original and new correction procedure: subcell mean values.}
    \label{fig_burgers_compar_sol}
  \end{center}
\end{figure}

While the previous results showed how those two approaches work and do have a quasi-equivalent computation cost, let us now compare the approximated solutions obtained and assess the benefit of the new correction. In Figure~\ref{fig_burgers_compar_sol}, the subcell mean values
obtained by means of a sixth-order DG scheme corrected through the original and the new procedures are displayed. Comparing Figures~\ref{corr1} and \ref{corr2}, we can see that the new correction has improved the accuracy of the scheme by a sharper and less oscillatory representation of the shock, despite the very coarseness of the mesh used.\\

\begin{figure}[!ht]
  \begin{center}
    \includegraphics[height=7.5cm]{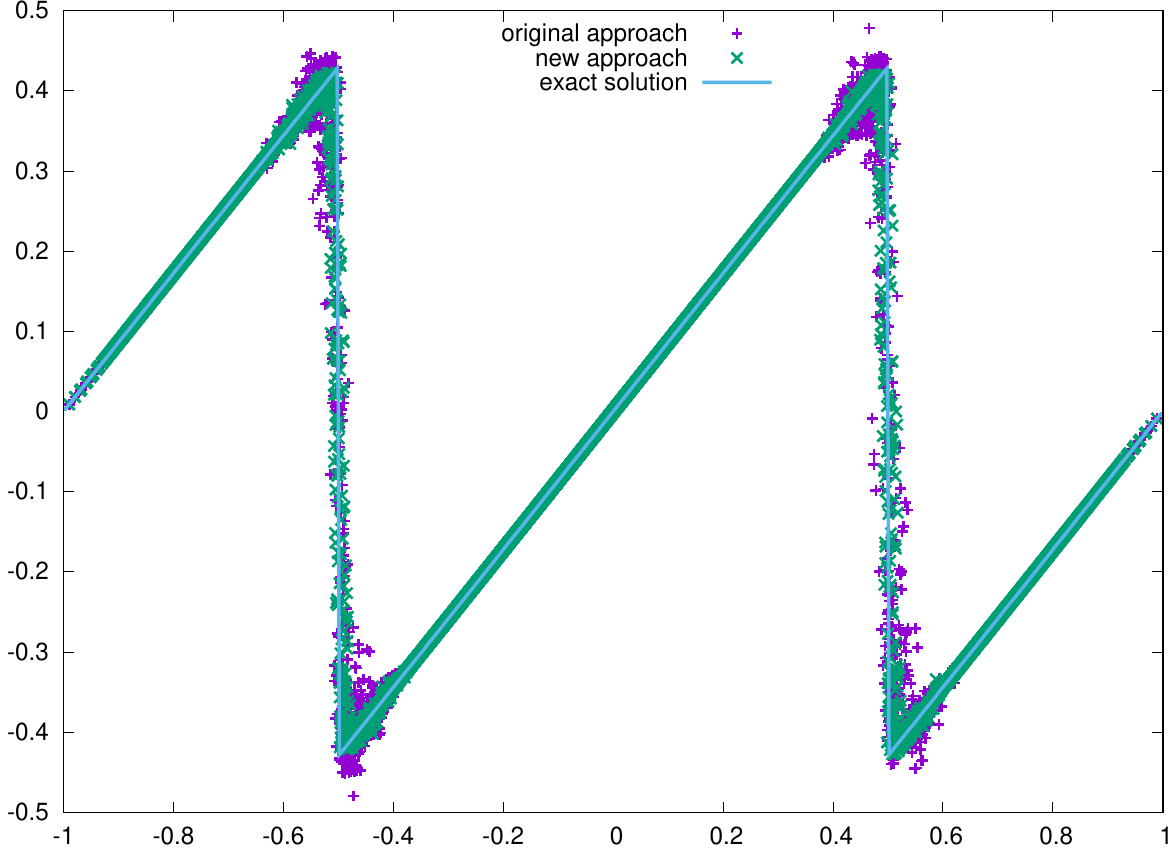}
    \caption{Comparison between original and new correction procedure: submean values versus $(x+y-1)$ coordinate.}
    \label{fig_burgers_compar_profile}
  \end{center}
\end{figure}

Finally, we plot in Figure~\ref{fig_burgers_compar_profile} for both corrections, as well as for the exact entropic solution, the subcell mean values $\ov{u}_m^c$ versus $r_m^c=x_c+y_c-1$, where $(x_c,\,y_c)$ stands for the barycenter of subcell $S_m^c$. These results demonstrates once more how the new correction leads to a better representation of the shocks, with a less oscillatory profile. For that reason, the new correction based on a convex combination of first-order FV numerical flux and high-order DG reconstructed flux, with decreasing weight in the vicinity of a troubled subcell, will be adopted for the numerical applications presented in the next section.

\newpage
\section{Numerical results}
\label{sect_results}

In this numerical results section, we make use of several widely addressed and challenging test cases to demonstrate the performance and robustness of this \textit{a posteriori} local subcell correction for discontinuous Galerkin schemes. In all following test cases, the simple case of global Lax-Friedrichs numerical flux will be used for both the DG scheme and the first-order finite volume reconstructed flux correction. Regarding the cell decomposition into subcells, while this has no impact on the reformulation of DG schemes into subcell finite volume method, see Section~\ref{sect_DG_as_FV}, it may have a slight impact on the results obtained by means of the present correction. In the one-dimensional case, see \cite{vilar_aplsc_1D}, it has been demonstrated that the use of a non-uniform subdivision, for instance by means of the Gauss-Lobatto points, leads to better results compared to a uniform subdivision. In this two-dimensional framework the two types of subdivisions introduced respectively in Figure~\ref{fig_tri1} and \ref{fig_tri3}.\\

Regarding the time integration, we make use of the classical third-order SSP Runge-Kutta scheme, see for instance \cite{Osher}. As the correction described earlier combines both DG scheme on the primal cells $\omega_c$ and FV scheme on the subcells $S_m^c$, the time step is computed adaptively using the following CFL condition
\begin{align}
  \label{cfl}
  \Dt=\Frac{\Min_c \(\Frac{d_c}{2\,k+1},\;\Min_m d_m^{\,c}\)}{\gamma},
\end{align}

where $\gamma=\Max_{c,\, m}\(||\bs{F}'(\ov{u}_m^c)||_2\)$, and where the cell and subcell characteristic lengths $d_c$ and $d_m^{\,c}$ are defined as follows
\begin{align}
  \label{length}
  d_c=\Frac{|\omega_c|}{\Sum_{\omega_v\,\in \,\mc{V}_c}l_{cv}}\qquad\text{ and }\qquad d_m^{\,c}= \Frac{|S_m^c|}{\Sum_{S_{p}^v\,\in \,\mc{V}_m^{\,c}}l_{mp}^c}.
\end{align}

We recall that $l_{mp}^c$ stands for the length of face $f_{mp}^c$ separating subcell $S_m^c$ and its neighbor $S_p^v$, while $l_{cv}$ is the length of the interface between cell $\omega_c$ and its neighbor $\omega_v$. Let us note that in cases where we compute rates of convergence, a time step $\Dt\leq\Min_cd_c^{\frac{k+1}{3}}$ is used in order to make the time error negligible in comparison to the spatial discretization error.\\

Let us emphasize that in all figures to come, if not stated otherwise, the solution subcell mean values are displayed. Consequently, for a mesh made of $N_c$ cells, the numerical solution will be represented on $N_k\, N_c$ subcells, where $N_k=\frac{(k+1) (k+2)} {2}$.

\subsection{Linear case}
\label{subsect_linear}

Let us first assess the performance and accuracy of the present \textit{A Posteriori} Local Subcell Corrected DG (APLSC-DG) scheme in the case of 2D linear conservation laws. In Section~\ref{subsect_detection}, when the NAD criterion based on a discrete maximum principle has been introduced, we did not specified the set $\mc{N}(S_m^c)$ of subcell $S_m^c$ characterizing the DMP. For linear problems, we make use of a cell-wise DMP, meaning $\mc{N}(S_m^c)$ will be constituted by all the subcells of cell $\omega_c$, as well as the subcells of $\omega_v$ the face neighboring cells of $\omega_c$. By means of notations previously introduced, this definition can be rewritten as

\begin{align}
  \label{NAD_set}
  \mc{N}(S_m^c)=\Big\{S_q^v;\;\;\forall\; \omega_v\in\mc{V}_c\cup\{\omega_c\},\; \forall\, q\in\bk{1,\,N_k}\Big\}.
\end{align}

\begin{figure}[!ht]
  \begin{center}
    \subfigure[Structured subdivision.]{\includegraphics[height=7.cm]{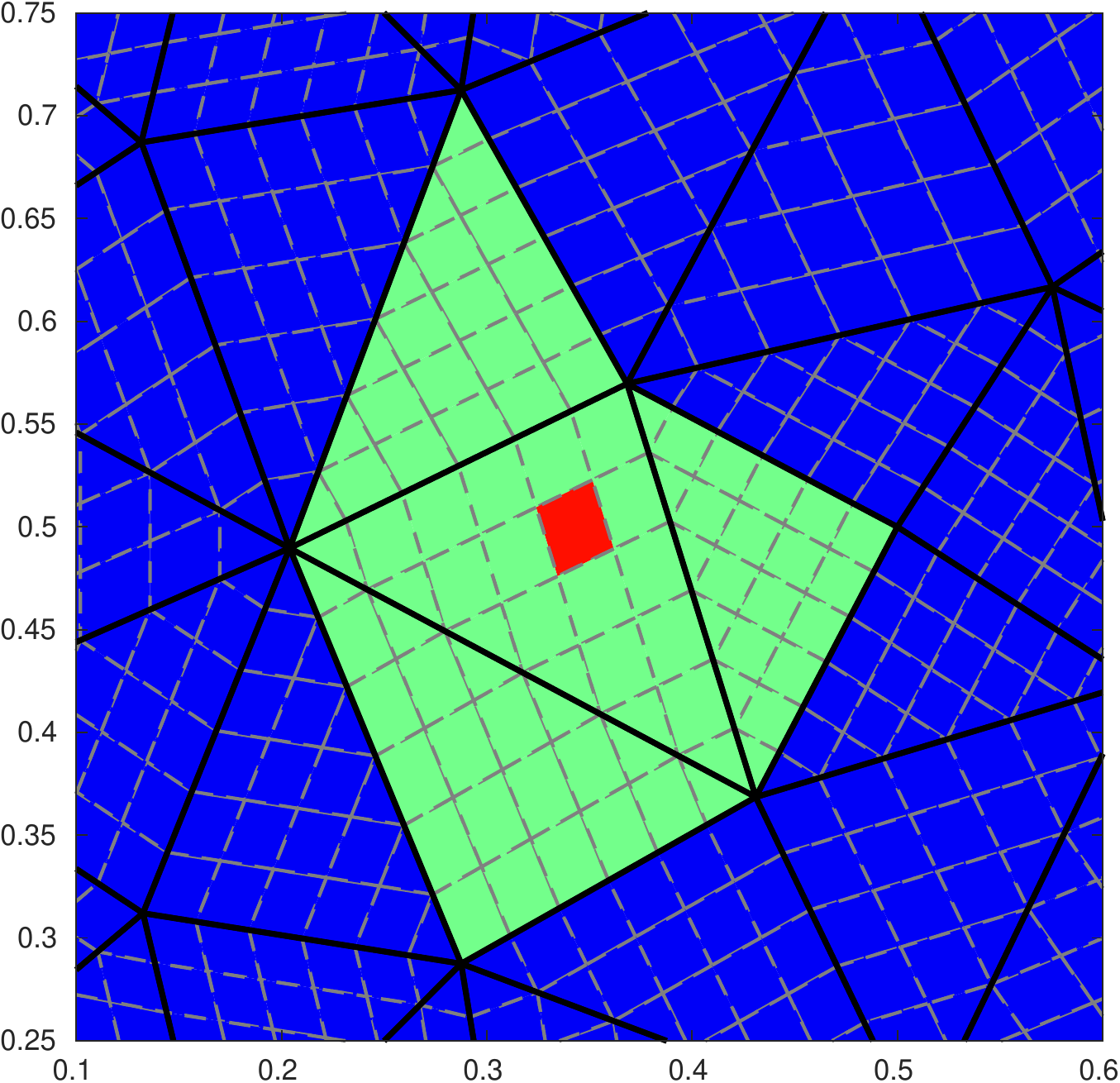}}\hspace*{12.mm}
    \subfigure[Voronoi-type subdivision.]{\includegraphics[height=7.cm]{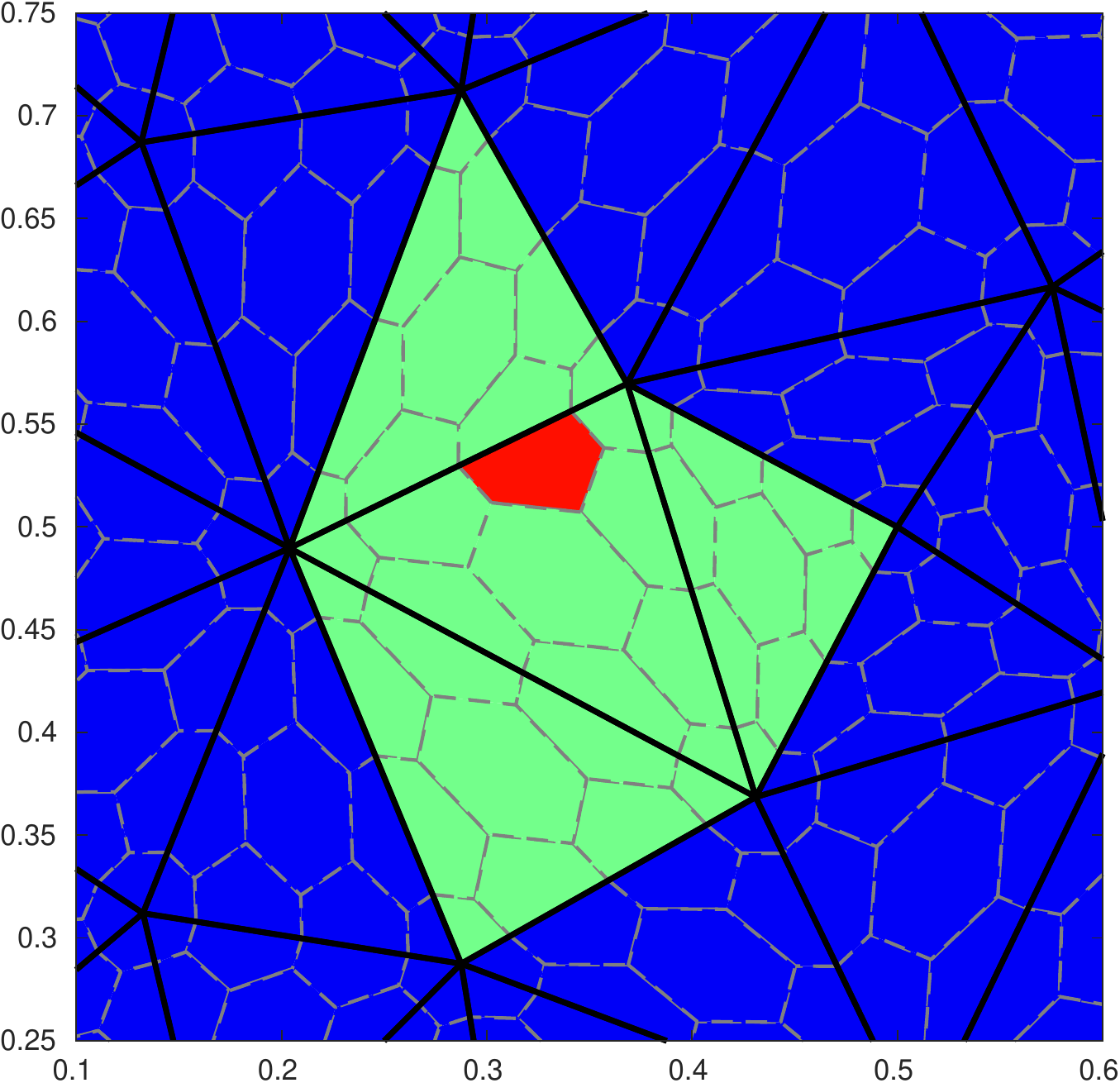}}
    \caption{Neighboring subcells set $\mc{N}(S_m^c)$ for the NAD criterion in the linear case: subcell $S_m^c$ is colored red while the subcells in $\mc{N}(S_m^c)$ are colored green.}
    \label{NAD_linear}
  \end{center}
\end{figure}

This particular set is depicted in Figure~\ref{NAD_linear}, for both the simple structured subdivision as well as the polygonal Voronoi-type one. In this figures, the subcell $S_m^c$ under consideration would be colored red, while the subcells constituting $\mc{N}(S_m^c)$ would be colored green. Let us emphasize that subcell $S_m^c$ is also part of $\mc{N}(S_m^c)$.

\subsubsection{Linear advection}
\label{subsubsect_advect_smooth}

To display the efficiency of DG schemes plus correction, let us first assess how the APLSC behaves in the linear advection case. To this end, we consider the following equation

\begin{subequations}
  \label{eq_advect}
\begin{empheq}[left = \empheqlbrace\,]{align}
&\vdt u(\bs{x},t) + \bs{A}(\bs{x})\pds\gradx{\,u(\bs{x},t)}=0, && (\bs{x},t)\in\,[0,\, 1]^2\times[0,T], \label{adv}\\[2mm]
&u(\bs{x},0)=u_0(\bs{x}), &&\bs{x}\in\,[0,\, 1]^2, \label{adv_ini}
\end{empheq}\\[-4mm]
\end{subequations}

where the transport velocity will be set to $\bs{A}=\(1,\,1\)\tra$.

\vspace*{2mm}
\paragraph{\textbf{Linear advection of a smooth signal}}

\begin{figure}[!ht]
  \begin{center}
    \subfigure[Solution map.]{\includegraphics[height=6.cm]{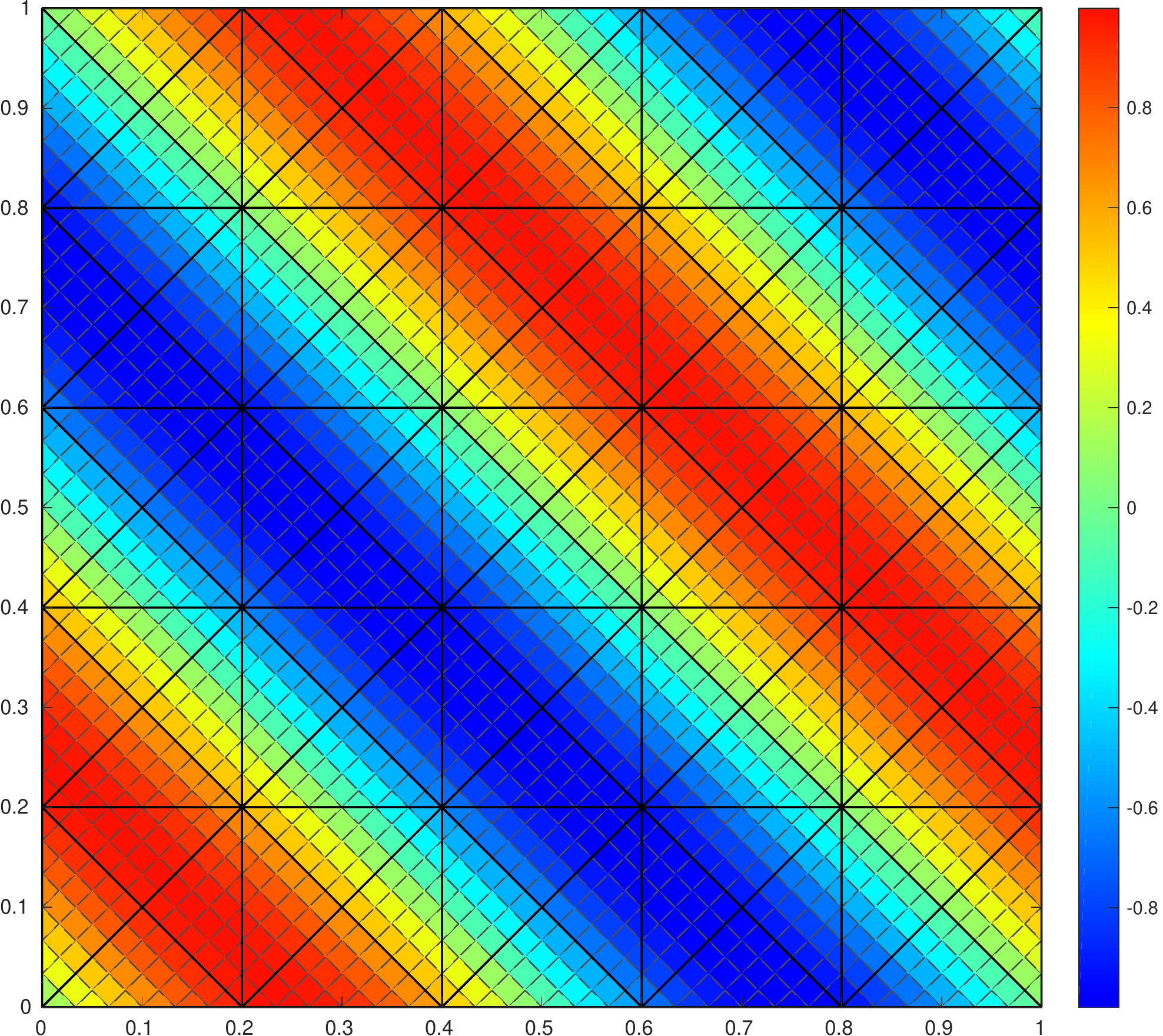}\label{fig_advection_smooth_map}}\hspace*{5.mm}
    \subfigure[Solution profile.]{\includegraphics[height=6.cm]{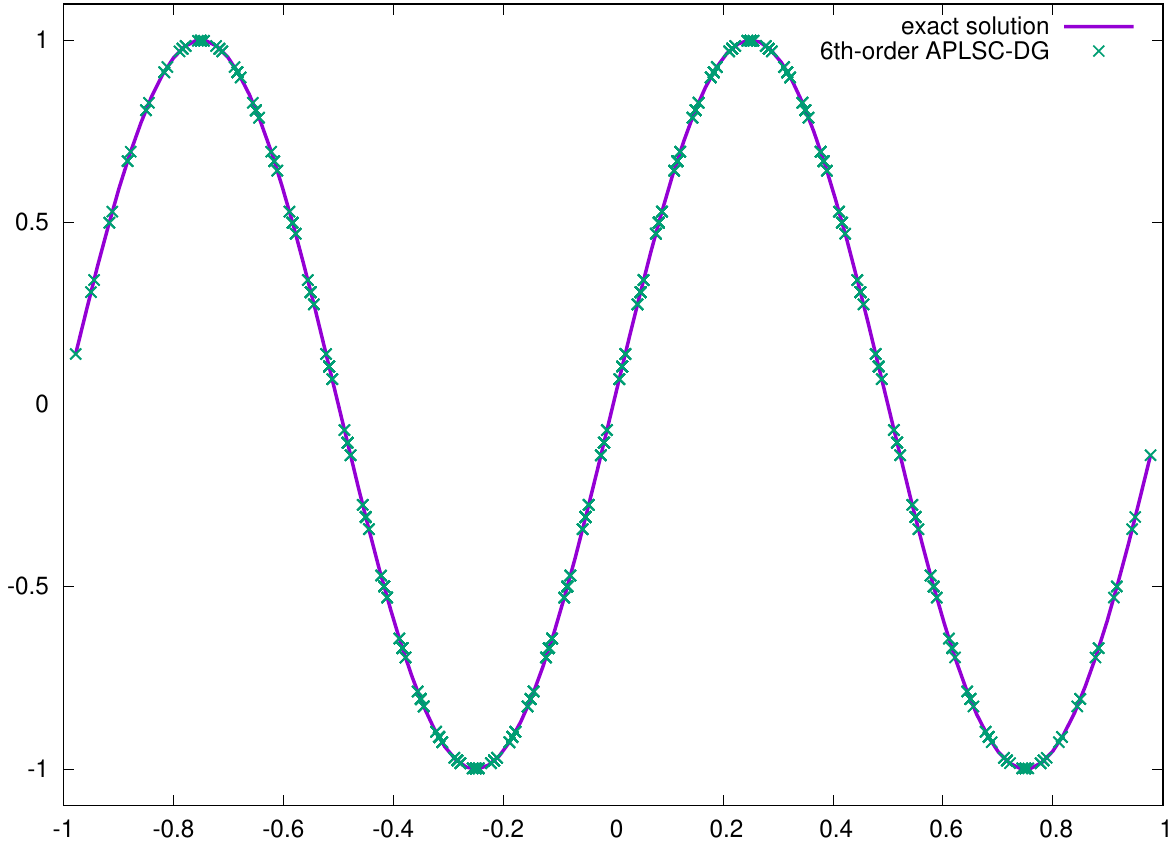}\label{fig_advection_smooth_profile}}
    \caption{Linear advection a smooth signal with a 6th corrected DG scheme, on a 5x5x4 grid after 1 period.}
    \label{fig_advection_smooth}
  \end{center}
\end{figure}

We start from a smooth initial datum $u_0(x,y)=\sin(2 \pi \,(x+y))$, and consider periodic boundary conditions. We assess the scheme accuracy after one period, namely at time $t=1$. In Figure~\ref{fig_advection_smooth}, the numerical solution of the APLSC of sixth-order DG scheme, obtained on grid made of only 100 cells, is depicted. One can see that with only 100 cells, the corrected DG scheme is extremely accurate. Actually, the correction procedure does not activate in this case, which proves that the relaxation criterion on smooth extrema works properly. The rates of convergence are gathered in Table~\ref{table_order_adv_2D} and do exhibit a convergence to six.

\vspace*{10mm}
\begin{table}[!ht]
  \begin{center}
\vspace*{1cm}
    \begin{tabular}{|c||c|c||c|c||c|c|}
      \hline & \multicolumn{2}{c||}{$L_1$} & \multicolumn{2}{c|}{$L_2$} & \multicolumn{2}{c|}{$L_\infty$}\\ \hline\hline $h$ & $E^h_{L_1}$
      & $q^h_{L_1}$ & $E^h_{L_2}$ & $q^h_{L_2}$ & $E^h_{L_\infty}$ & $q^h_{L_\infty}$\\
      \hline $\frac{1}{10}$ & 1.62E-7 & 6.00 & 1.81E-7 & 6.00 & 3.98E-7 & 5.96\\
      \hline $\frac{1}{20}$ & 2.53E-9 & 5.97 & 2.82E-9 & 5.96 & 6.38E-9 & 5.41\\
      \hline $\frac{1}{40}$ & 4.03E-11 & - & 4.52E-11 & - & 1.50E-10 & -\\
      \hline 
    \end{tabular}
  \end{center}
    \caption{Convergence rates for the linear advection case for a 6th-order APLSC-DG scheme}
  \label{table_order_adv_2D}
\end{table}

\paragraph{\textbf{Linear advection of a crenel signal}}

To assess the efficiency of the correction presented in the presence of discontinuity, let us start with the simple case of the advection of a crenel signal, where the initial solution $u_0$ is defined as follows

\begin{align}
  \label{u0_crenel}
  u_0(\bs{x})=\left\{\begin{array}{ll}
                       1 \quad &\text{if }\; (x+y)\in\,[\inv{4}, \,\inv{2}]\cup[\frac{5}{4}, \,\frac{3}{2}],\\[3mm]
                       0 &\text{if }\; (x+y)\in\,[\frac{3}{4}, \,1]\cup[\frac{7}{4}, \,2],\\[3mm]
                       \demi &\text{otherwise}.
                     \end{array}\right.
\end{align}

\begin{figure}[!ht]
  \begin{center}
    \subfigure[Uncorrected DG: min=-0.06, max=1.06.]{\includegraphics[height=6.5cm]{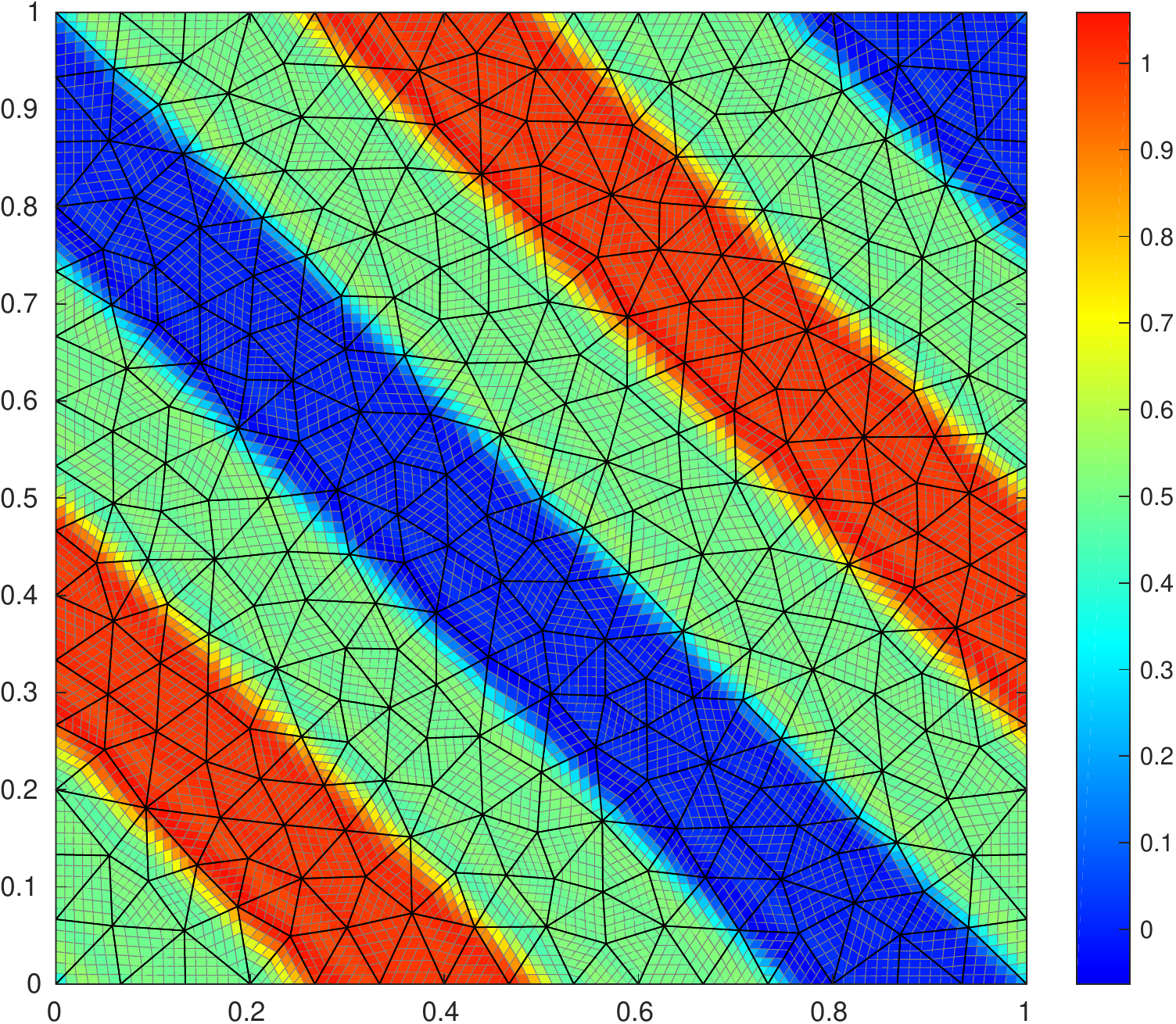}\label{fig_crenel_order6_nocorr}}\hspace*{5mm}
    \subfigure[APLSC-DG: min=8.5E-6, max=1.]{\includegraphics[height=6.5cm]{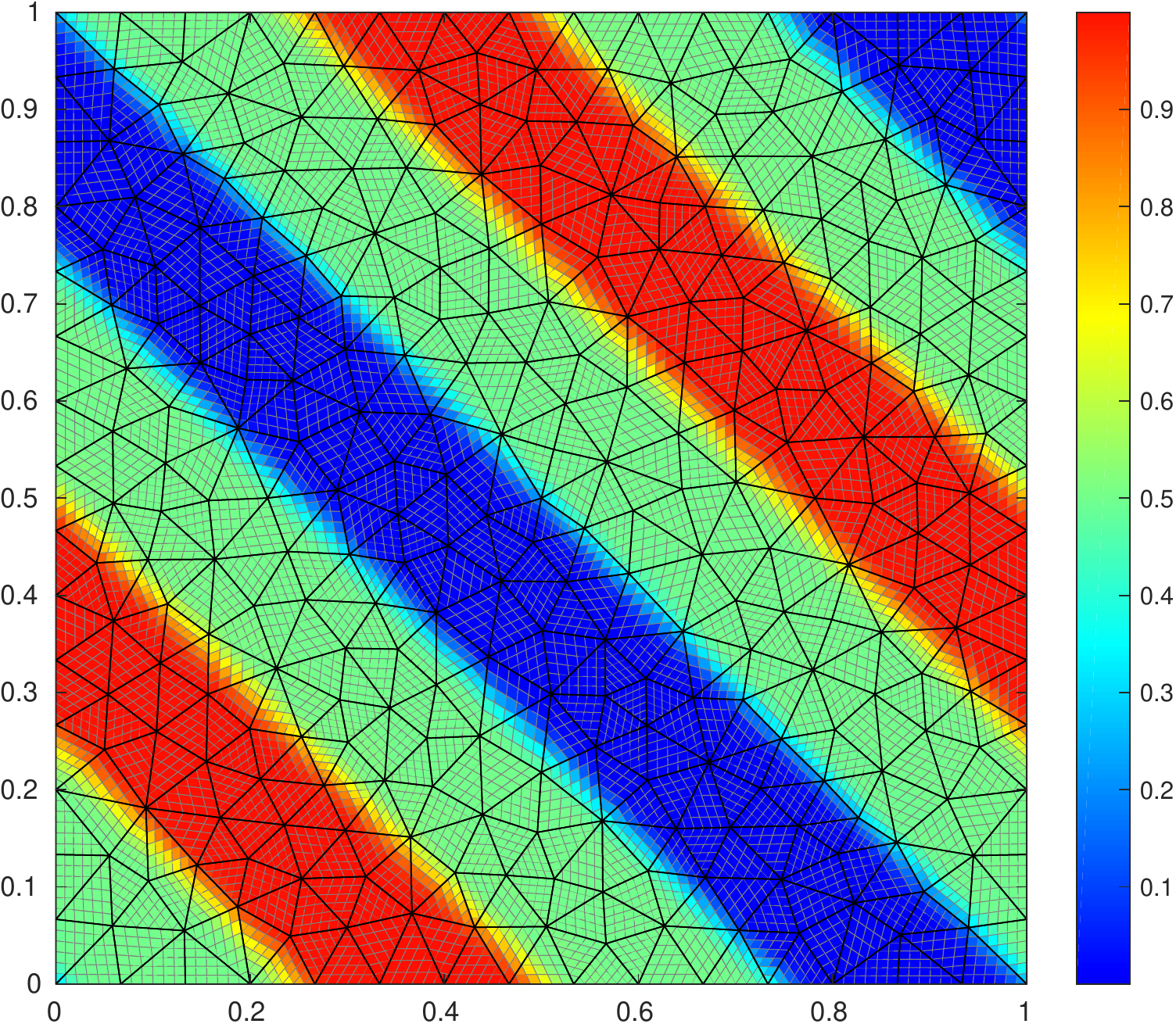}\label{fig_crenel_order6}}
    \caption{6th-order DG solutions for the linear advection of a crenel signal on 576 cells after one period.}
  \label{fig_crenel_6th}
  \end{center}
\end{figure}

In Figure~\ref{fig_crenel_6th}, we compare the uncorrected and corrected versions of the 6th-order DG scheme on an unstructured grid made of 576 cells, after one period. One can see that in both cases, the numerical solution obtained is very accurate, but the one obtained through the APLSC-DG method respects the maximum principle as the final solution remains in the bounds of the initial one.

\begin{figure}[!ht]
  \begin{center}
    \includegraphics[height=7.5cm]{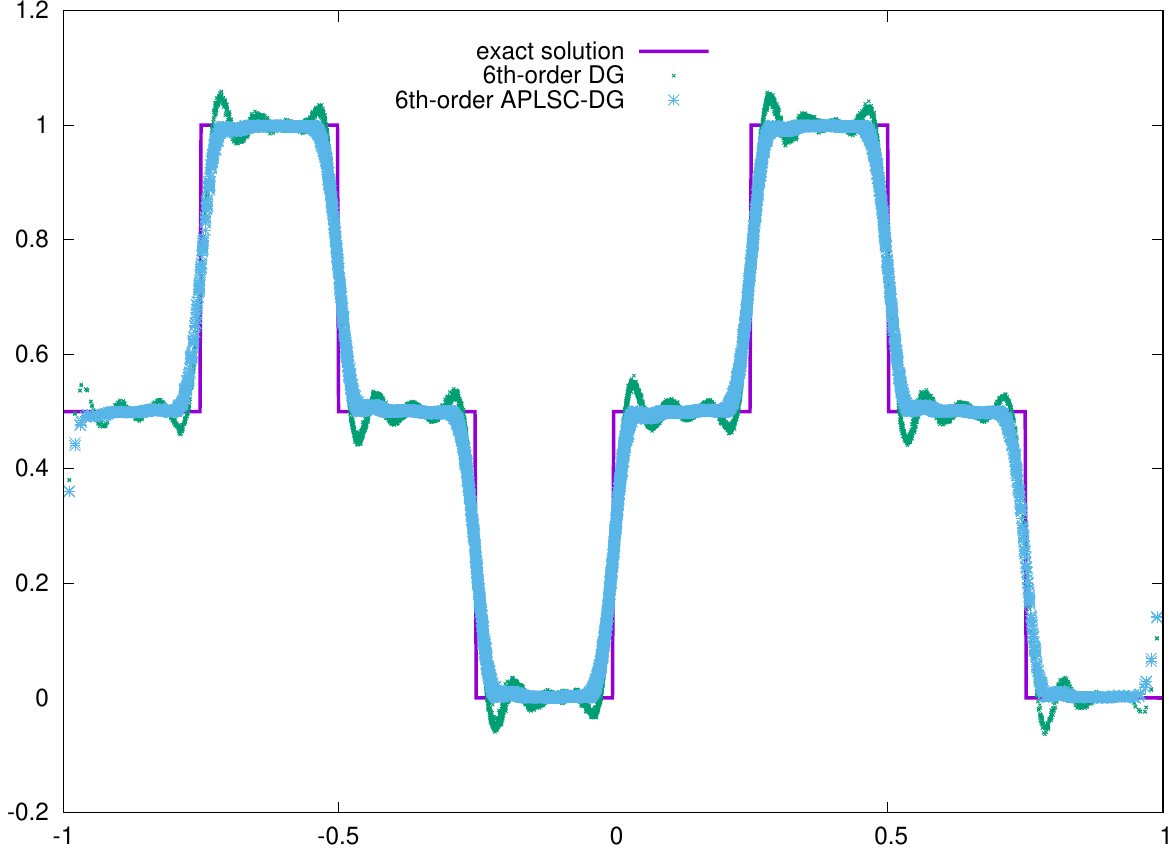}
    \caption{6th-order solutions for the crenel advection case on 576 cells: submean values versus $(x+y-1)$ coordinate.}
    \label{fig_crenel_6th_profile}
  \end{center}
\end{figure}

In Figure~\ref{fig_crenel_6th_profile}, the submean values versus $(x+y-1)$ coordinate are compared for both the uncorrected and the APLSC-DG schemes. We note how the spurious oscillations in the vicinity of the discontinuities have been removed, while preserving the very precise resolution of the 6th-order DG scheme.\\

In the previous examples, the simple uniform structured cell subdivision has been used. Let us emphasize that the cell subdivision does not  theoretically impact the equivalency between DG scheme and subcell FV-like scheme defined through the reconstructed fluxes introduced previously, see Theorems~\ref{main_thm_residual} and \ref{main_thm_flux}. This choice, if it does, should only influence the corrected scheme. To illustrate such statement, let us run the linear advection of the crenel signal on five periods for different types of cell subdivision.\\

\begin{figure}[!ht]
  \begin{center}
    \subfigure[Equidistant boundary points.]{\includegraphics[height=4.5cm]{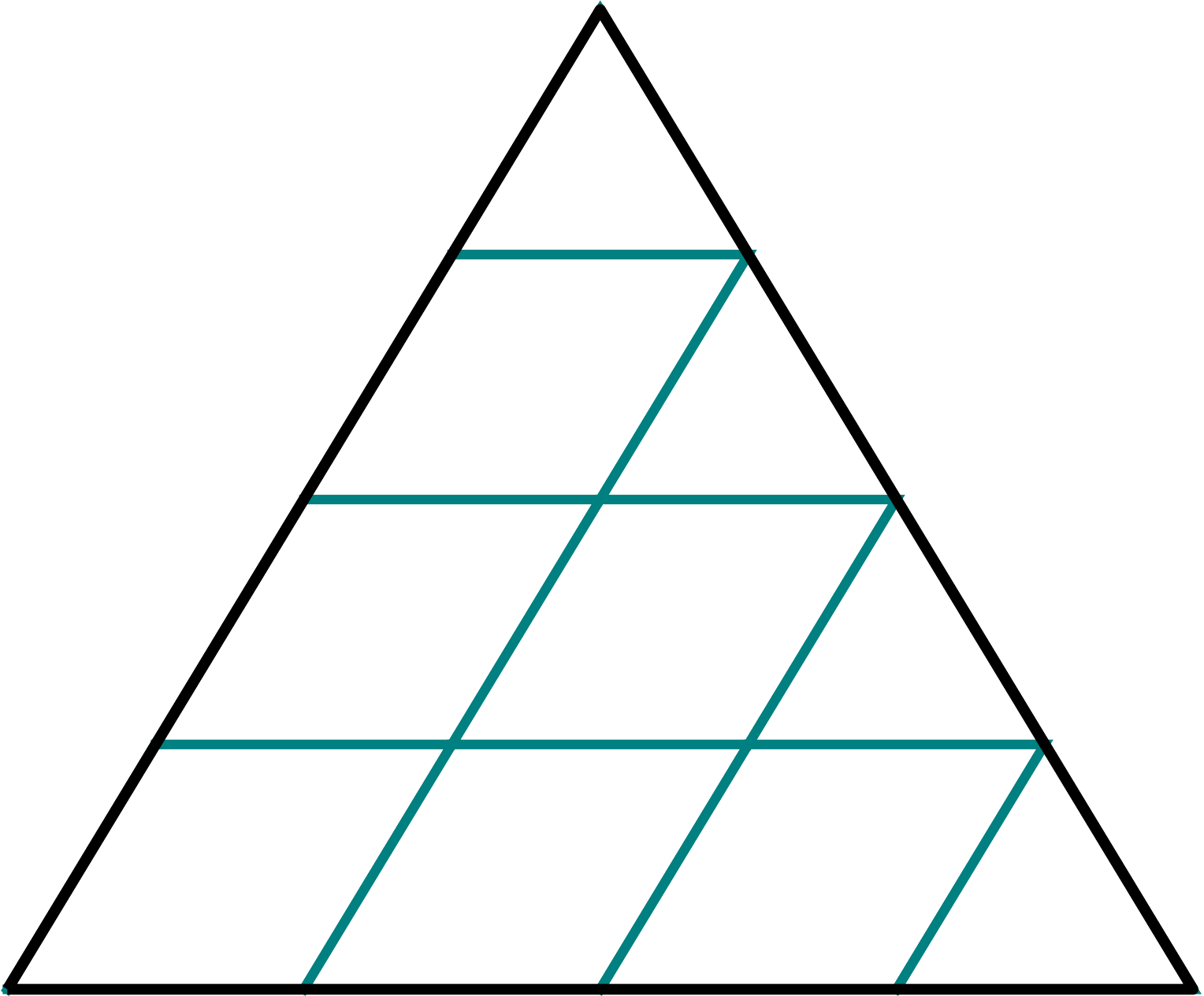}\label{fig_subdvi_cart1}}\hspace*{12mm}
    \subfigure[Gauss-Lobatto boundary points.]{\includegraphics[height=4.5cm]{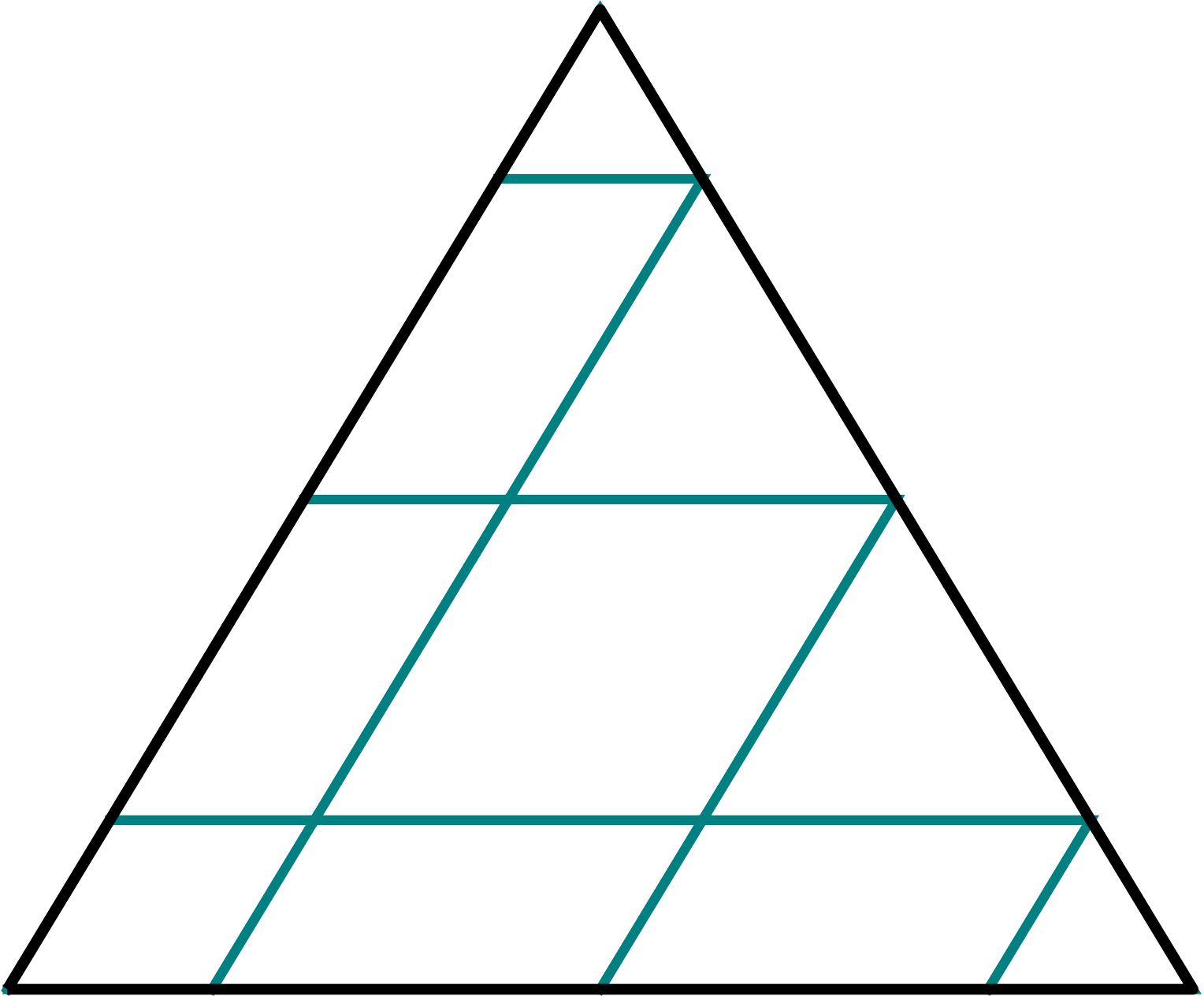}\label{fig_subdvi_cart3}}
    \caption{Examples of structured subdivisions for a triangular cell and a $\P^3$ DG scheme.}
  \label{fig_subdiv_cart}
  \end{center}
\end{figure}

\begin{figure}[!ht]
  \begin{center}
    \subfigure[Equidistant boundary points.]{\includegraphics[height=4.5cm]{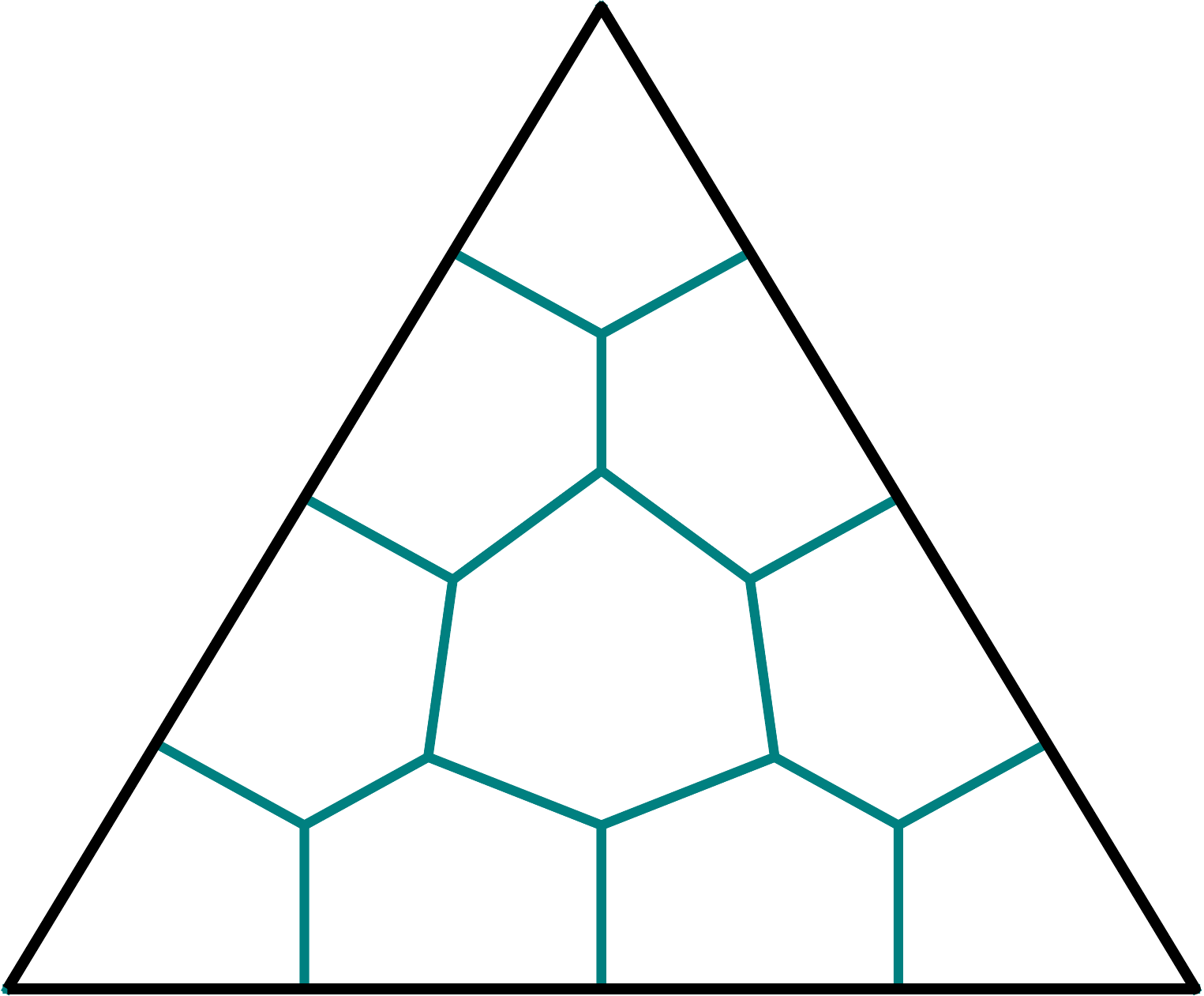}\label{fig_subdvi_poly1}}\hspace*{12mm}
    \subfigure[$\P^3$ Lagrangian mid-points.]{\includegraphics[height=4.5cm]{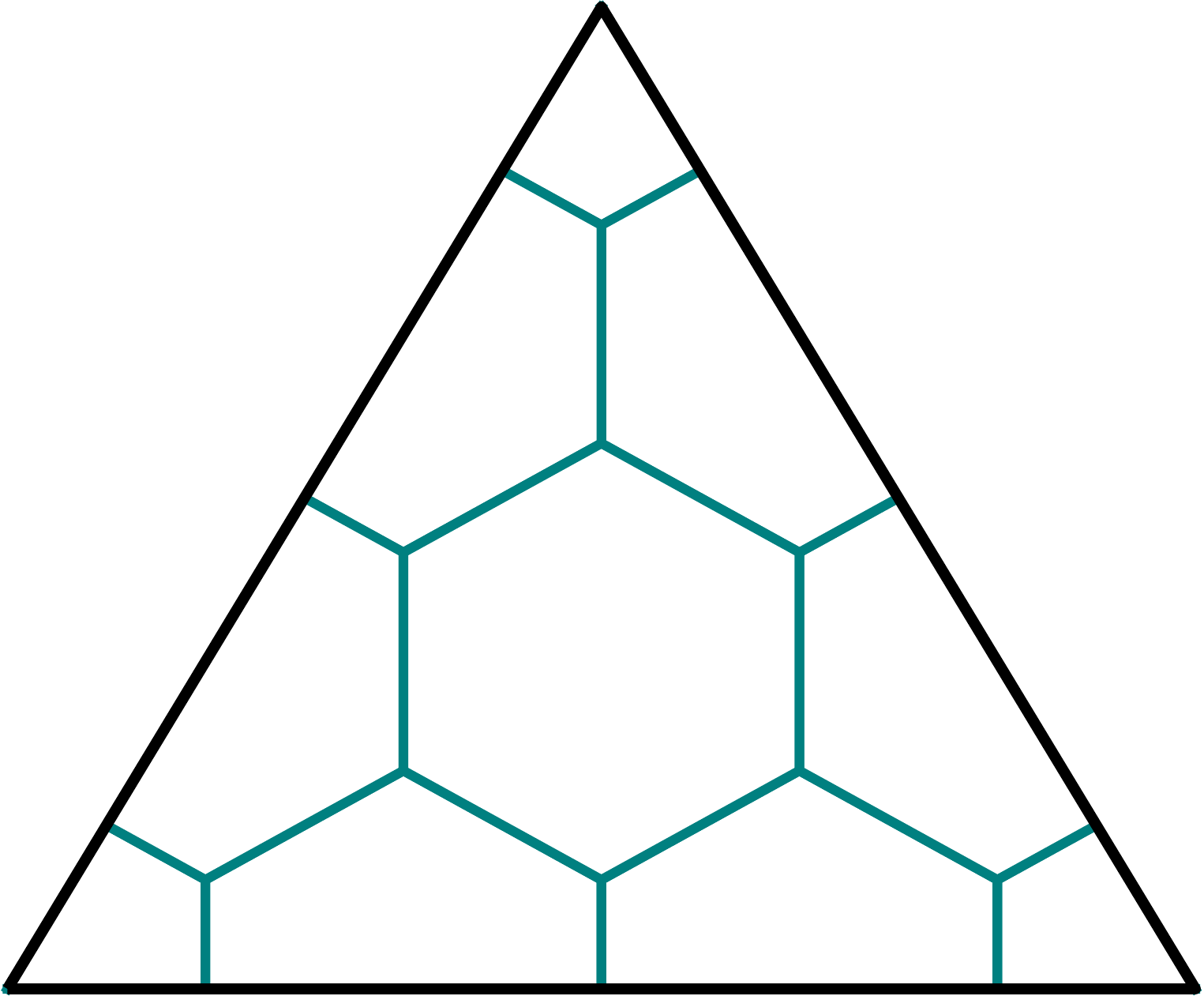}\label{fig_subdvi_poly2}}
    \caption{Examples of Voronoi-type subdivision for a triangular cell and a $\P^3$ DG scheme.}
  \label{fig_subdiv_poly}
  \end{center}
\end{figure}

In Figure~\ref{fig_subdiv_cart}, the simple case of structured subdivision is depicted where the cell is partitioned as a quadrilateral cell would be and then split into two to fit the triangular shape. In this case, one can note that the subdivision is not rotation invariant and that a choice has to be made as one corner subcell is a quadrilateral while the other two are triangles. In this work, we made the choice to start the structured cell subdivision from the wider angle corner, which induces the quadrilateral subcell to stand as this particular corner. In Figure~\ref{fig_subdvi_cart1}, the cell boundary points defining the subdivision are distributed in an uniform manner, while in Figure~\ref{fig_subdvi_cart3} there are defined as Gauss-Lobatto quadrature points. Let us emphasize that in both cases, those subdivisions are very simple to implement and to generalize to any order of accuracy. Furthermore, the subcells normals are nothing but the ones to the original triangular cell, and the subcells are either triangles or parallelograms. In Figure~\ref{fig_subdiv_poly}, polygonal subdivisions are displayed. Those Voronoi type subdivisions are widely used in subcell techniques and can be found for instance in \cite{Chen2006PartitionsOA,peraire_2012} and references within. In Figure~\ref{fig_subdvi_poly2}, the subdivision is obtained as follows: first, we define a triangulation of the element by joining the $\P^k$ Lagrangian nodes. Then, the centroid of each triangle and the midpoint of the edges form a set of $N_k$ polygonal subcells. In Figure~\ref{fig_subdvi_poly1}, the previous procedure is modified to yield a more uniform subdivision by setting equidistant cell boundary points.\\

\begin{figure}[!ht]
  \begin{center}
    \subfigure[Uniform structured subdivision.]{\includegraphics[height=6.5cm]{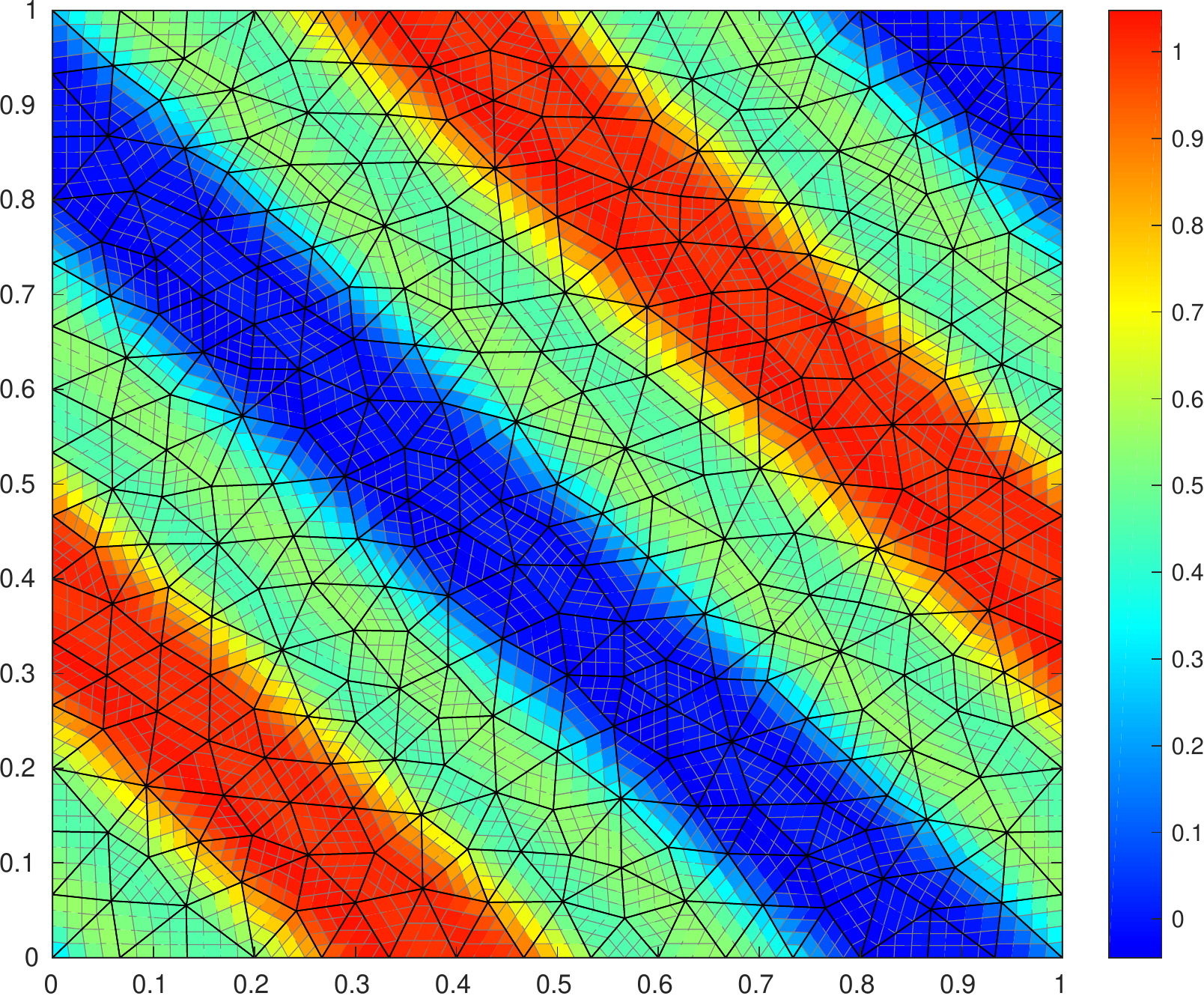}\label{fig_crenel_order4_cart1_nocorr}}\hspace*{5mm}
    \subfigure[Non-uniform structured subdivision.]{\includegraphics[height=6.5cm]{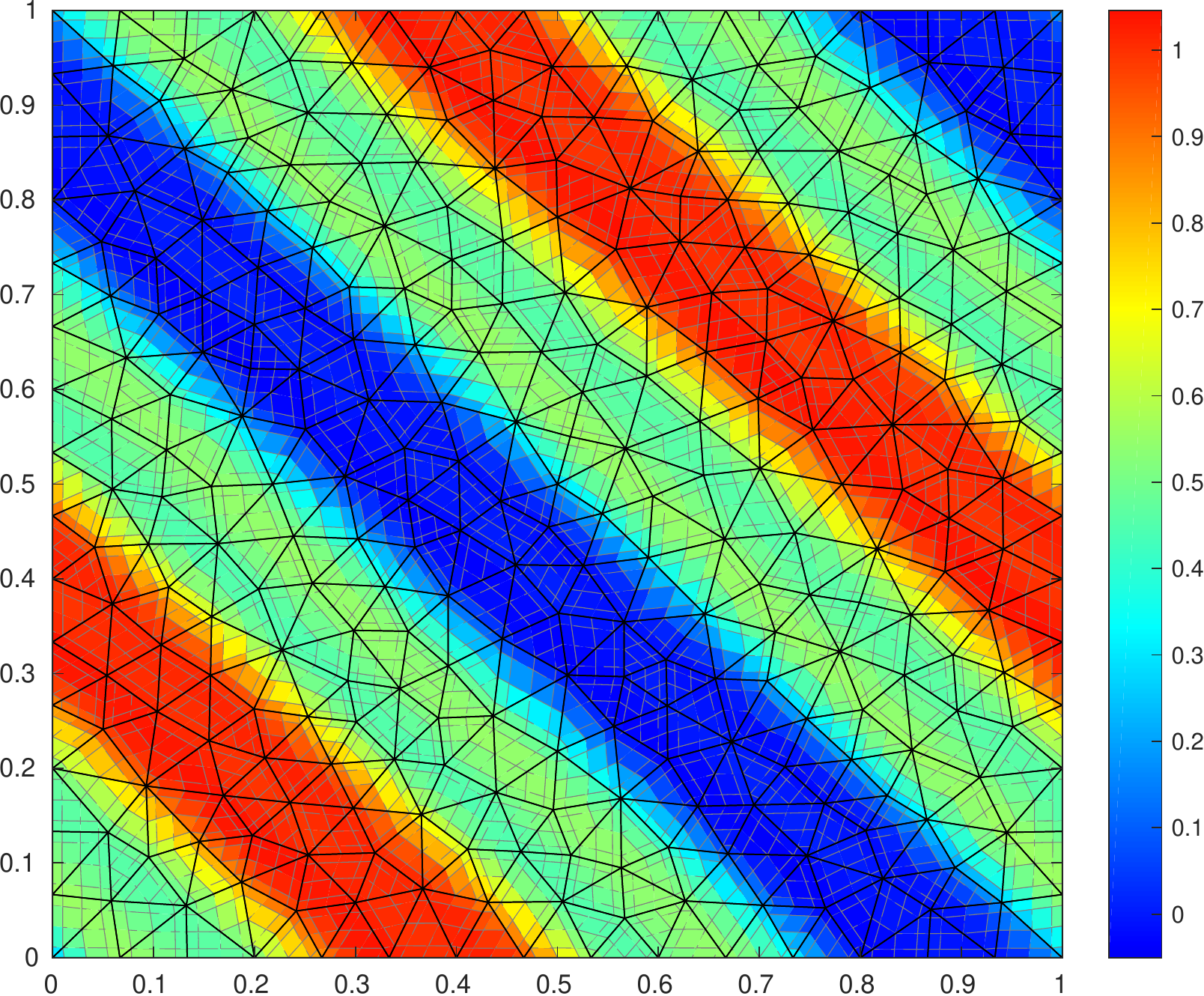}\label{fig_crenel_order4_cart3_nocorr}}
    \caption{4th-order DG solutions for the crenel signal advection on 576 cells after five periods: structured subdivision.}
  \label{fig_crenel_4th_cart_nocorr}
  \end{center}
\end{figure}

\begin{figure}[!ht]
  \begin{center}
    \subfigure[Uniform polygonal subdivision.]{\includegraphics[height=6.5cm]{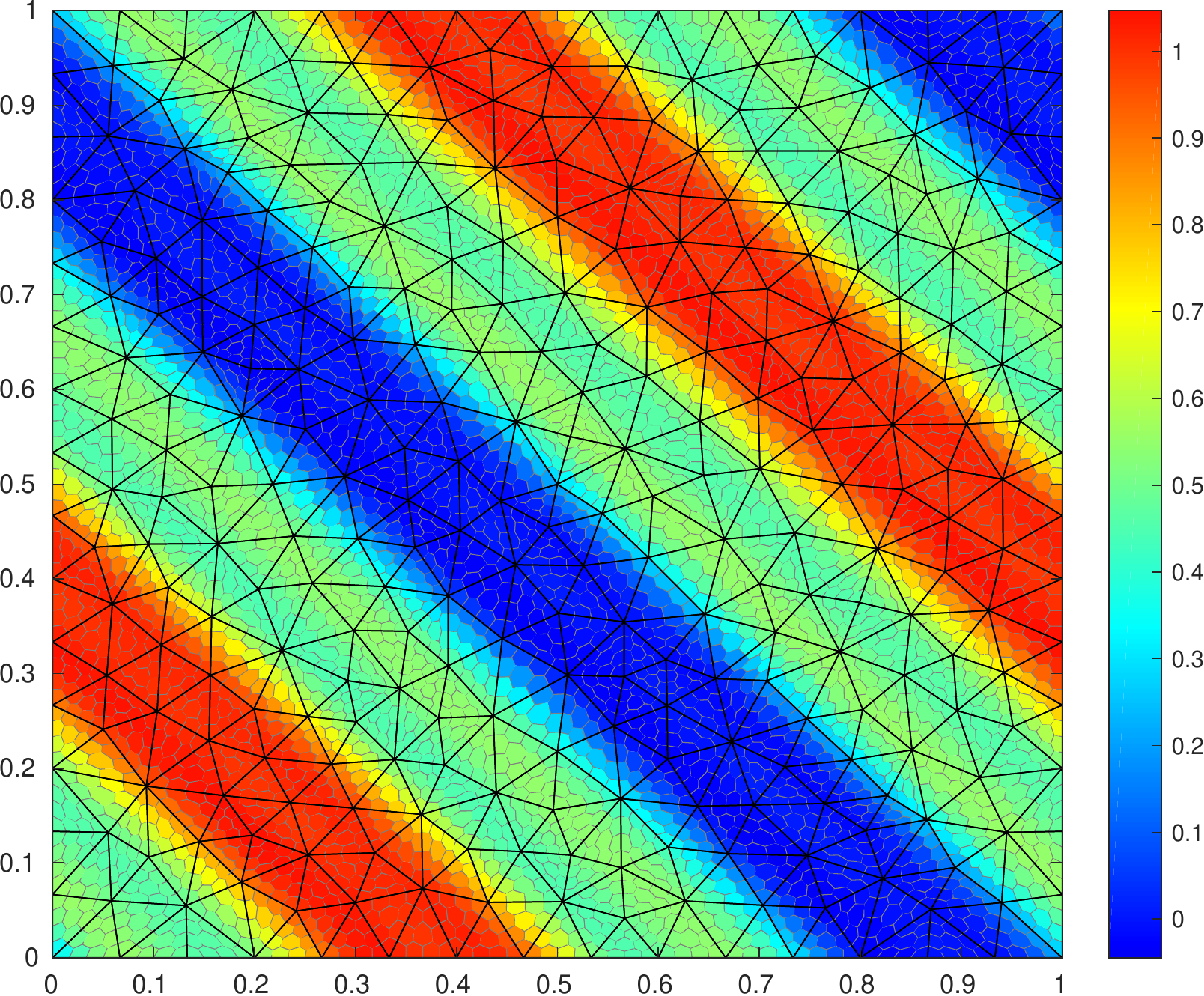}\label{fig_crenel_order4_poly1_nocorr}}\hspace*{5mm}
    \subfigure[Non-uniform polygonal subdivision.]{\includegraphics[height=6.5cm]{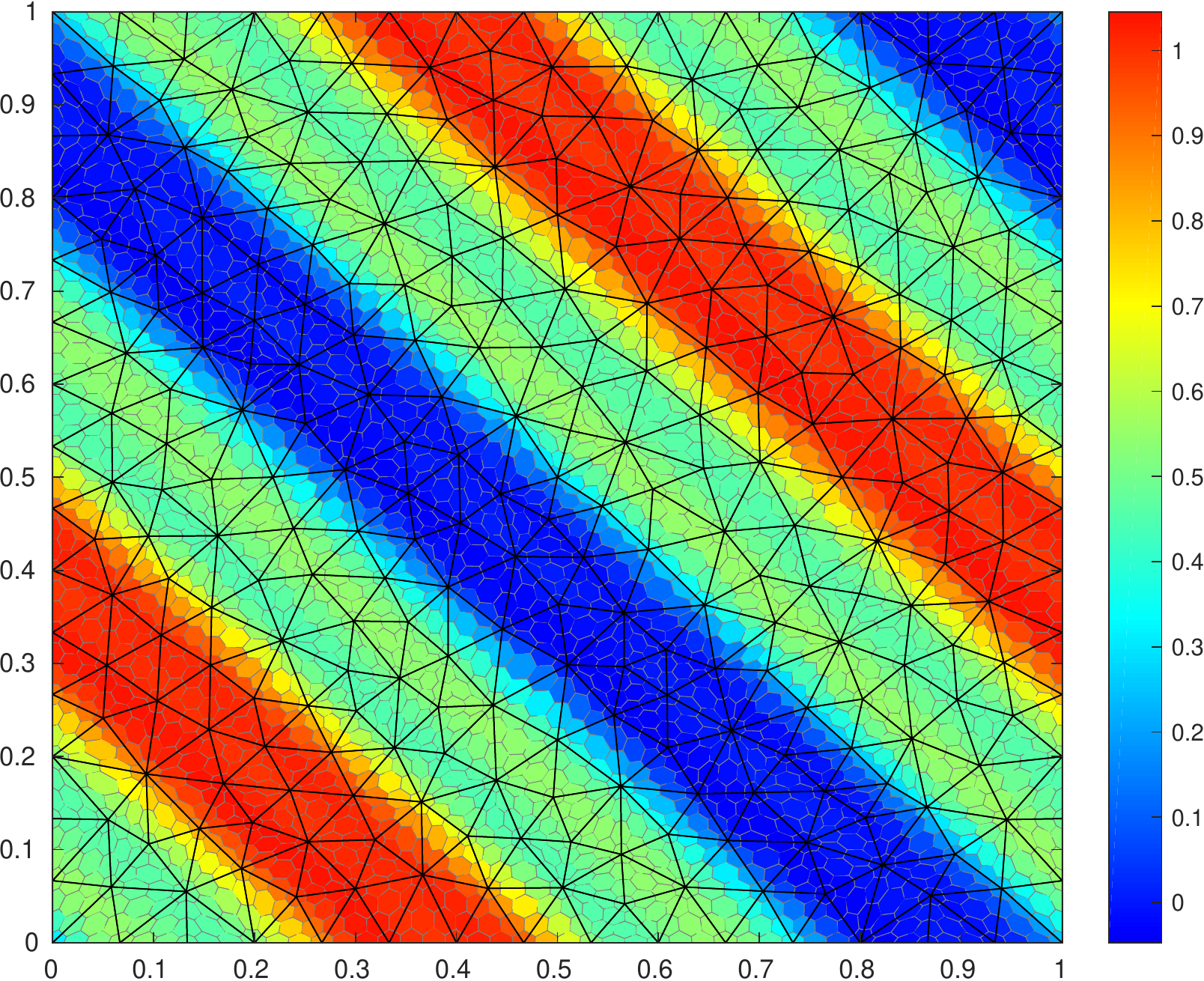}\label{fig_crenel_order4_poly2_nocorr}}
    \caption{4th-order DG solutions for the crenel signal advection on 576 cells after five periods: polygonal subdivision.}
  \label{fig_crenel_4th_poly_nocorr}
  \end{center}
\end{figure}

\begin{figure}[!ht]
  \begin{center}
    \includegraphics[height=7.5cm]{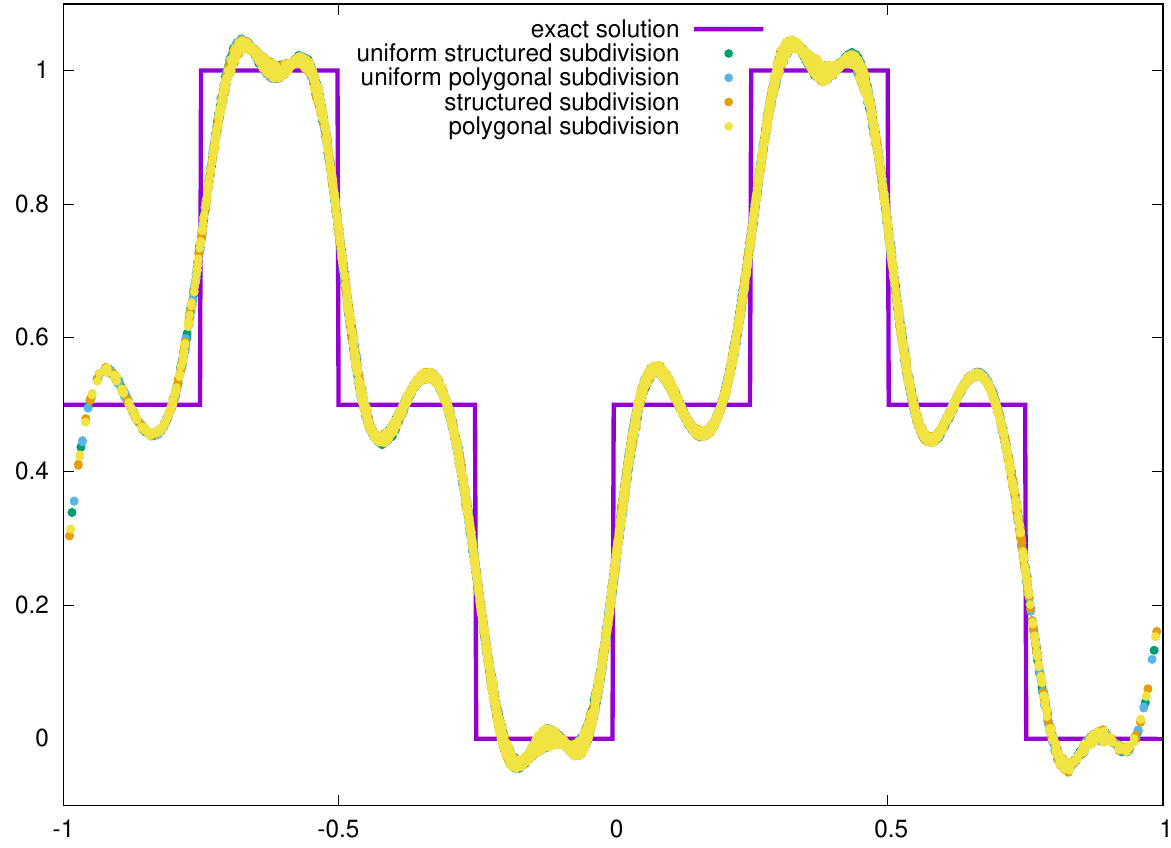}
    \caption{4th-order DG solutions for the crenel signal advection on 576 cells using different cell subdivisions: submean values versus $(x+y-1)$ coordinate.}
    \label{fig_crenel_4th_nocorr_profile}
  \end{center}
\end{figure}

Now, to make sure that cell subdivision does not have any impact on the numerical solution, we display in Figures~\ref{fig_crenel_4th_cart_nocorr}, \ref{fig_crenel_4th_poly_nocorr} and \ref{fig_crenel_4th_nocorr_profile} the solutions obtained through subcell finite volume scheme with high-order reconstructed fluxes as numerical fluxes, equation \eqref{sub_FV_thm_flux}, with the four different cell subdivisions previously introduced. 
In Figure~\ref{fig_crenel_4th_nocorr_profile}, where submean values versus $(x+y-1)$ coordinate are displayed, it is clear that the four calculations lead to exactly the same numerical solution.\\

That being said, the type of subdivision does have an impact on the condition number of the projection matrix, which may lead for very high-order of accuracy to greater numerical errors. In Table~\ref{table_cond_number}, the projection matrix condition numbers for the four different subdivisions displayed before and for different orders of approximations are gathered.

\begin{table}[!ht]
  \begin{center}
    \begin{tabular}{|c||c|c|c|c|c|}
      \hline & $\P^0$ & $\P^1$ & $\P^2$ & $\P^3$\\
      \hline\hline Unif. struct. subdiv. & 1 & 4 & 10.91 & 31.75\\
      \hline Non-unif. struct. subdiv. & 1 & 4 & 9.52 & 29.28\\
      \hline Unif. polyg. subdiv. & 1 & 2.87 & 8.73 & 27.89\\
      \hline Non-unif. polyg. subdiv. & 1 & 2.87 & 8.19 & 26.94\\
      \hline 
    \end{tabular}
  \end{center}
    \caption{Projection matrix condition number for different orders and subdivisions.}
  \label{table_cond_number}
\end{table}

These different condition numbers have been computed through the $||.||_\infty$ consistent norm, such that $||P||_\infty=\max_{\,i}\big(\sum_{\,j} P_{ij}\big)$. It appears that, even if the condition number are closed for these subdivisions, the structured ones do produce higher condition number than the polygonal ones, as well as the uniform ones compared to the non-uniform ones. One can see that higher orders of approximation yield greater condition number. This may ultimately lead to potential issue in the case of extremely high-orders, but in all the cases studied nothing symptomatic has been observed.\\

\begin{figure}[!ht]
  \begin{center}
    \subfigure[Uniform structured subdivision.]{\includegraphics[height=6.5cm]{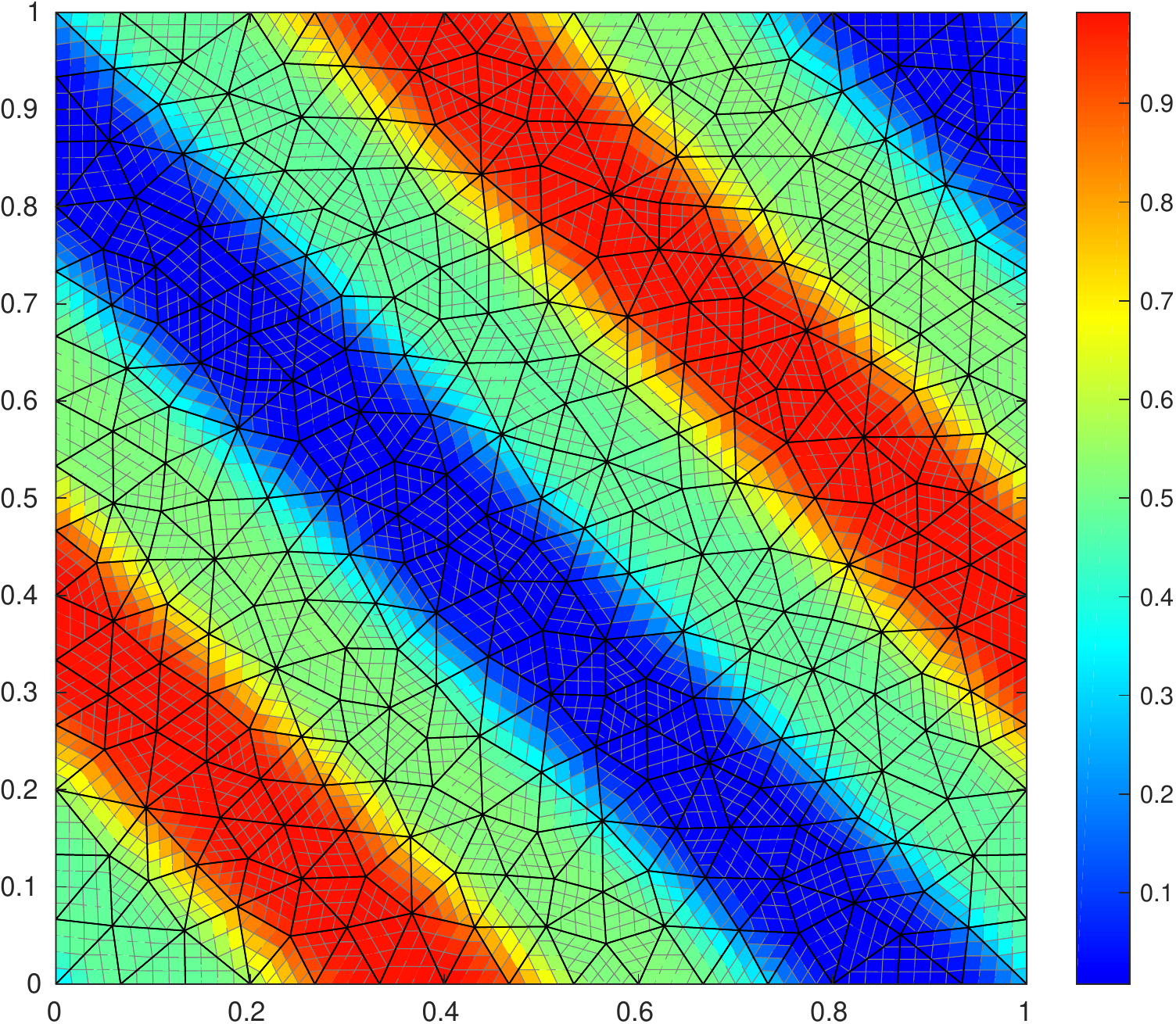}\label{fig_crenel_order4_cart1}}\hspace*{5mm}
    \subfigure[Non-uniform structured subdivision.]{\includegraphics[height=6.5cm]{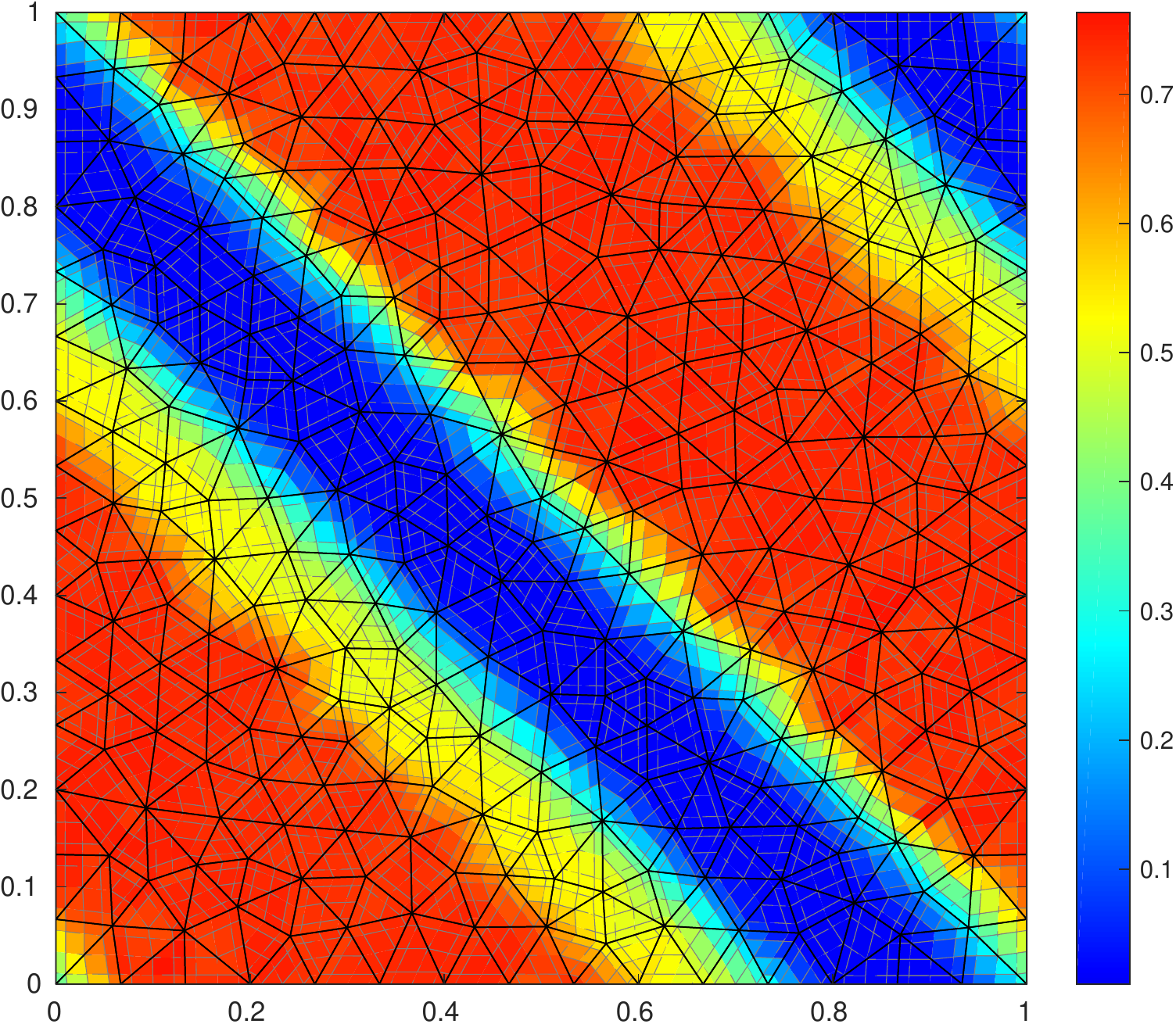}\label{fig_crenel_order4_cart3}}
    \caption{4th-order APLSC-DG solutions for crenel advection on 576 cells after five periods: structured subdivision.}
  \label{fig_crenel_4th_cart}
  \end{center}
\end{figure}

\begin{figure}[!ht]
  \begin{center}
    \subfigure[Uniform polygonal subdivision.]{\includegraphics[height=6.5cm]{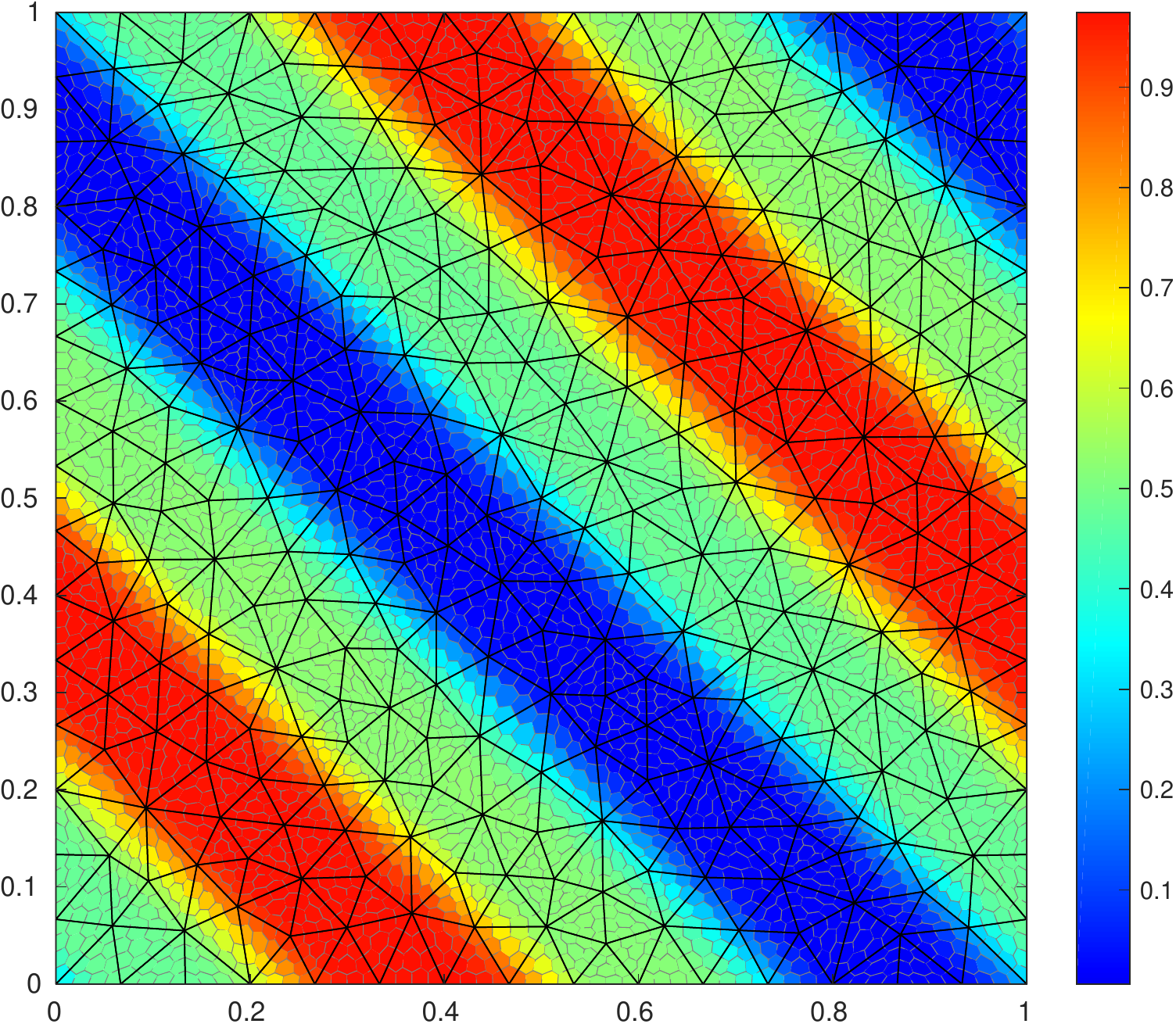}\label{fig_crenel_order4_poly1}}\hspace*{5mm}
    \subfigure[Non-uniform polygonal subdivision.]{\includegraphics[height=6.5cm]{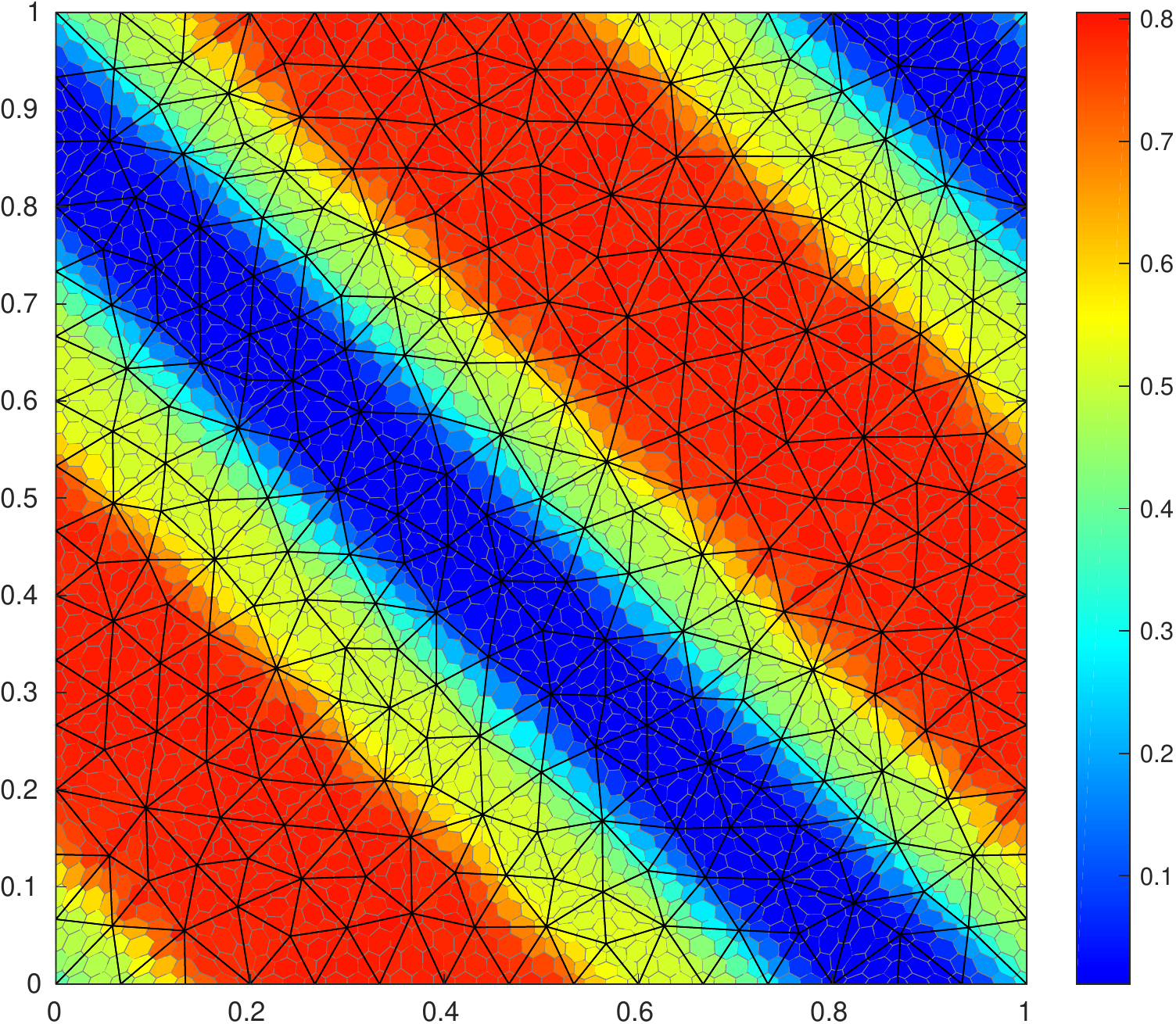}\label{fig_crenel_order4_poly2}}
    \caption{4th-order APLSC-DG solutions for crenel advection on 576 cells after five periods: polygonal subdivision.}
  \label{fig_crenel_4th_poly}
  \end{center}
\end{figure}

\begin{figure}[!ht]
  \begin{center}
    \includegraphics[height=7.5cm]{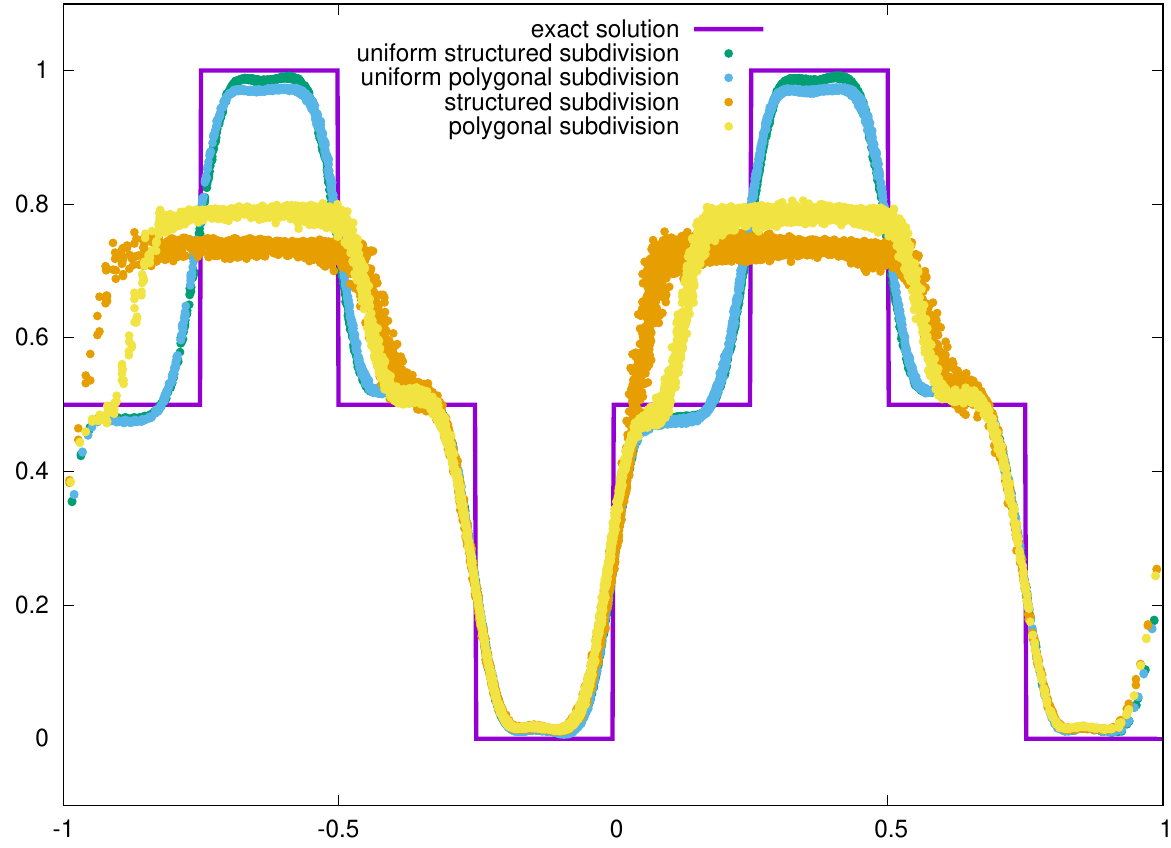}
    \caption{4th-order APLSC-DG solutions for the crenel advection case on 576 cells using different cell subdivisions: submean values versus $(x+y-1)$ coordinate.}
    \label{fig_crenel_4th_profile}
  \end{center}
\end{figure}

While the cell subdivision does not affect the uncorrected DG numerical solution, see Figure~\ref{fig_crenel_4th_nocorr_profile}, it has been proved in the 1D case \cite{vilar_aplsc_1D} that it does have an impact on the quality of the results when the correction is used, especially in the linear advection case. To assess if any similar phenomenon exists in the 2D unstructured case, we use the same setup as before but this time with the full APLSC-DG method, see Figures~\ref{fig_crenel_4th_cart}, \ref{fig_crenel_4th_poly} and \ref{fig_crenel_4th_profile}.
In the light of those results, it appears that uniform subdivisions lead to a way better resolution of the problem under consideration. This huge difference mainly derives from the application of the NAD troubled subcell detector, as if we make just use of the PAD criterion to only enforce the global maximum principle, the difference in the results is a lot more slight. NAD criterion based on a discrete maximum principle is here highly dependent of the subcell aspect ratio. Similar conclusion can be obtained considering other type of problem, even if the difference in the quality of the results would be a lot less significative.

\vspace*{5mm}
\subsubsection{Solid body rotation}
\label{subsubsect_rotation}

\begin{figure}[!ht]
  \begin{center}
    \subfigure[$t=0$: min=0, max=1.]{\includegraphics[height=6.5cm]{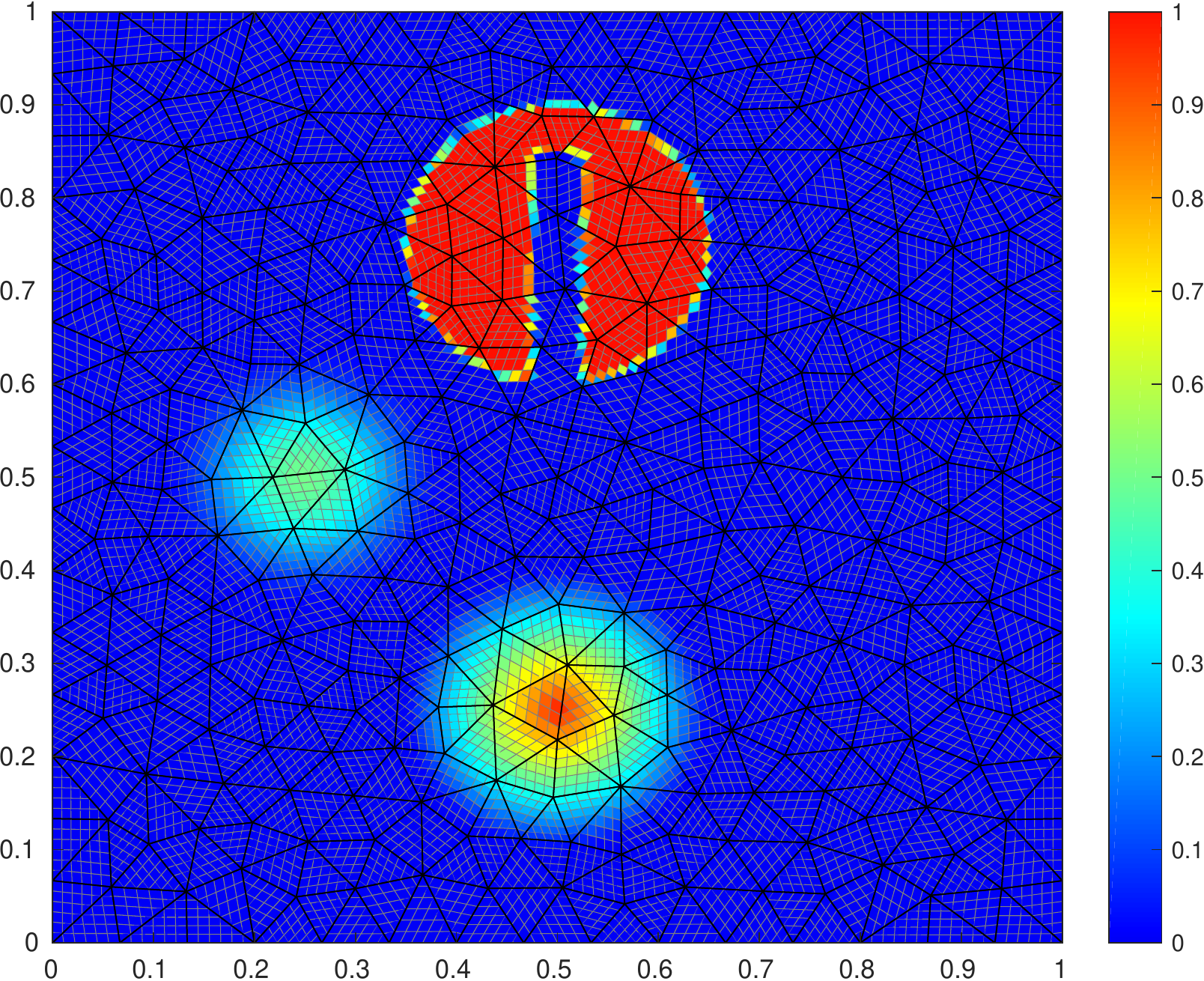}\label{fig_rotate_6th_ini}}\hspace*{4mm}
    \subfigure[$t=2\,\pi$: min=6.7E-10, max=0.97.]{\includegraphics[height=6.5cm]{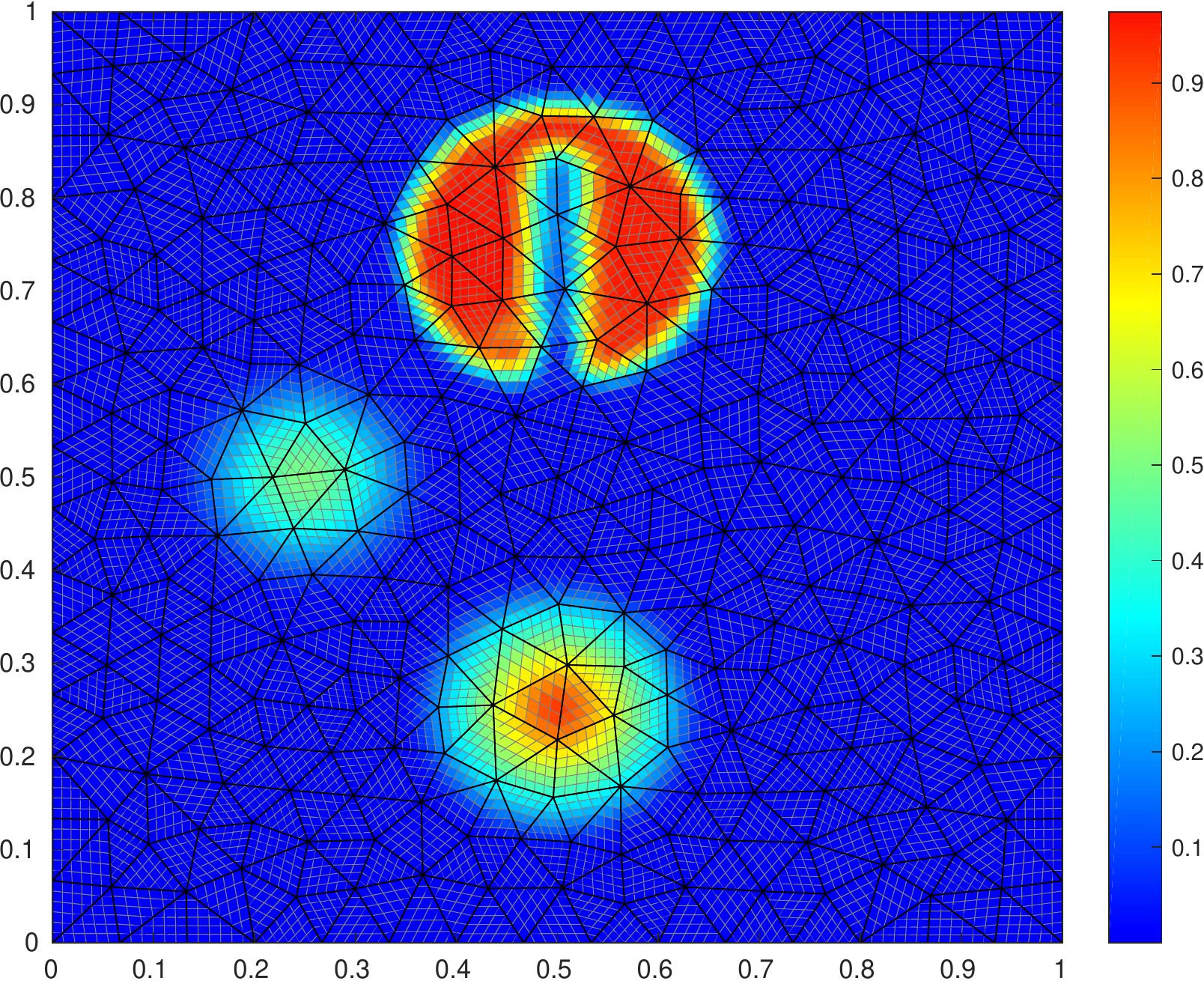}\label{fig_rotate_6th_1T}}
    \caption{6th-order APLSC-DG solution for the rigid rotation case on 576 cells.}
  \label{fig_rotate}
  \end{center}
\end{figure}

We make use of the classical test case taken from \cite{LeVeque2}. Let us then consider \eqref{eq_advect} with a divergence-free velocity field corresponding to a rigid rotation, defined by $\bs{A}(\bs{x})=(\demi-y,\, x-\demi)\tra$. We apply this solid body rotation to the initial data displayed in Figure~\ref{fig_rotate_6th_ini}, which includes both a plotted disk, a cone and a smooth hump. 
In Figure~\ref{fig_rotate}, the submean values obtained through the 6th-order APLSC-DG scheme on a 576 cells coarse grid are displayed. The uniform structured cell subdivision has been used. One can see on Figure~\ref{fig_rotate_6th_1T} how the corrected DG scheme produces a very accurate solution, even using a quite coarse mesh, while still ensuring a global maximum principle as well as a mainly non-oscillatory behavior.

\begin{figure}[!ht]
  \begin{center}
    \subfigure[Solution profile for $y=0.25$.]{\includegraphics[height=5.8cm]{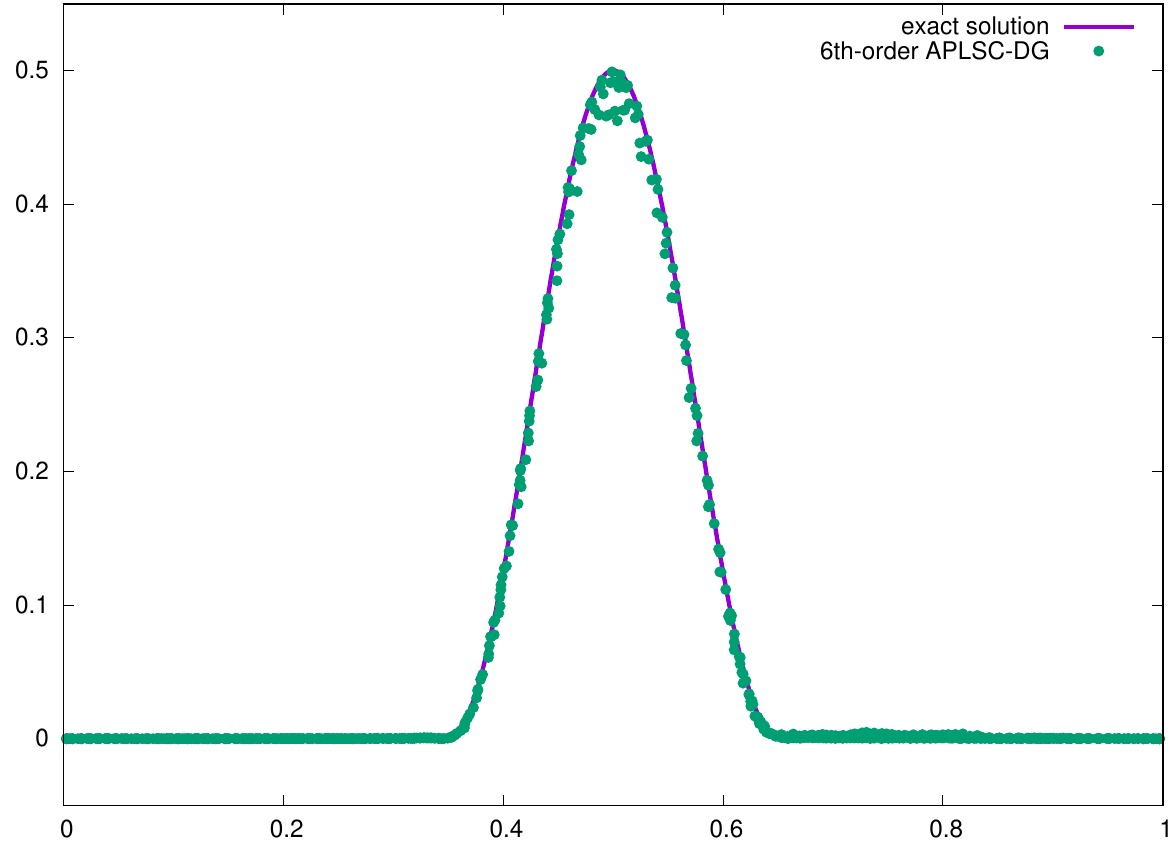}\label{fig_rotate_6th_profile_a}}\hspace*{4mm}
    \subfigure[Solution profile for $x=0.25$.]{\includegraphics[height=5.8cm]{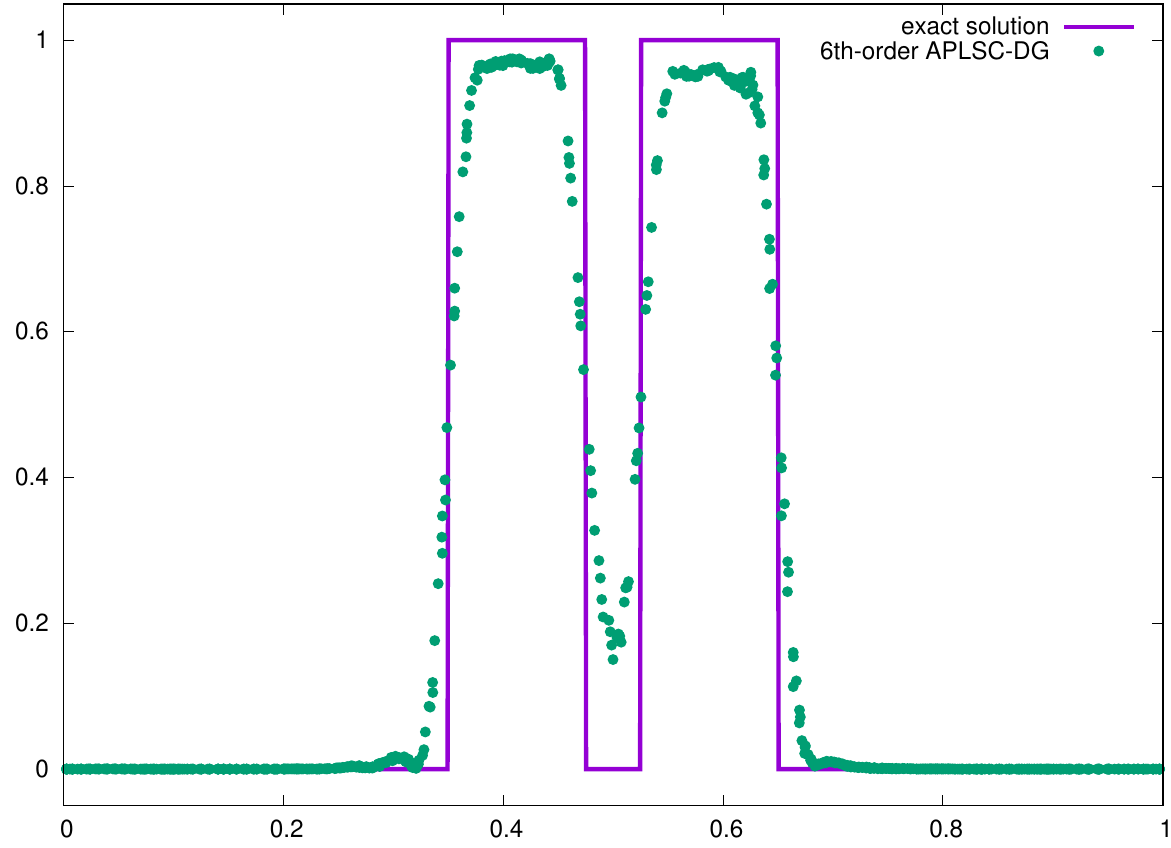}\label{fig_adv_rotate_6th_profile_c}}

    \caption{6th-order APLSC-DG solution for rigid rotation on 576 cells after one full rotation: solution profiles.}
  \label{fig_rotate_6th_profile}
  \end{center}
\end{figure}

In Figure~\ref{fig_rotate_6th_profile}, cross-sections of the solution along lines $y=0.25$  
and $x=0.25$ have been plotted. Those results further demonstrates the very high capability of the correction procedure presented.\\

Now, similarly to the linear advection equation, we want to understand how the choice in the cell subdivision impacts the quality of the results now considering a solid body rotation problem. To do so, we compare the solutions obtained by means of the 4th-order APLSC-DG scheme and the four different subdivisions introduced before, see Figures~\ref{fig_subdiv_cart} and \ref{fig_subdiv_poly}, after five full rotations.
\begin{figure}[!ht]
  \begin{center}
    \subfigure[Uniform structured subdivision.]{\includegraphics[height=6.5cm]{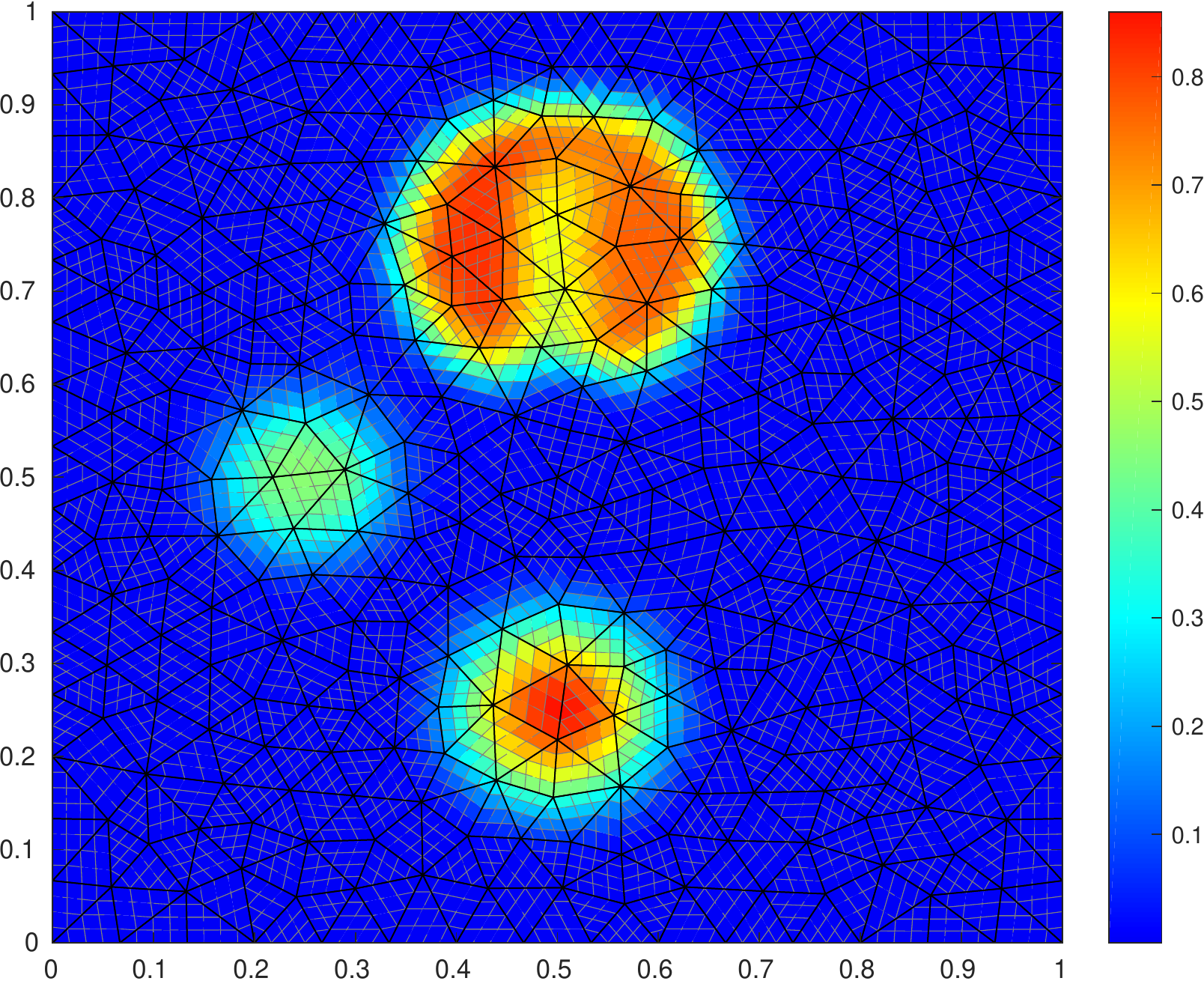}\label{fig_rotate_order4_cart1}}\hspace*{5mm}
    \subfigure[Non-uniform structured subdivision.]{\includegraphics[height=6.5cm]{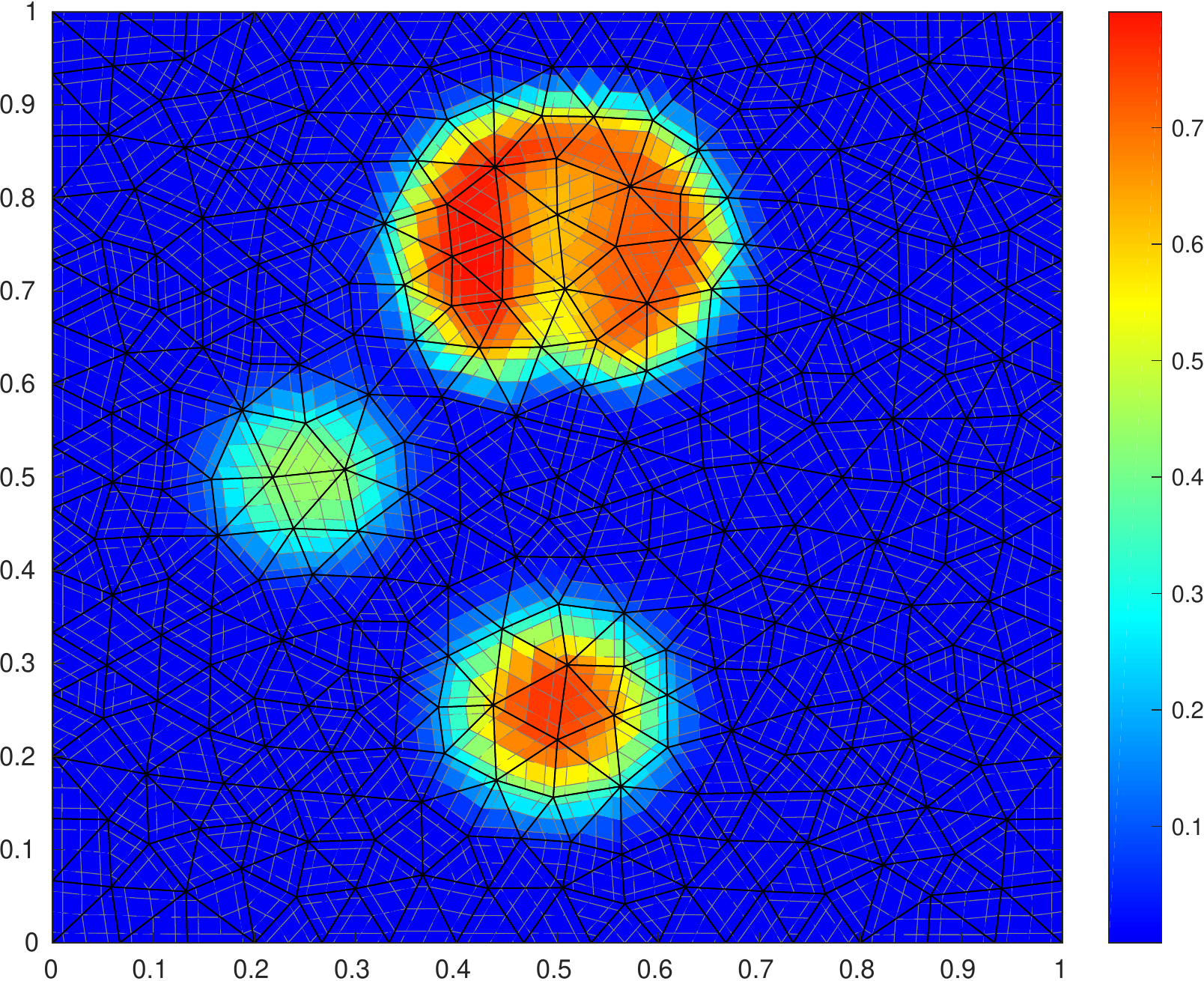}\label{fig_rotate_order4_cart3}}
    \caption{4th-order APLSC-DG solutions for rigid rotation on 576 cells after five full rotations: structured subdivision.}
  \label{fig_rotate_4th_cart}
  \end{center}
\end{figure}
\begin{figure}[!ht]
  \begin{center}
    \subfigure[Uniform polygonal subdivision.]{\includegraphics[height=6.5cm]{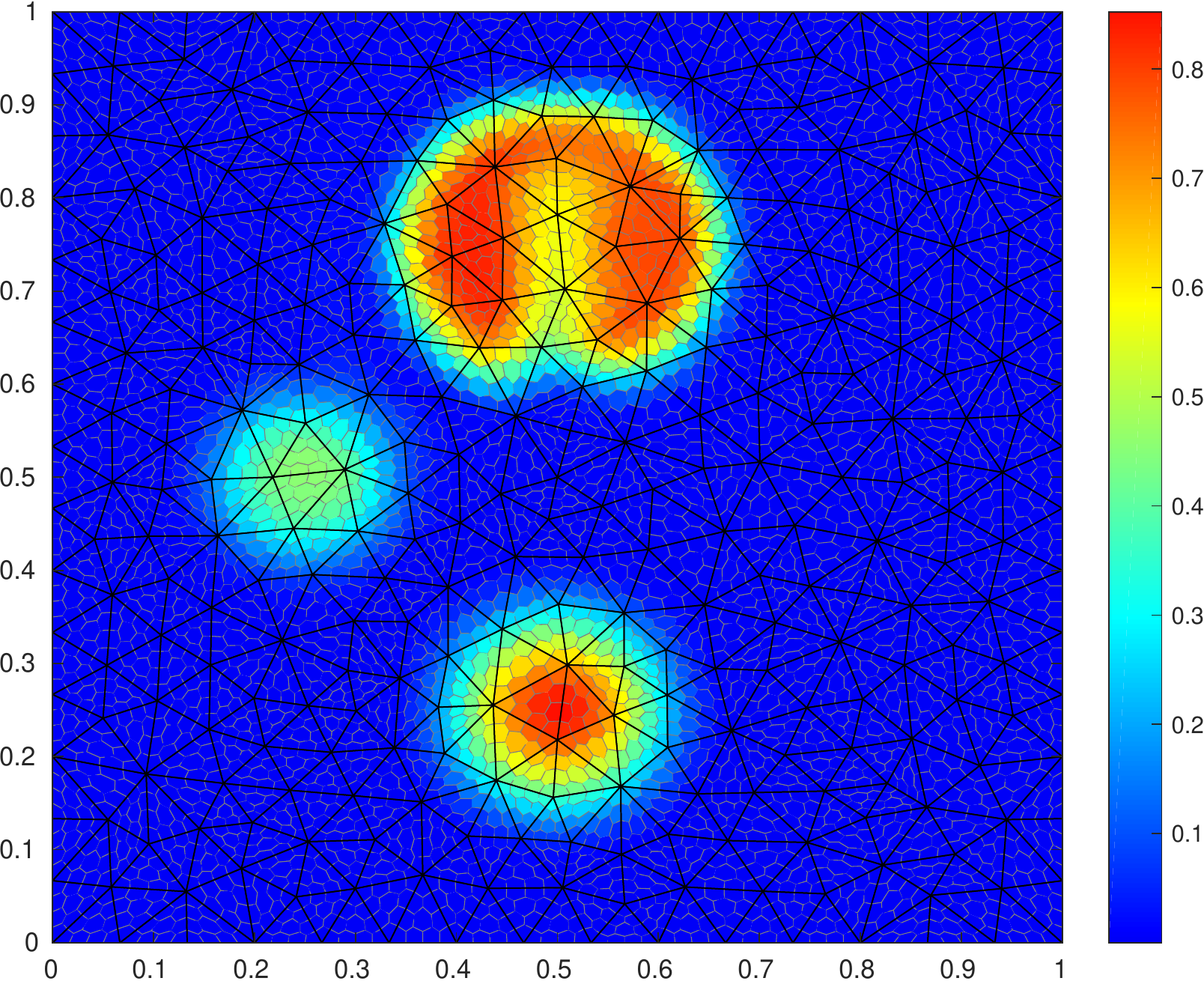}\label{fig_rotate_order4_poly1}}\hspace*{5mm}
    \subfigure[Non-uniform polygonal subdivision.]{\includegraphics[height=6.5cm]{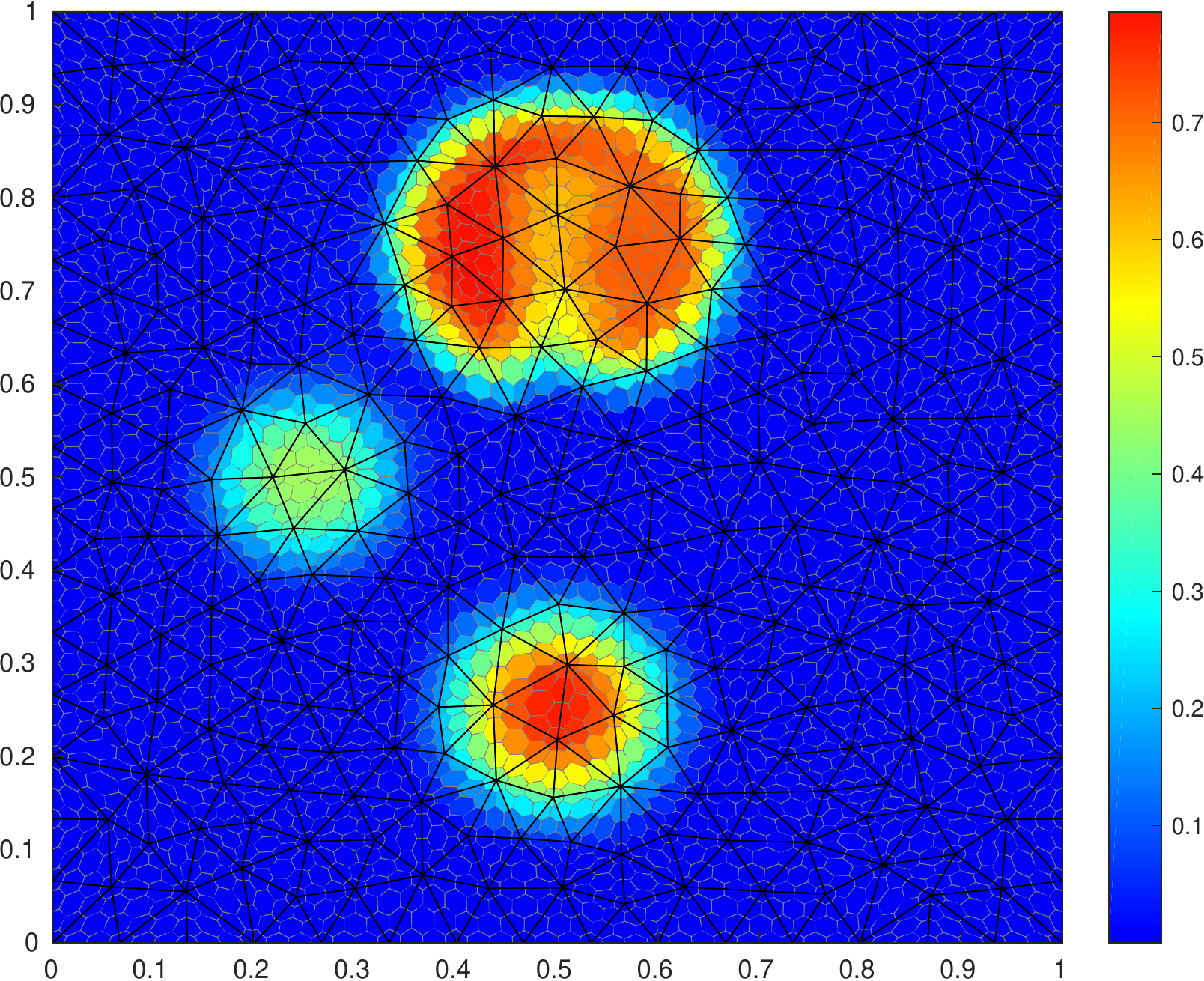}\label{fig_rotate_order4_poly2}}
    \caption{4th-order APLSC-DG solutions for rigid rotation on 576 cells after five full rotations: polygonal subdivision.}
  \label{fig_rotate_4th_poly}
  \end{center}
\end{figure}
In the light of Figures~\ref{fig_rotate_4th_cart} and \ref{fig_rotate_4th_poly}, we once again note how more uniform cell subdivision lead to less diffused solutions.
\begin{figure}[!ht]
  \begin{center}
    \subfigure[Solution profile for $y=0.25$.]{\includegraphics[height=5.8cm]{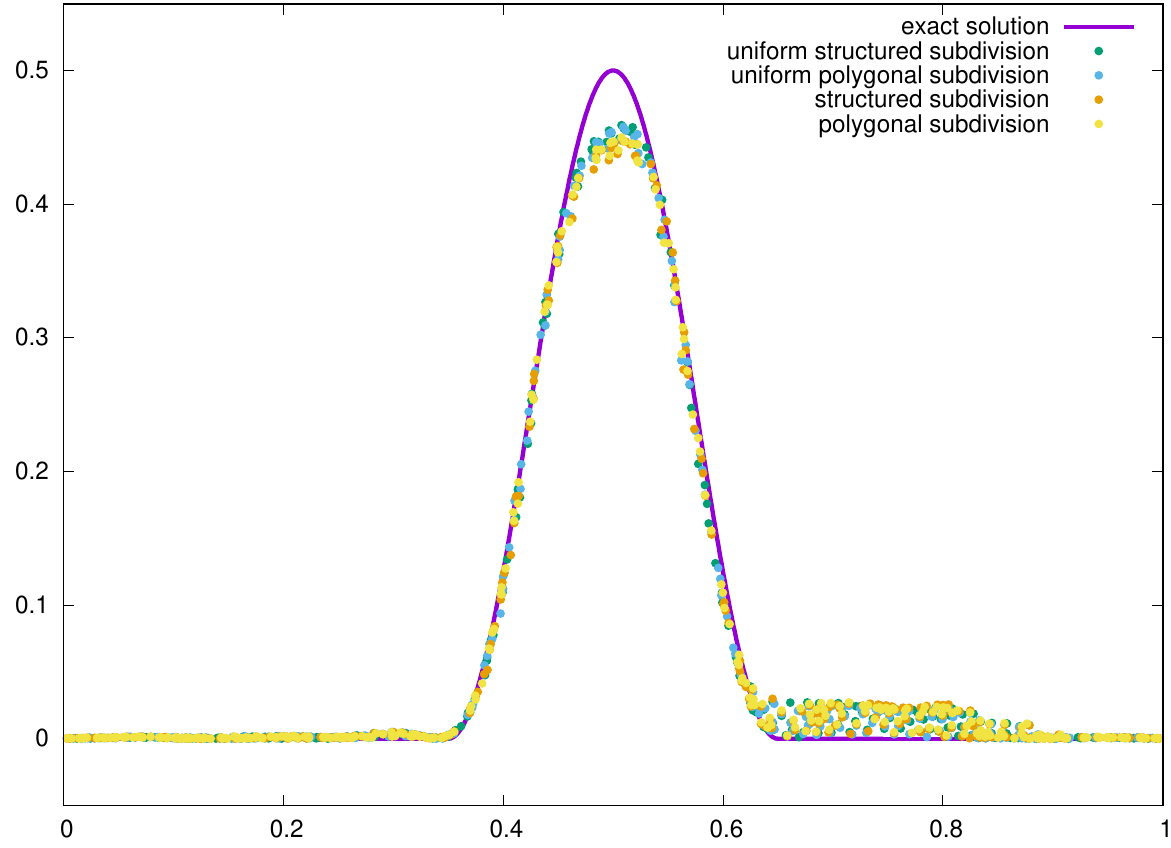}\label{fig_rotate_4th_subdiv_profile_a}}\hspace*{2mm}
    \subfigure[Solution profile for $x=0.25$.]{\includegraphics[height=5.8cm]{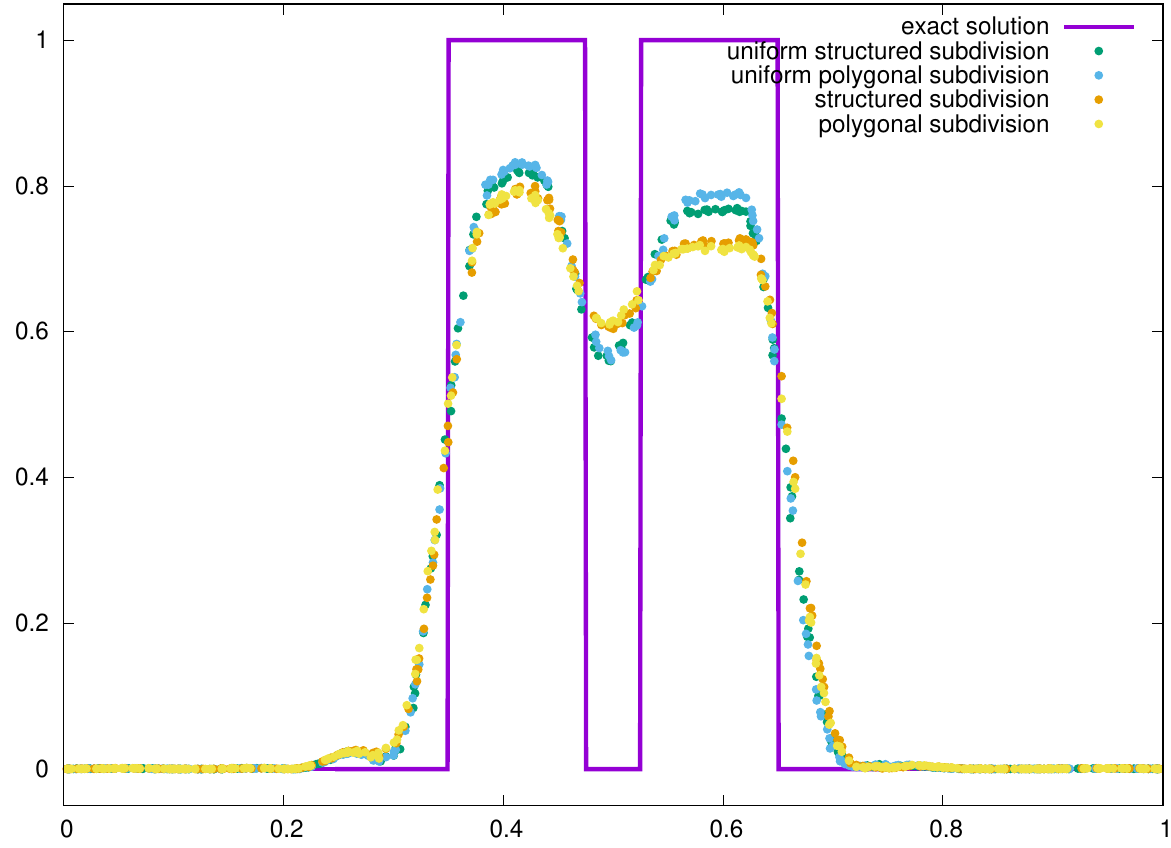}\label{fig_rotate_4th_subdiv_profile_c}}

    \caption{4th-order APLSC-DG solutions for rigid rotation on 576 cells after five full rotations: solution profiles.}
  \label{fig_rotate_subdiv_profile}
  \end{center}
\end{figure}
The solutions cross-sections displayed in Figure~\ref{fig_rotate_subdiv_profile} confirm this statement. However, it is worth mentioning that the difference in the result is not as critical as it was in the linear advection case. Furthermore, even if the uniform structured subdivision is not rotation invariant, it led to comparable results to the uniform Voronoi-type subdivision, while being a lot more simpler to implement and to generalize to any order of accuracy.

\vspace*{5mm}
\subsection{Non-linear case}
\label{subsect_nonlinear}

Let us now assess the performance and accuracy of our APLSC-DG technique in the 2D non-linear case. Both cases of scalar conservation laws as well as system of conservation laws will be addressed. Similarly to the subsection devoted to the linear case, let us defined set $\mc{N}(S_m^c)$ when considering the NAD criterion on subcell $S_m^c$. Unlike the linear case, we make use here of a subcell-wise DMP, meaning $\mc{N}(S_m^c)$ will be constituted by subcell $S_m^c$, as well as all its face and node neighboring subcells of $S_q^v$, either they belong to the same cell or not. By introducing $\mc{P}_m^c$ the set of vertices of subcell $S_m^c$ as well as $\mc{N}_p$ the set of subcells that share $\bs{x}_p$ as a vertex, this definition can be rewritten as

\begin{align}
  \label{subNAD_set}
  \mc{N}(S_m^c)=\Bigcup_{ \bs{x}_p\,\in\, \mc{P}_m^c}\,\mc{N}_p.
\end{align}

\begin{figure}[!ht]
  \begin{center}
    \subfigure[Structured subdivision.]{\includegraphics[height=7.cm]{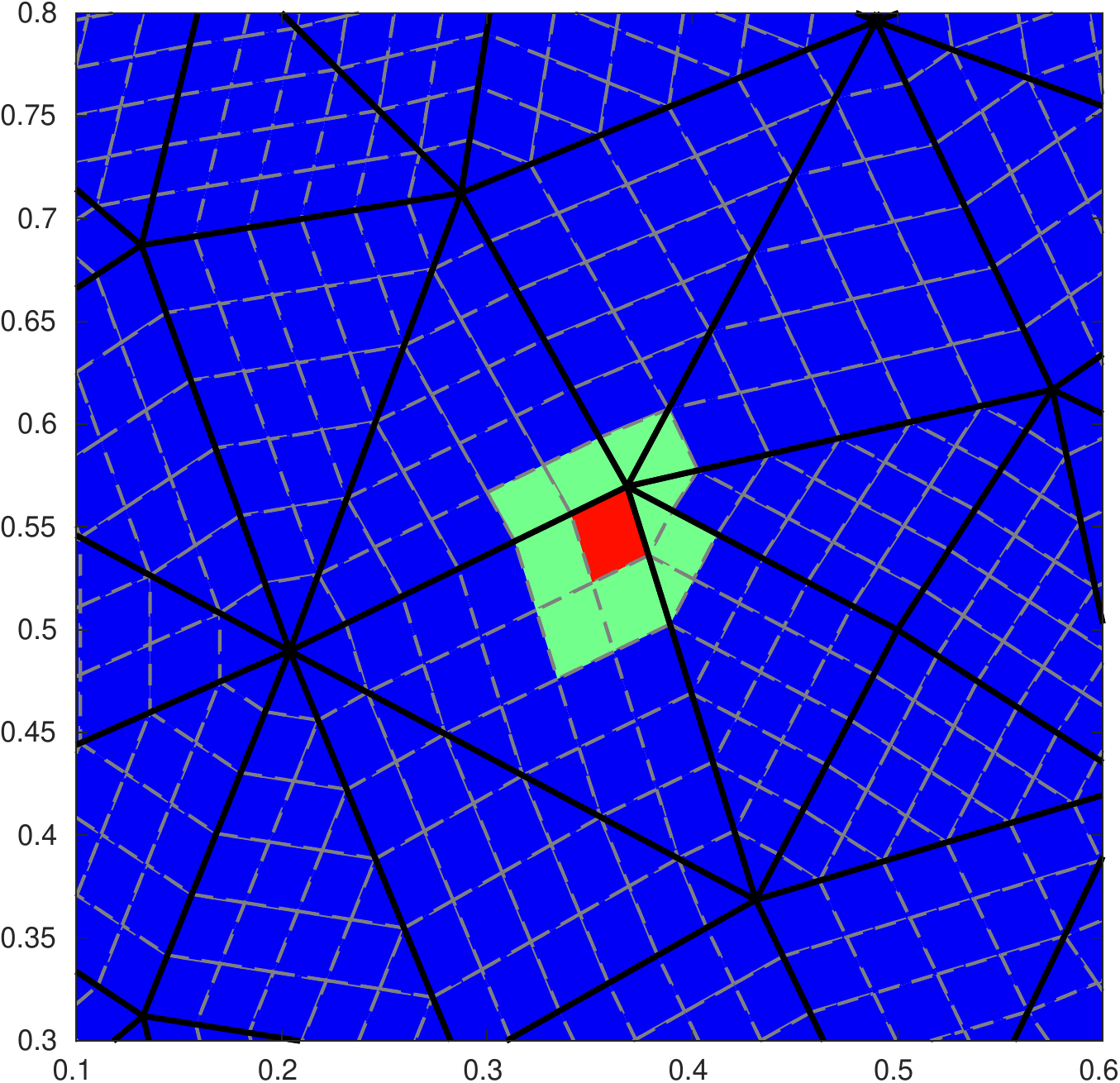}}\hspace*{12.mm}
    \subfigure[Voronoi-type subdivision.]{\includegraphics[height=7.cm]{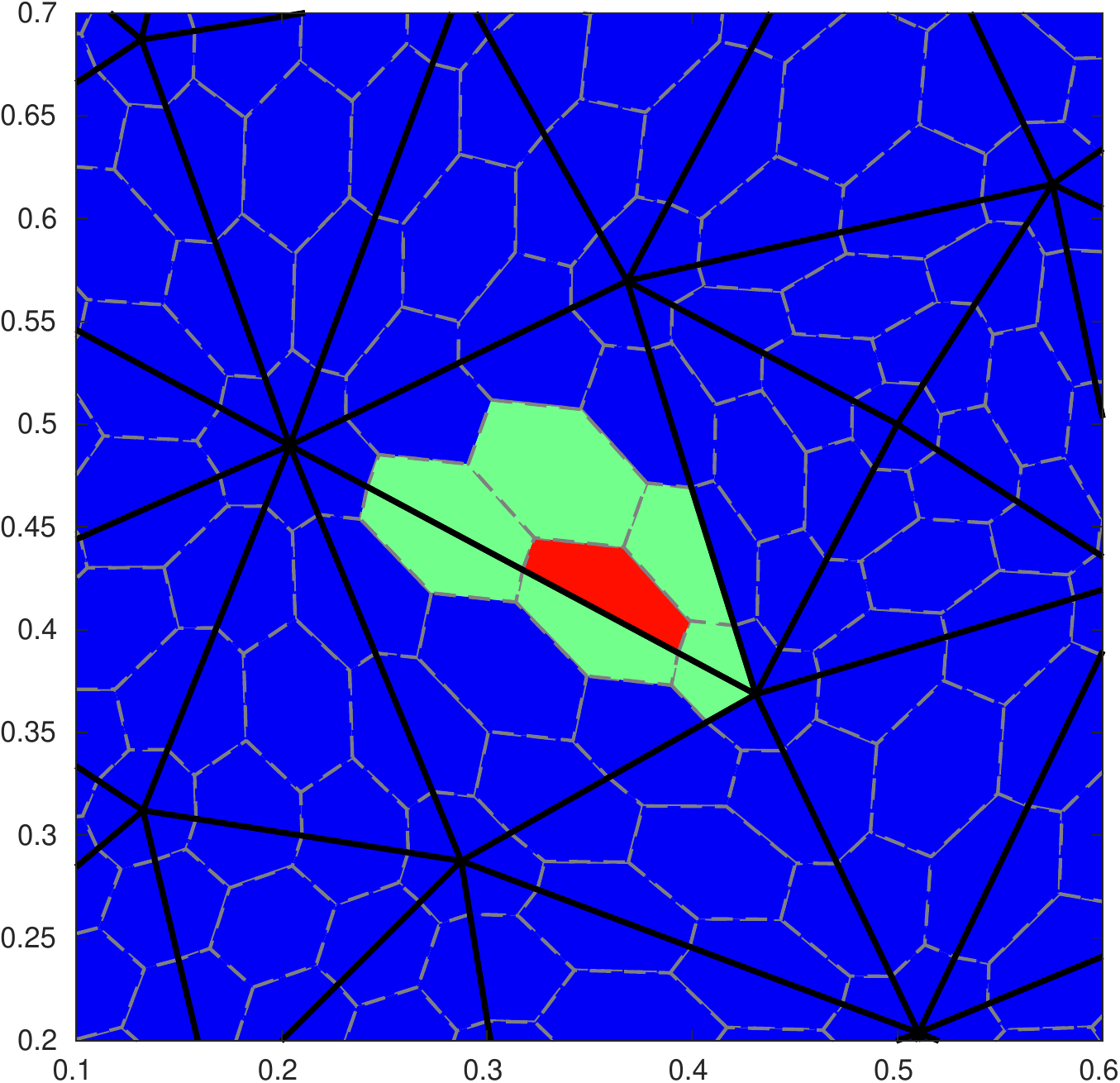}}
    \caption{Neighboring subcells set $\mc{N}(S_m^c)$ for the NAD criterion in the non-linear case: subcell $S_m^c$ is colored red while the subcells in $\mc{N}(S_m^c)$ are colored green.}
    \label{NAD_non-linear}
  \end{center}
\end{figure}

This particular set is depicted in Figure~\ref{NAD_linear}, for both the simple structured subdivision as well as the polygonal Voronoi-type one. In this figures, the subcell $S_m^c$ under consideration would be colored red, while the subcells constituting $\mc{N}(S_m^c)$ would be colored green. Let us emphasize that subcell $S_m^c$ is also part of $\mc{N}(S_m^c)$.

\vspace*{3mm}
\subsubsection{Burgers equation with a smooth initial solution}
\label{subsubsect_2D_burgers_sinus}

\begin{figure}[!ht]
  \begin{center}
    \subfigure[Solution map.]{\includegraphics[height=6.cm]{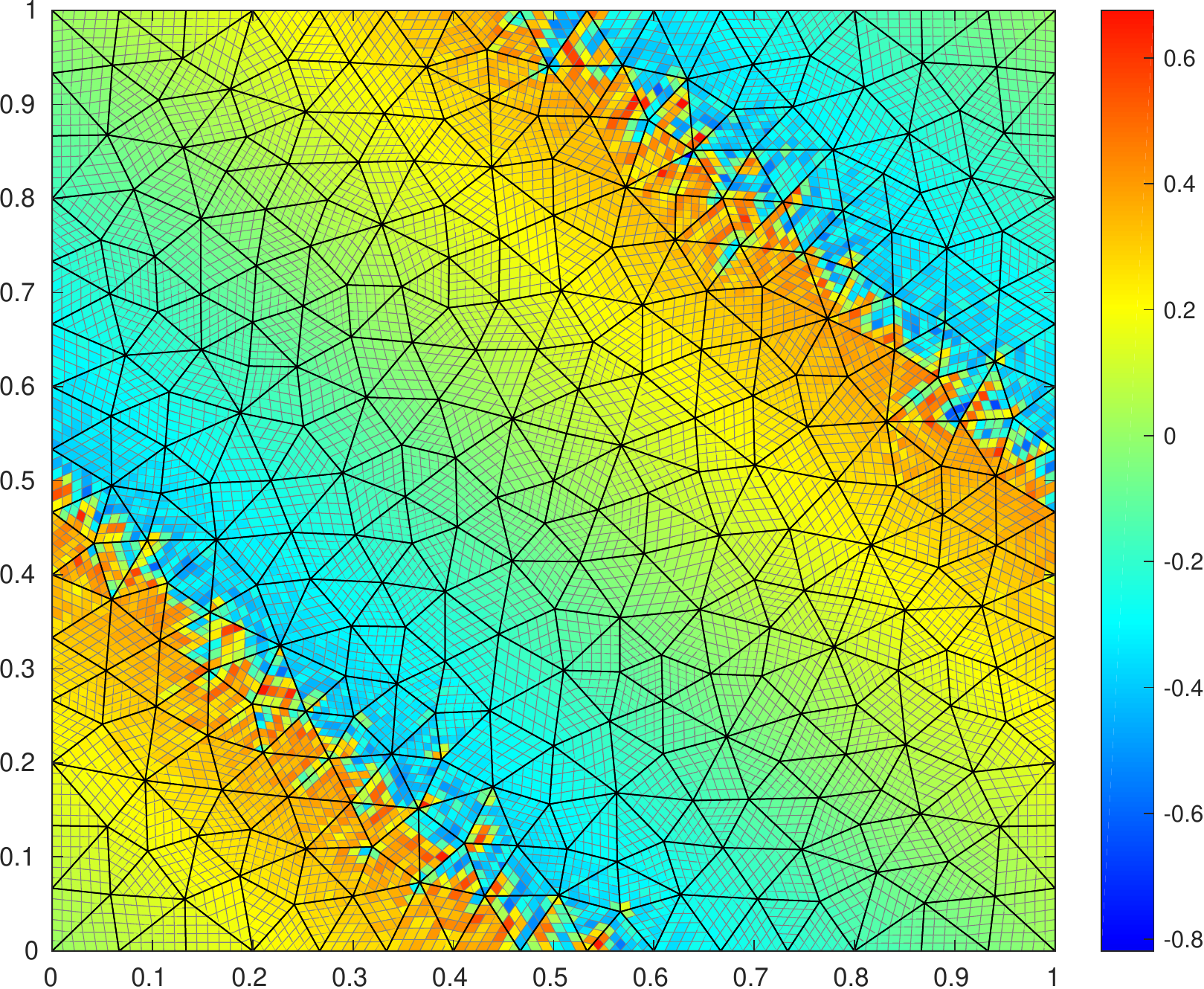}}\hspace*{5.mm}
    \subfigure[Submean values versus $(x+y-1)$ coordinate.]{\includegraphics[height=6.cm]{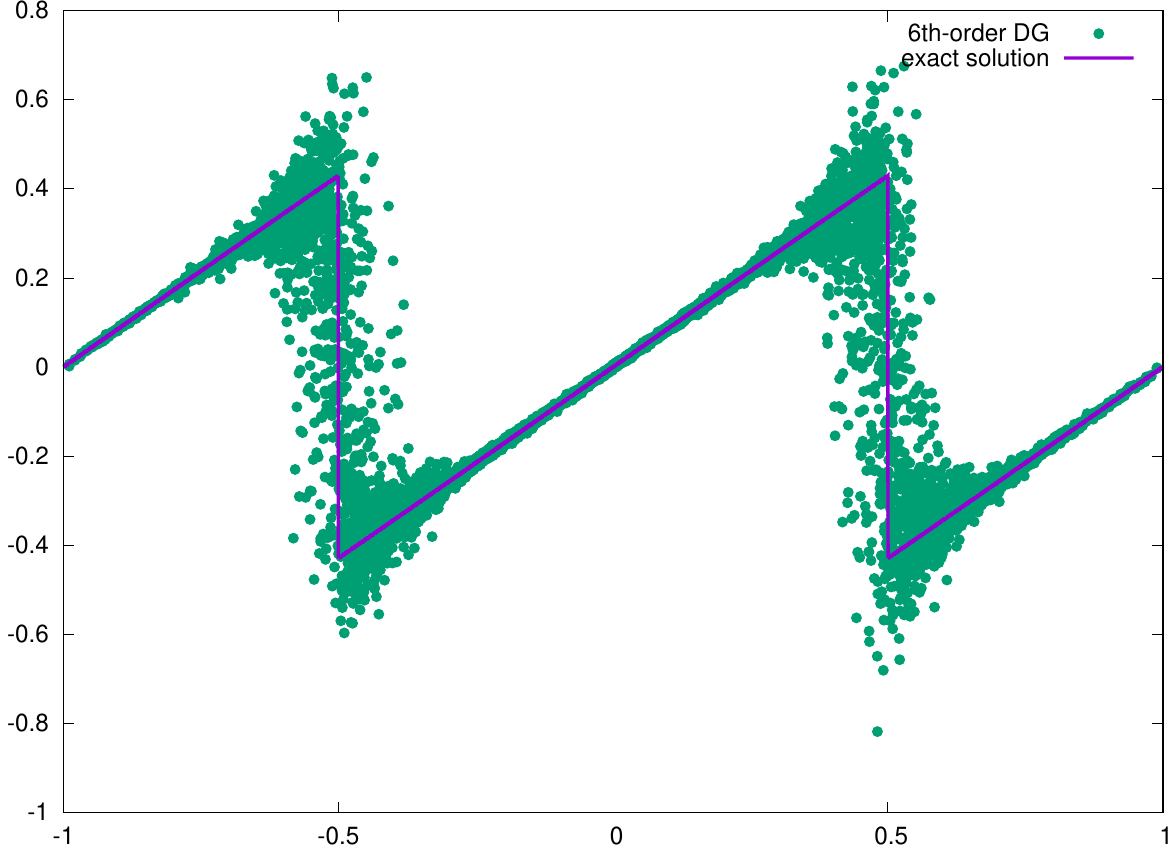}}
    \caption{6th-order uncorrected DG solution for 2D Burgers equation on a \ncell\, cells mesh at $t=0.5$.}
    \label{fig_burgers_2D_sinus_nolim}
  \end{center}
\end{figure}

To highlight the efficiency of the developed APLSC-DG scheme in the non-linear case, let us first consider the 2D Burgers equation
\begin{subequations}
  \label{eq_nonlinear}
\begin{empheq}[left = \empheqlbrace\,]{align}
&\vdt u(\bs{x},t) + \divx{\bs{F}\(u(\bs{x},t)\)}=0, && (\bs{x},t)\in\,[0,\, 1]^2\times[0,T], \label{pde_nonlinear}\\[2mm]
&u(\bs{x},0)=u_0(\bs{x}), &&\bs{x}\in\,[0,\, 1]^2, \label{pde_ini}
\end{empheq}\\[-4mm]
\end{subequations}

where the convex flux function writes $\bs{F}(u)=\demi\(u^2,\,u^2\)\tra$. As seen previously, starting from the smooth initial condition $u_0(x)=sin(2\pi\,(x,y))$ on $[0,1]^2$, two stationary discontinuities form along the lines $\{(x,y)\in\,[0,1]^2,\,x+y=\demi\}$ and $\{(x,y)\in\,[0,1]^2,\,x+y=\frac{3}{2}\}$. To emphasize how important a limiter or a correction technique is needed in this non-linear context, we first represent the numerical solution obtained by means of the 6th-order uncorrected DG on a very coarse grid made of 242 cells, see Figure~\ref{fig_burgers_2D_sinus_nolim}.
One can see how oscillating the numerical solution is. Furthermore, the two shocks are absolutely not well resolved.\\

\begin{figure}[!ht]
  \begin{center}
    \subfigure[Solution map.]{\includegraphics[height=6.cm]{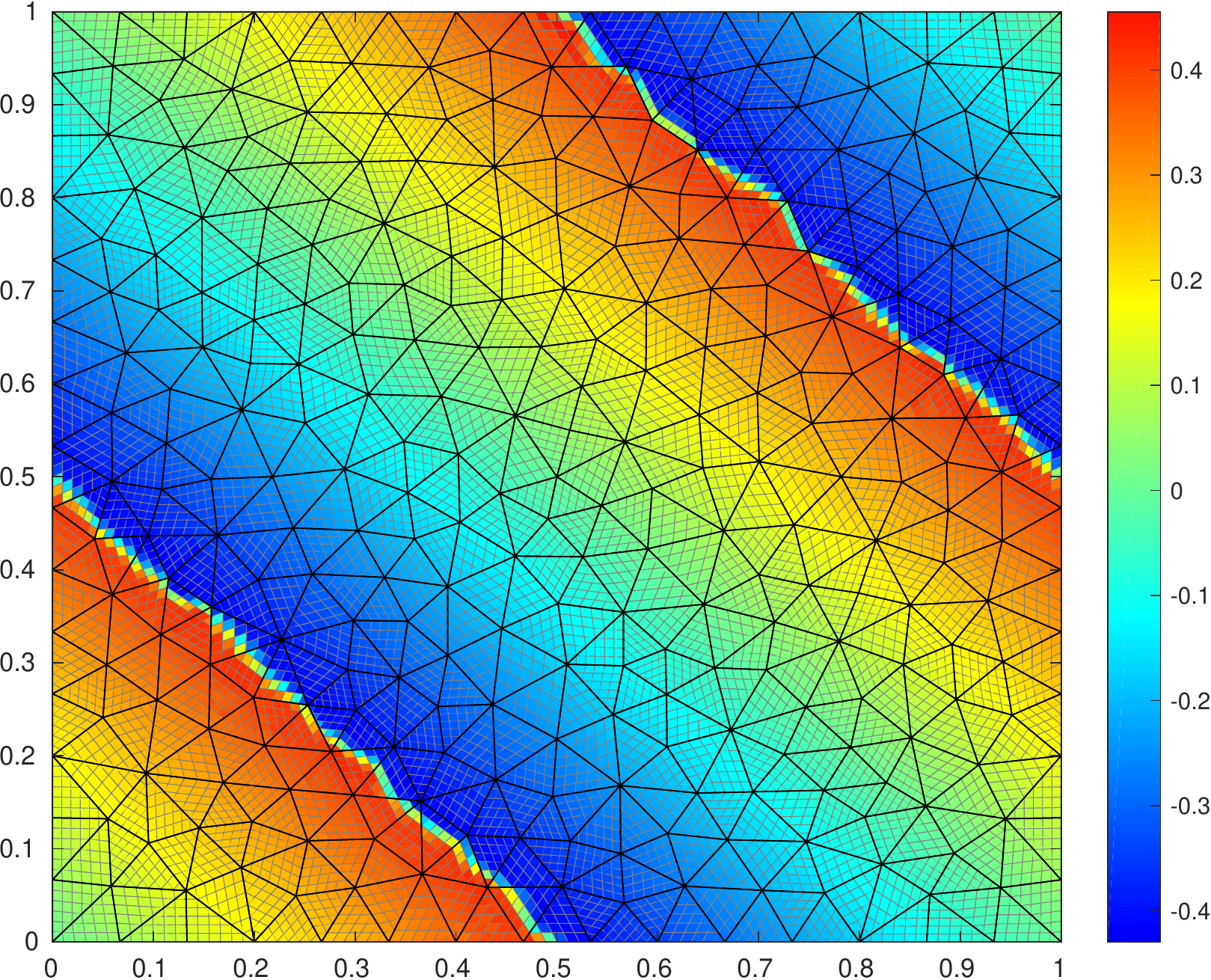}}\hspace*{5.mm}
    \subfigure[Submean values versus $(x+y-1)$ coordinate.]{\includegraphics[height=6.cm]{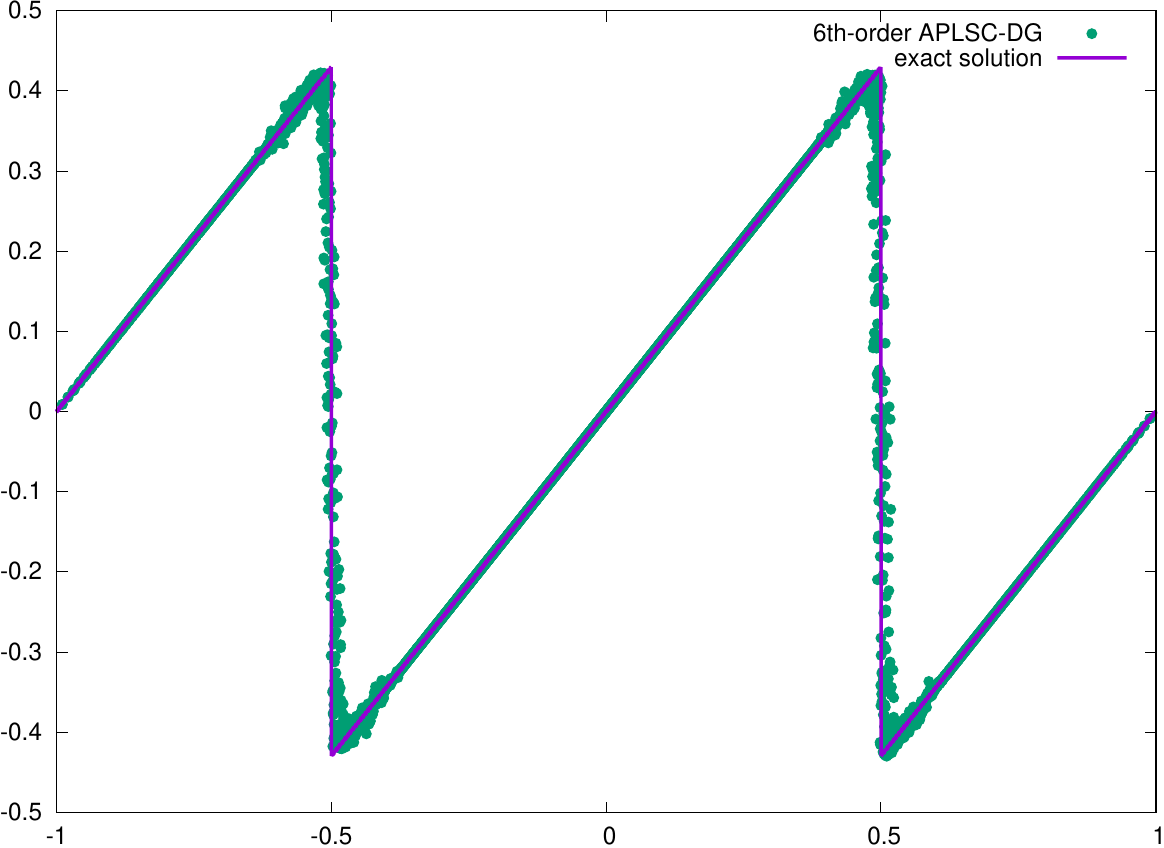}}
    \caption{6th-order APLSC-DG solution for 2D Burgers equation on a \ncell\, cells mesh at $t=0.5$.}
    \label{fig_burgers_2D_sinus}
  \end{center}
\end{figure}

In Figure~\ref{fig_burgers_2D_sinus}, the numerical solution obtained with the 6th-order APLSC-DG scheme is illustrated at time $t=0.5$. We can see in this totally anisotropic triangles coarse grid, the corrected scheme quite accurately recovers the two straight line shocks, while ensuring a robust low oscillatory behavior. Now, to investigate once more if the cell subdivision type has any influence on the quality of the numerical solution obtained, we simulate this two-shocks Burgers test case with an even coarser grid with the four subdivisions depicted in Figures~\ref{fig_subdiv_cart} and \ref{fig_subdiv_poly}.

\begin{figure}[!ht]
  \begin{center}
    \subfigure[Uniform structured subdivision.]{\includegraphics[height=6.3cm]{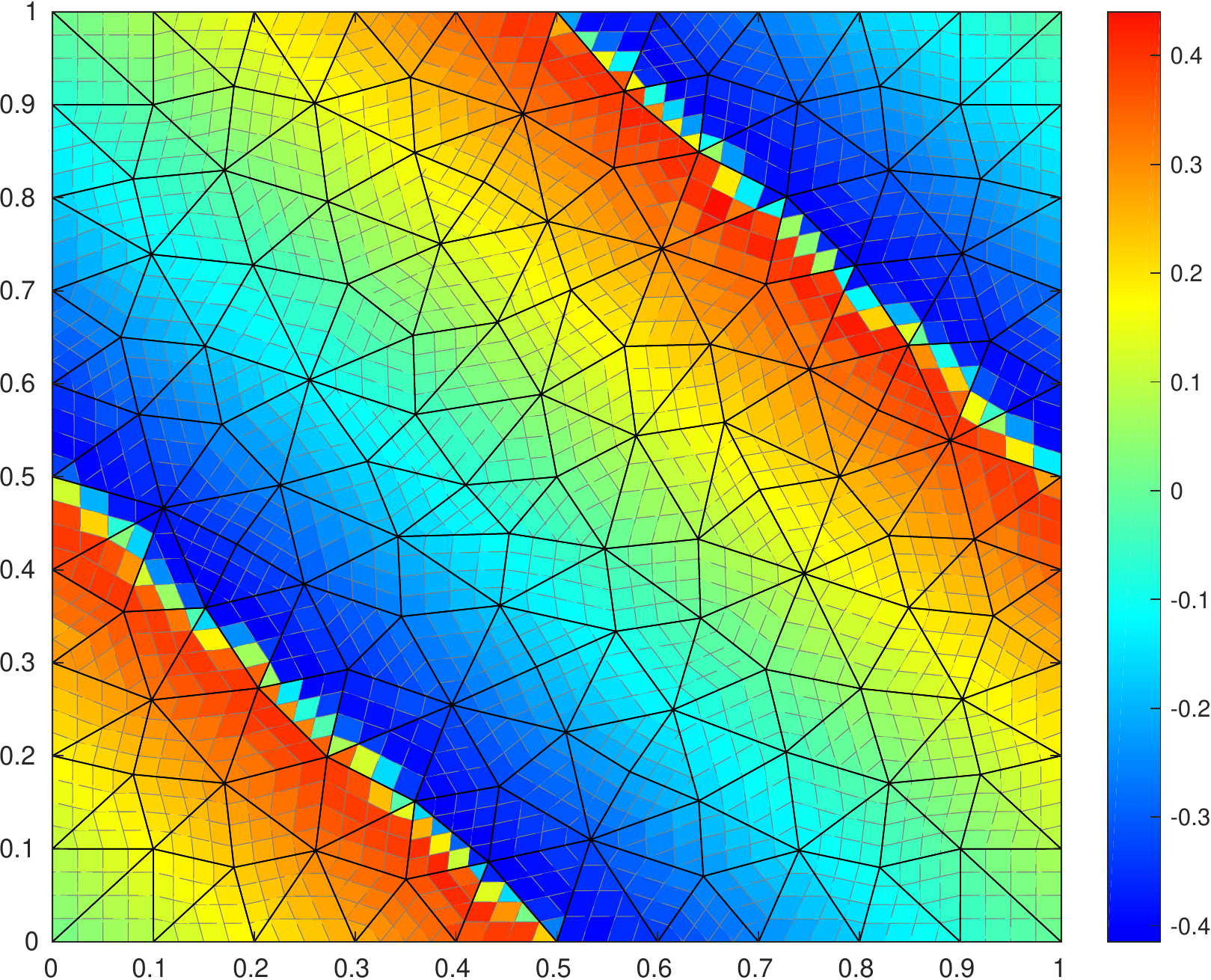}\label{fig_burgers_order4_cart1}}\hspace*{5mm}
    \subfigure[Non-uniform structured subdivision.]{\includegraphics[height=6.3cm]{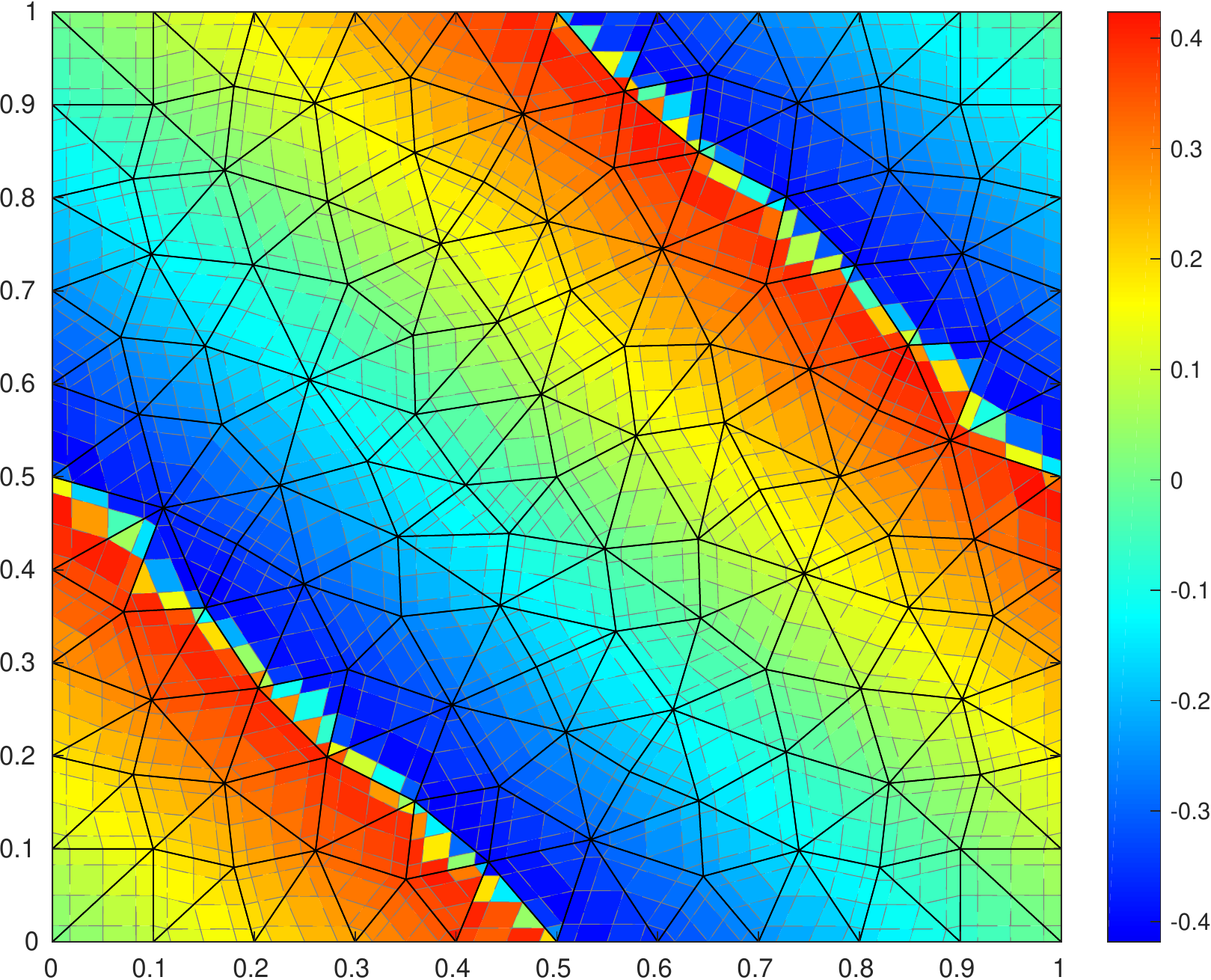}\label{fig_burgers_order4_cart3}}
    \caption{4th-order APLSC-DG solutions for 2D Burgers equation on 242 cells at $t=0.5$: structured subdivision.}
  \label{fig_burgers_4th_cart}
  \end{center}
\end{figure}
\begin{figure}[!ht]
  \begin{center}
    \subfigure[Uniform polygonal subdivision.]{\includegraphics[height=6.3cm]{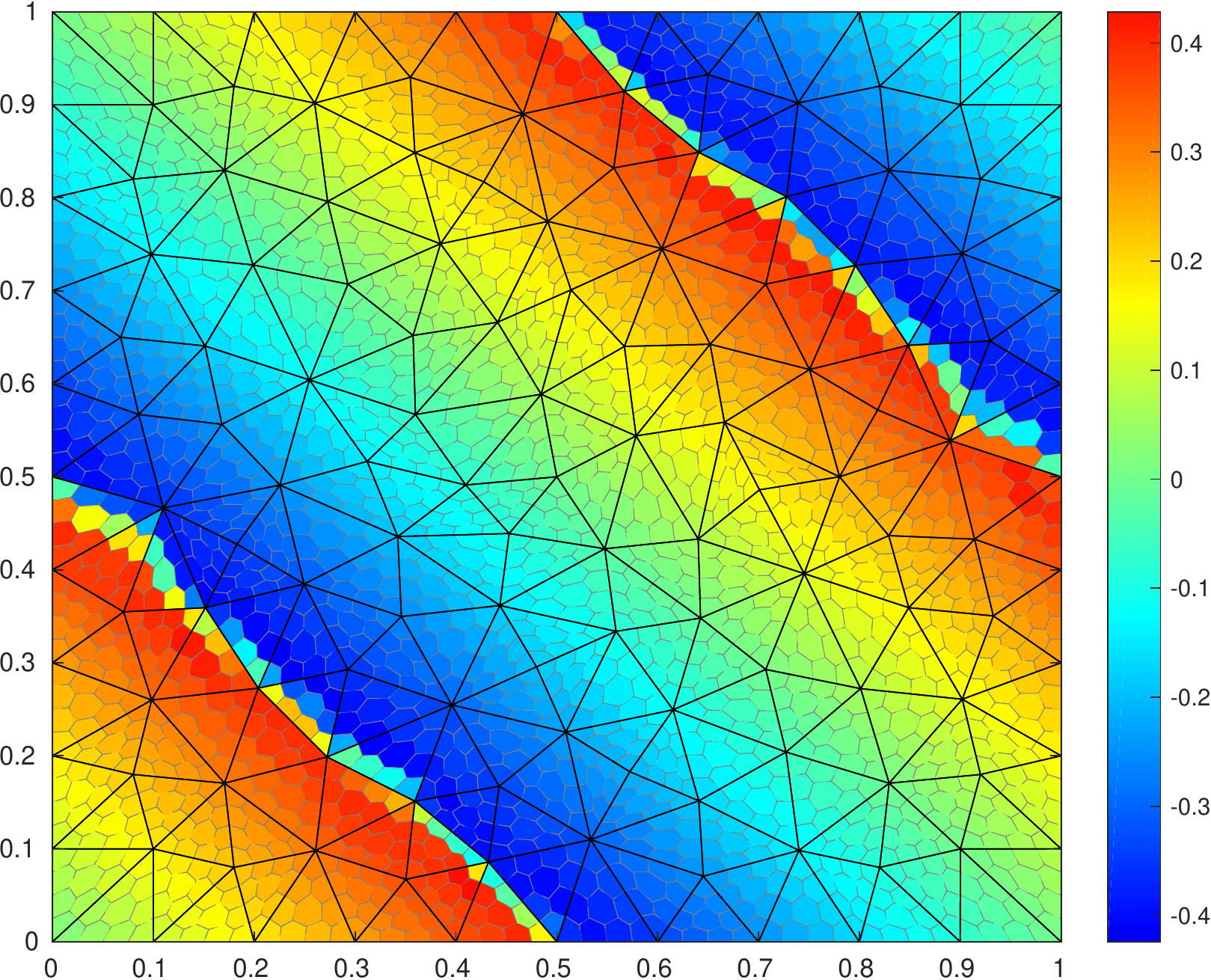}\label{fig_burgers_order4_poly1}}\hspace*{5mm}
    \subfigure[Non-uniform polygonal subdivision.]{\includegraphics[height=6.3cm]{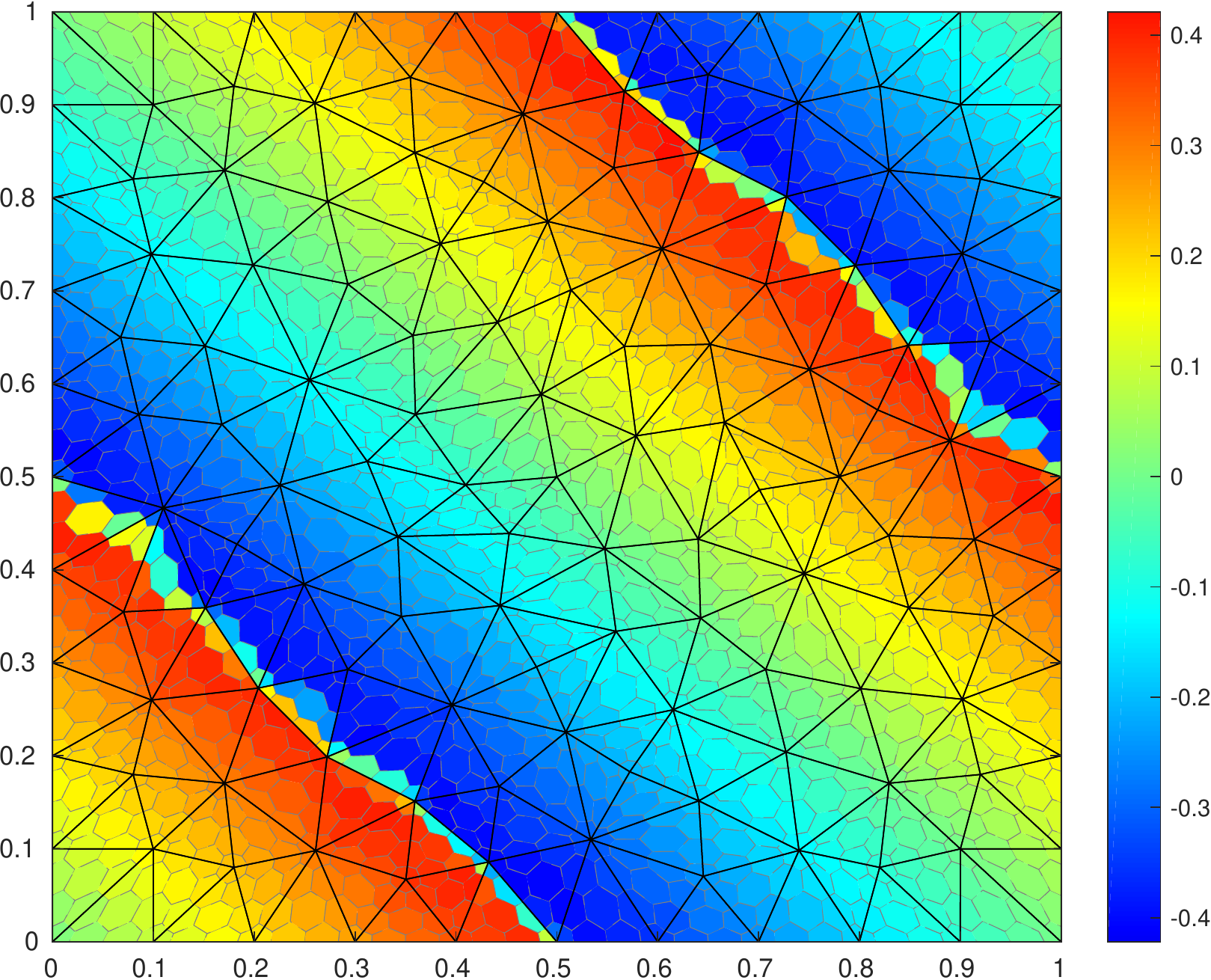}\label{fig_burgers_order4_poly2}}
    \caption{4th-order APLSC-DG solutions for 2D Burgers equation on 242 cells at $t=0.5$: polygonal subdivision.}
  \label{fig_burgers_4th_poly}
  \end{center}
\end{figure}

\begin{figure}[!ht]
  \begin{center}
    \includegraphics[height=7.5cm]{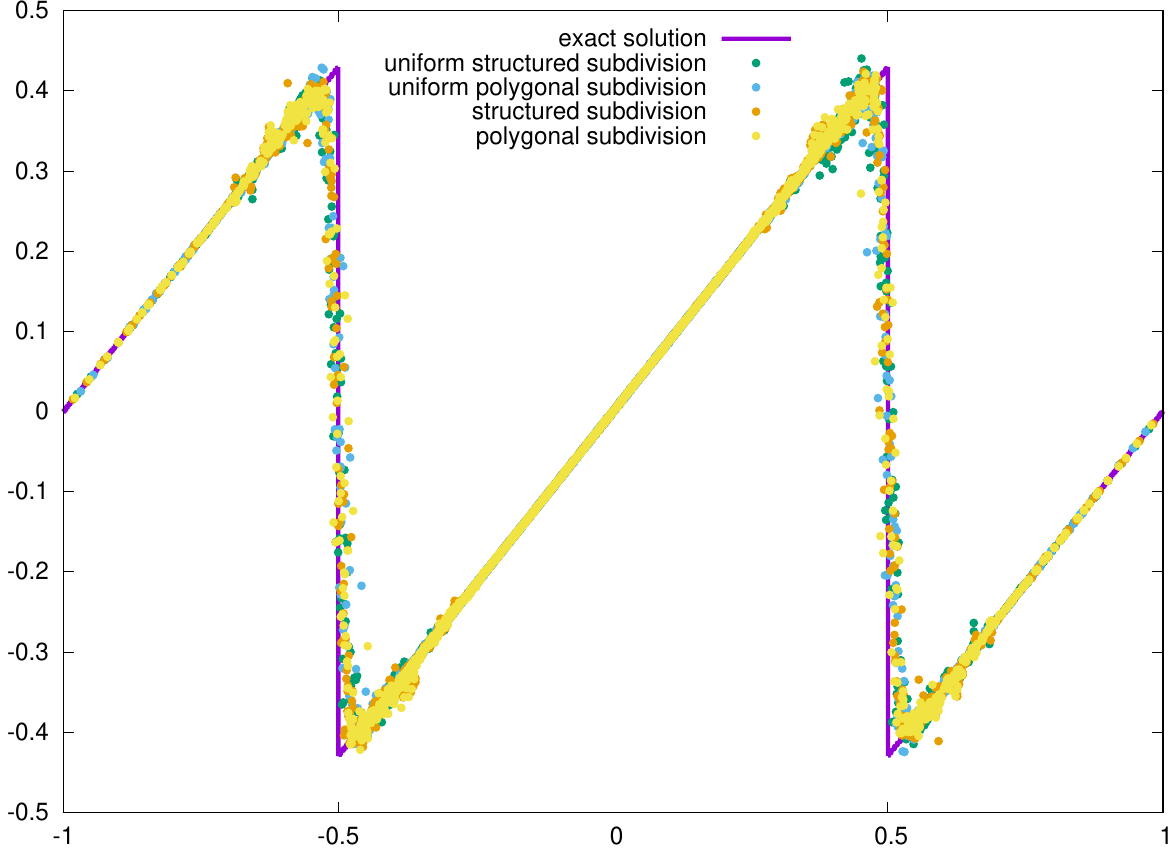}
    \caption{4th-order APLSC-DG solutions for 2D Burgers equation on 242 cells at $t=0.5$: submean values versus $(x+y-1)$ coordinate.}
    \label{fig_burgers_4th_profile}
  \end{center}
\end{figure}

In the light of Figure~\ref{fig_burgers_4th_profile}, the four different subdivisions seem to produce equivalent results, which are further quite satisfactory considering the extremely coarse grid used. However, the use of the uniform structured cell partition, Figure~\ref{fig_burgers_order4_cart1}, appears to capture in a sharper fashion the two straight line shocks. As we have shown that uniform subdivision, structured or Voronoi-type, lead to better results when APLSC-DG scheme is used, only those two subdivisions will be used in the remainder of the article.

\subsubsection{KPP problem}
\label{subsubsect_2D_kpp}

Before investigating the non-linear system case, we now turn our attention to non-linear conservation laws with non-convex fluxes. To this end, we consider the KPP problem proposed by Kurganov, Petrova, Popov (KPP) in \cite{Popov} to test the convergence properties of some WENO schemes in the context of non-convex fluxes. For this particular problem, we study the non-linear problem \eqref{eq_nonlinear} where the flux function is given by $F(u)=\(\sin(u),\,\cos(u)\)\tra$. Considering the computational domain $[-2,2]\times[-2.5,1.5]$, the initial condition reads as follows
\begin{align*}
  u_0(x)=\left\{\begin{array}{ll}
  \Frac{7\,\pi}{2}\qquad &\text{if }\; x<\demi,\\[3mm]
  \Frac{\pi}{4}\qquad &\text{if }\; x>\demi.
  \end{array}\right.
\end{align*}

\begin{figure}[!ht]
  \begin{center}
    \subfigure[DG solution map.]{\includegraphics[height=6.5cm]{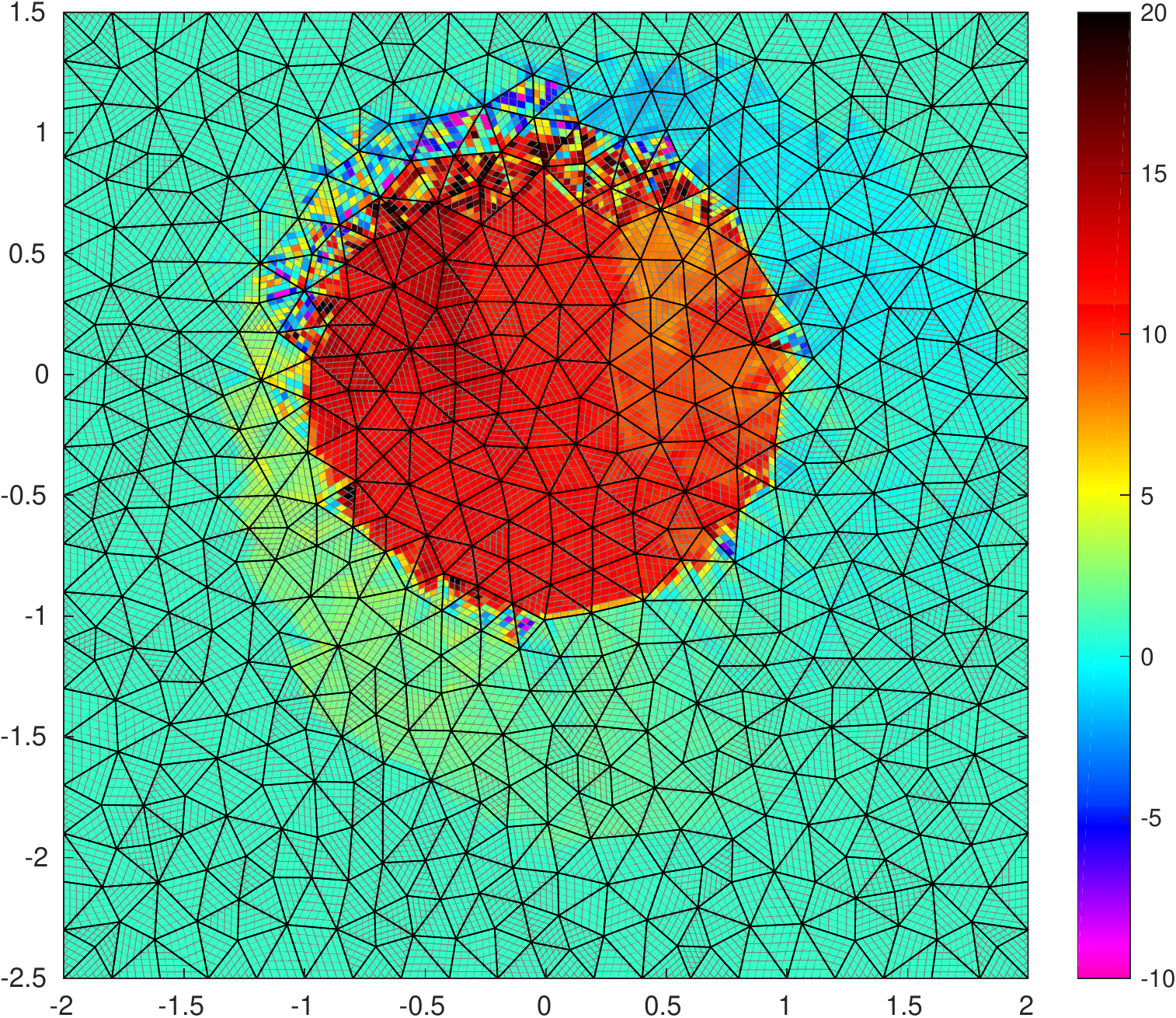}\label{fig_kpp_order6_nocorr}}\hspace*{4mm}
    \subfigure[APLSC-DG solution map.]{\includegraphics[height=6.5cm]{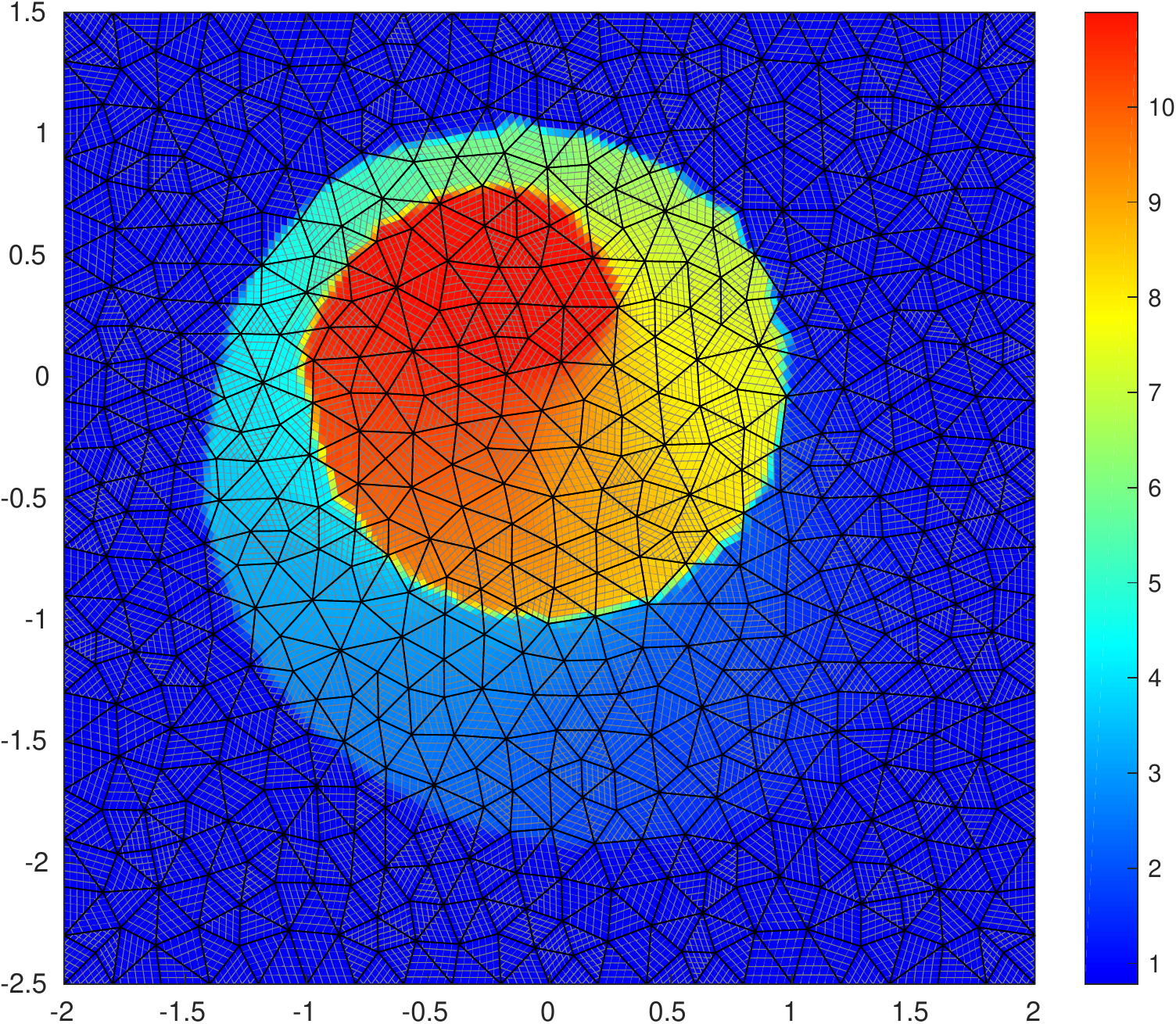}\label{fig_kpp_order6_sol}}
    
    \subfigure[Corrected subcells map.]{\includegraphics[height=6.5cm]{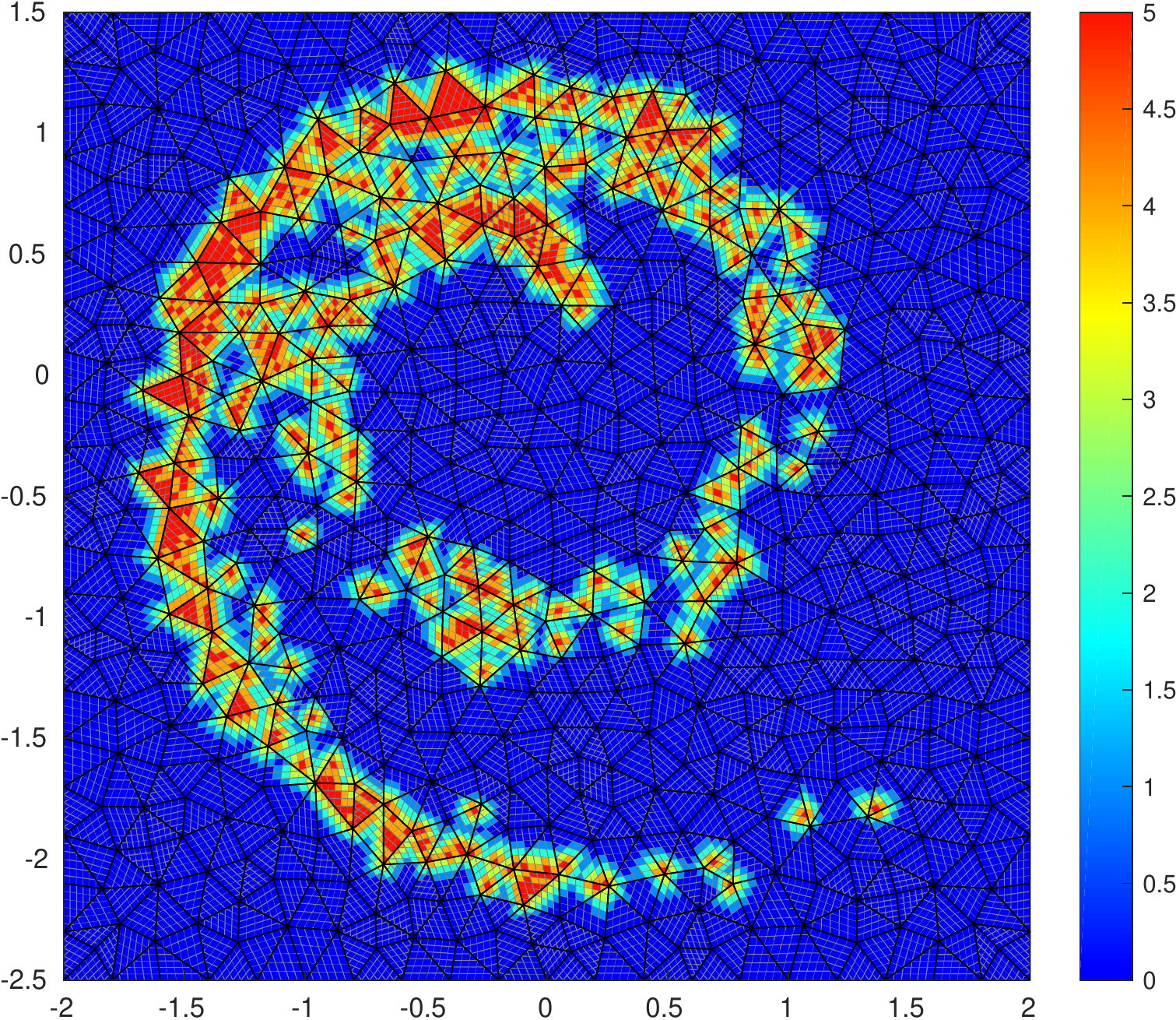}\label{fig_kpp_order6_detect}}
    \caption{6th-order APLSC-DG solution for the KPP problem on a 1054 cells mesh at time $t=1$.}
  \label{fig_kpp_order6}
  \end{center}
\end{figure}

This test is very challenging to many high-order schemes as the solution has a two-dimensional composite wave structure, and as generally numerical methods fail to converge to the unique entropic exact solution. In most cases, to be able to capture such rotation composite structure, very fine grids must be used. Here, by means of 6th-order uncorrected DG and then APLSC-DG scheme, we make use of an unstructured mesh made of 1054 triangular cells, which is very coarse in this quite complex situation. Results are displayed in Figure~\ref{fig_kpp_order6}.

\begin{figure}[!ht]
  \begin{center}
    \subfigure[Structured subdivision.]{\includegraphics[height=6.5cm]{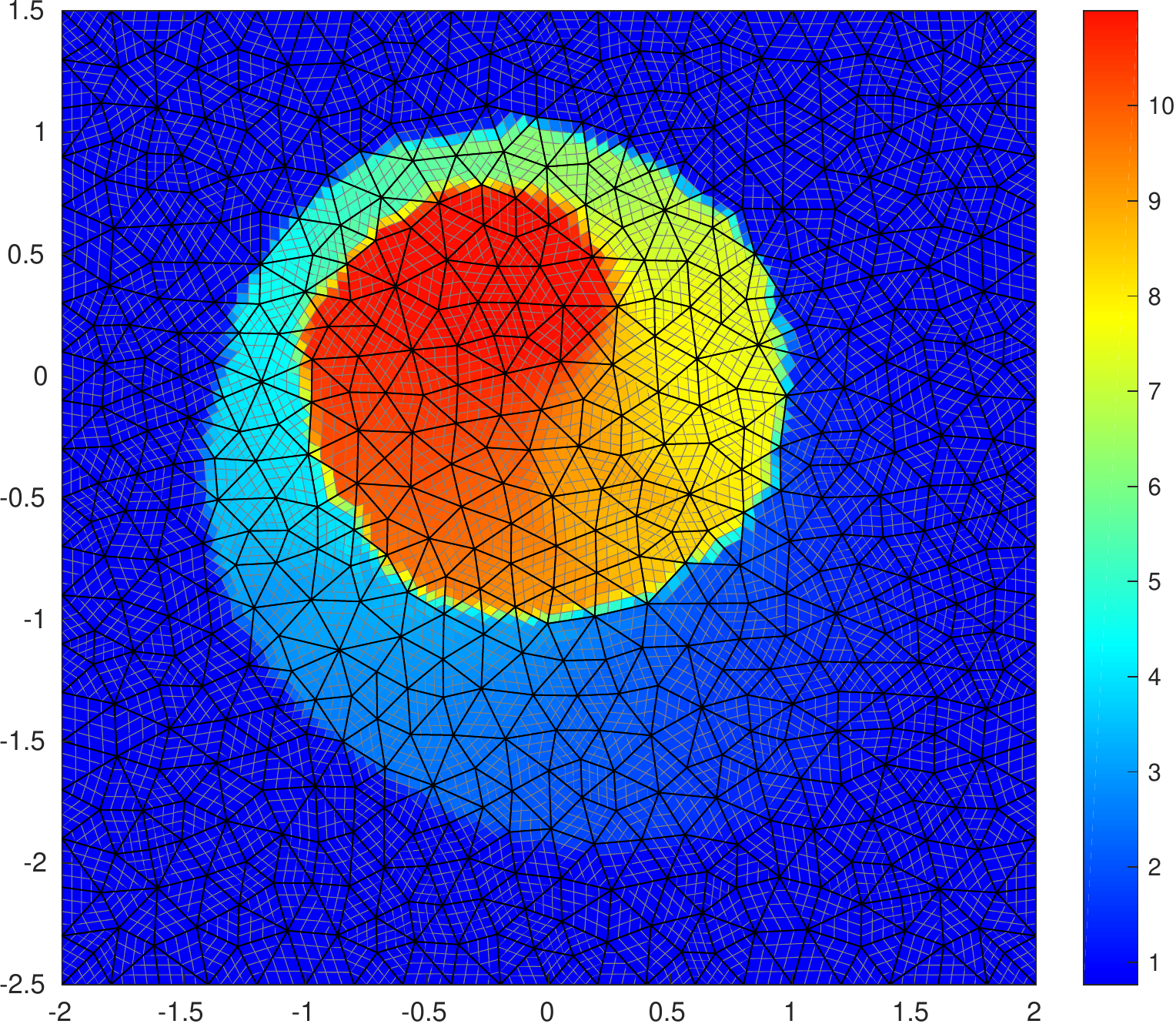}\label{fig_kpp_order4_cart}}\hspace*{4mm}
    \subfigure[Voronoi-type subdivision.]{\includegraphics[height=6.5cm]{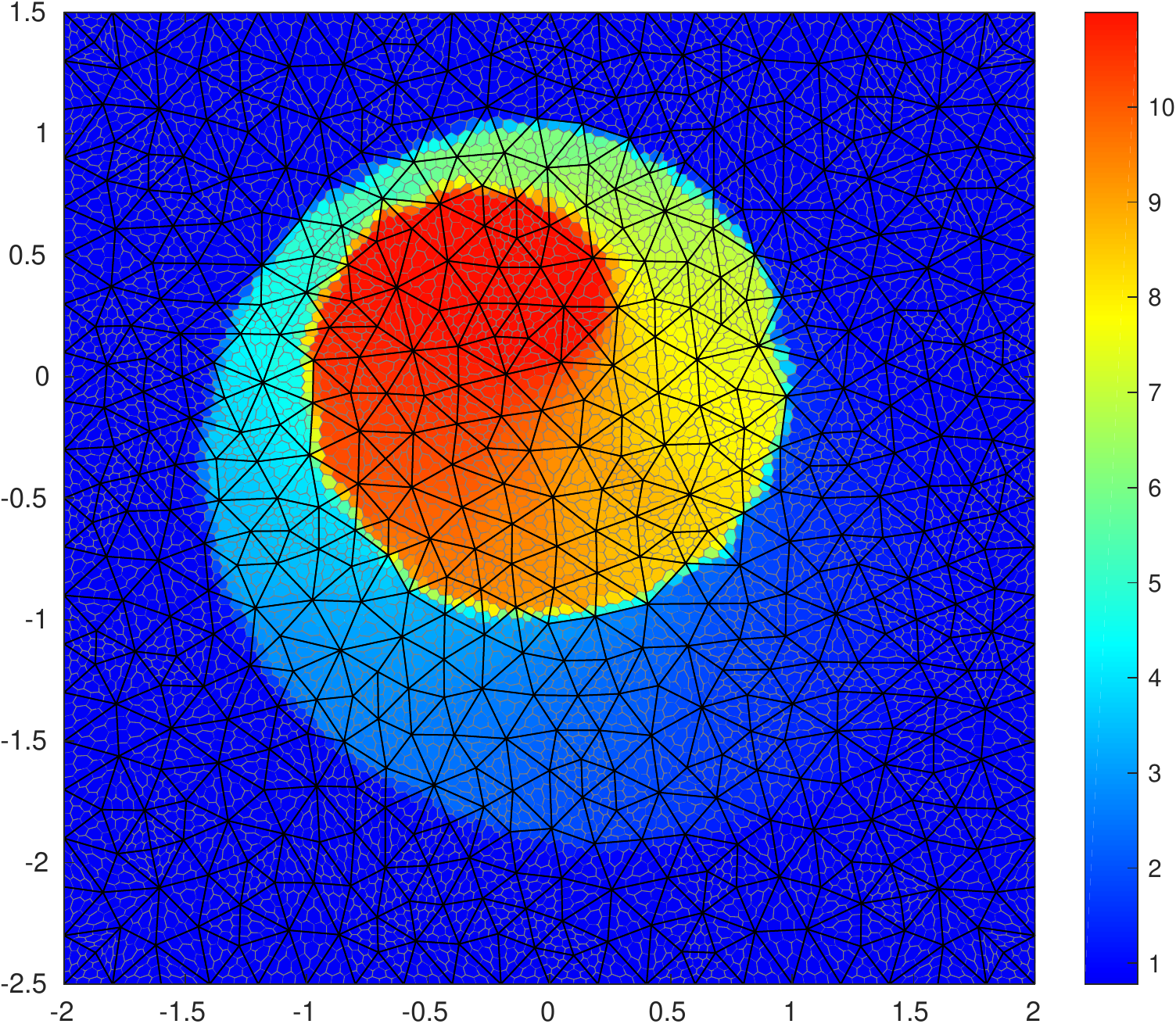}\label{fig_kpp_order4_poly}}
    \caption{4th-order APLSC-DG solution for the KPP problem on a 1054 cells mesh at time $t=1$: subdivisions comparison.}
  \label{fig_kpp_order4}
  \end{center}\vspace*{-5mm}
\end{figure}

Figure~\ref{fig_kpp_order6_nocorr} illustrates the non-entropic behavior of DG schemes in this non-convex flux case, as well as the high oscillatory solution produced if no additional mechanism is used. The application of our APLSC technique permits to capture to correct entropic solution, while avoiding the apparition of spurious oscillations, as displayed in Figure~\ref{fig_kpp_order6_sol}. Furthermore, we can observe that although the coarseness of the grid used, the APLSC-DG scheme allowed to recover the two-dimensional vortex-like wave structure of the solution. In Figure~\ref{fig_kpp_order6_detect}, subcells corrected during the different Runge-Kutta stages of the last time iteration are displayed, with different colors accordingly the amount of first-order correction applied. One can see how the NAD and PAD troubled subcell detection criteria accurately track the spiral discontinuity of the entropic exact solution. Here, a different colormap compared to Figure~\ref{fig_burgers_compar_detect} has been used here for a better readability of the results.
In Figure~\ref{fig_kpp_order4}, we once more assess the impact of the subdivision on the quality of the results obtained the APLSC-DG method. As the uniform ones have proved to yield better results, we only compare here the uniform structured subdivision with uniform Voronoi-type one. Anew, the 4th-order scheme is used here.
As one can see in Figure~\ref{fig_kpp_order4}, similar results have been obtained regardless the type of cell subdivision.\\[5mm]


\subsection{2D Euler system}
\label{subsect_Euler_2D}

To close this numerical application section, and assess once again the high capability of the \textit{a posteriori} local subcell correction technique presented here, the non-linear system case will be now addressed. To this end, let us consider the 2D Euler compressible gas dynamics system of equations
\begin{subequations}
\label{Euler_system}
\begin{empheq}[left = \empheqlbrace\,]{align}
  &\vdt\; \rho + \divx{\bs{q}}=0,\\[3mm]
  &\vdt\; \bs{q}+ \divx{\(\Frac{\bs{q}\otimes\bs{q}}{\rho}+p\,I_d\)}=\bs{0},\\[3mm]
  &\vdt\; E+ \divx{\((E+p)\,\Frac{\bs{q}}{\rho}\)}=0,
\end{empheq}\\[-4mm]
\end{subequations}

where the conserved variables $\rho$, $\bs{q}=\rho\,\bs{v}$ and $E$ respectively stand for the density, the momentum and the total energy, while $\bs{v}$ characterize the fluid velocity. The thermodynamic closure is given by the equation of state $p=p(\rho,\, \veps)$ where $\veps=E-\Frac{1}{2} \,\Frac{\|\bs{q}\|^2}{\rho}$ denotes the internal energy. In this paper, we make use of a gamma gas law, $i.e.$\; $p=(\gamma-1)\,\veps$, where $\gamma$ is the polytropic index of the gas.\\

Although the whole theory presented here has been introduced in the simple case of scalar conservation laws, the extension to the system case is perfectly straightforward. The only part which may differ is the troubled detection part. For the physical admissibility detection (PAD), we consider that a solution is physically admissible if the density and the internal energy are strictly positive. The use of other equations of state may lead a different convex set of admissibility, see \cite{vilar_lag_15_1D} for instance. For the numerical admissibility detection (NAD), the natural system counterpart would be to apply the previously introduced detection criteria to the Riemann invariants. However, in the non-linear system case, those quantities are not easy to get nor to manipulate. We could have use a linearized version of the Riemann invariants, as in \cite{vilar_lag_14} for instance, but for sake of simplicity we naively apply the NAD criterion to one of the conserved variable, similarly to what has been done in the 1D case \cite{vilar_aplsc_1D}. Here, we choose to work with the energy, as this physical quantity would be sensitive to any type of wave. Once again, the simple global Lax-Friedrichs numerical flux will be used in the remainder.

\newpage
\subsubsection{Sod shock tube problem.}
\label{subsubsect_Euler_sod}

We consider the extension of the classical Sod shock tube \cite{Sod} to the case of the cylindrical geometry. This problem consists of a cylindrical shock tube of unity radius. The interface is located at $r=0.5$. At the initial time, the states on the left and on the right sides of the interface are constant. The left state is a high pressure fluid characterized by $(\rho_0^L,p_0^L,\bs{v}_0^L)=(1,1,\bs{0})$, the right state is a low pressure fluid defined by $(\rho_0^R,p_0^R,\bs{v}_0^R)= (0.125,0.1,\bs{0})$. The gamma gas law is defined by $\gamma=\frac{7}{5}$. The computational domain is defined in polar coordinates by $(r,\theta)\in [0,1]\times [0,\frac{\pi}{4}]$. We prescribe symmetry boundary conditions at the boundaries $\theta=0$ and $\theta=\frac{\pi}{4}$, and an outflow condition at $r=1$. The exact solution consists of three circular waves, a shock followed by a contact discontinuity and rarefaction wave. The aim of this test case is then to assess the APLSC-DG scheme accuracy while ensuring a non-oscillatory behavior, and its ability to preserve the radial symmetry.

\begin{figure}[!ht]
  \begin{center}
    \subfigure[Subcell density values.]{\includegraphics[height=5.cm]{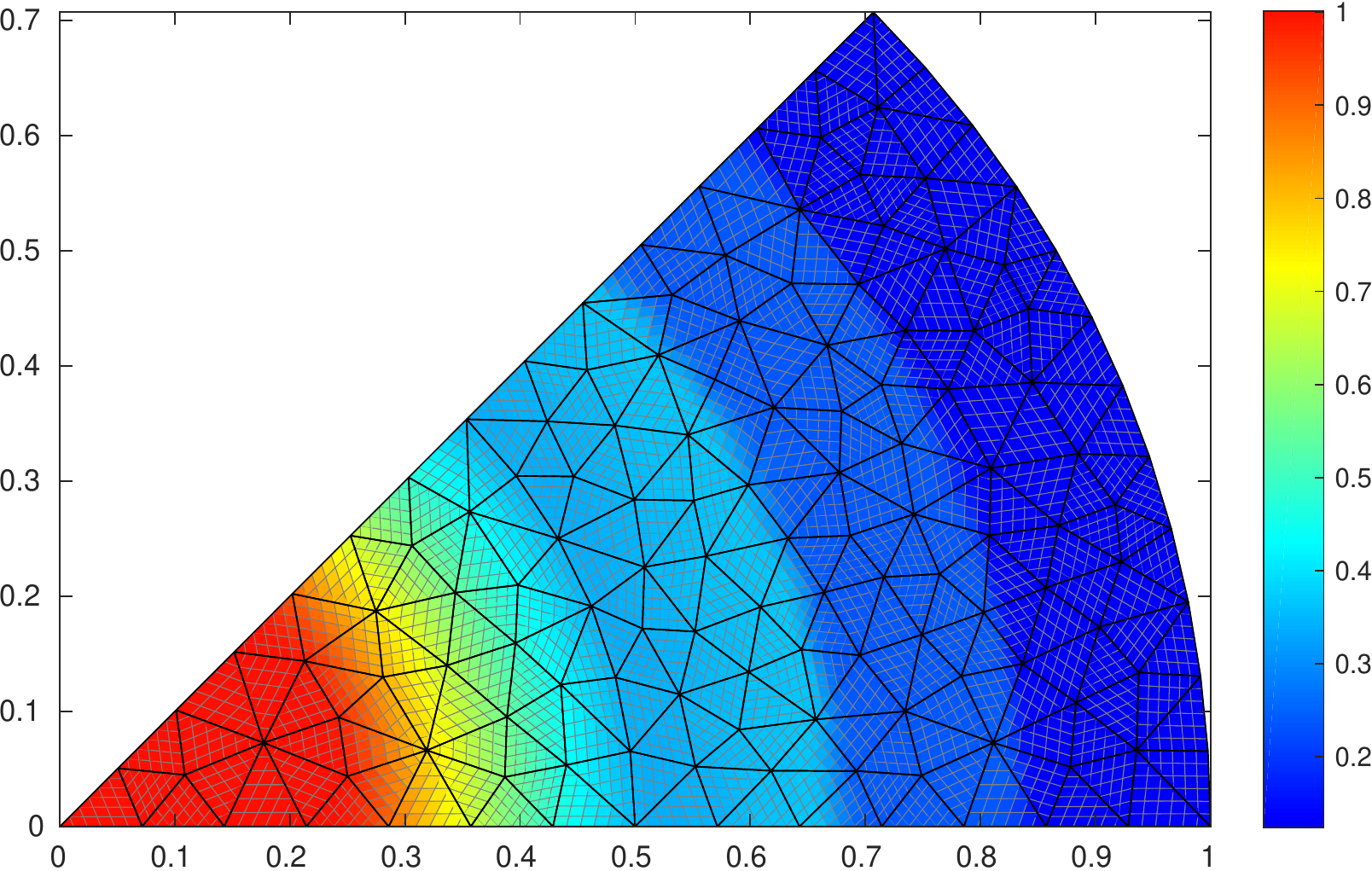}\label{fig_sod_den}}\hspace*{2mm}
    \subfigure[Corrected subcells map.]{\includegraphics[height=5.cm]{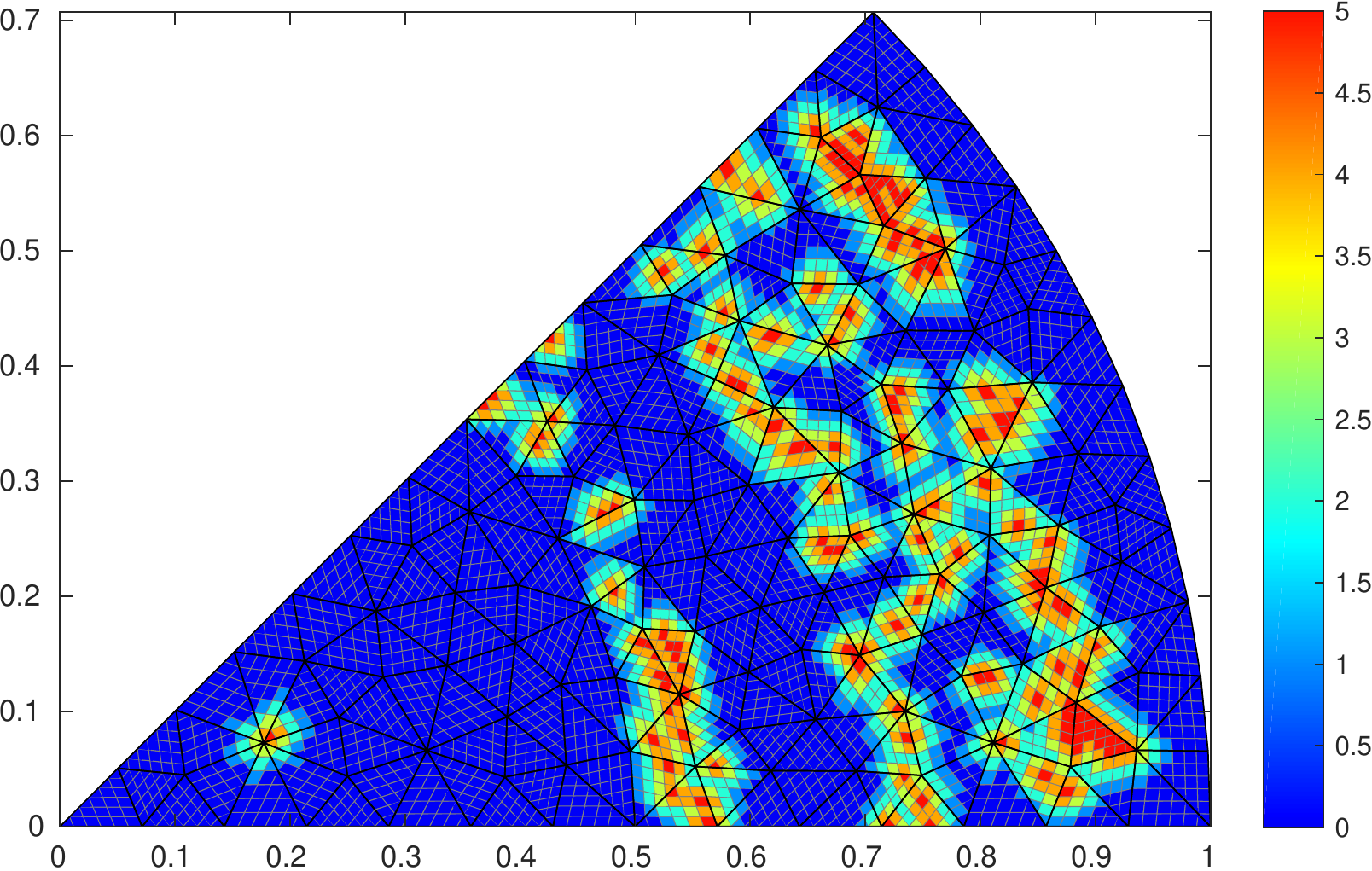}\label{fig_sod_detect}}
    \caption{6th-order APLSC-DG solution for the polar Sod shock tube problem on 230 cells.}
  \label{fig_sod_6th-order}
  \end{center}
\end{figure}

\begin{figure}[!ht]
  \begin{center}
    \includegraphics[height=6.5cm]{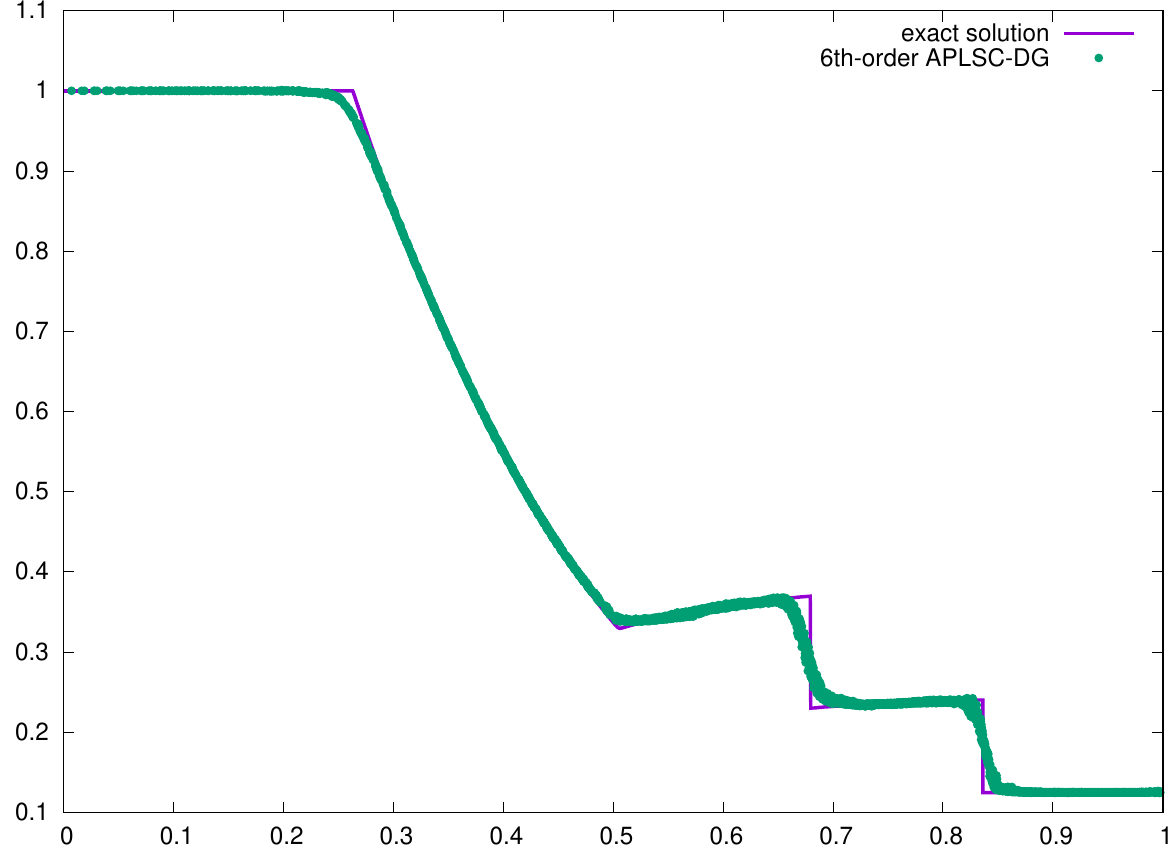}
    \caption{6th-order APLSC-DG solution for polar Sod tube problem on 230 cells:  submean density versus radius.}
  \label{fig_sod_profile}
  \end{center}
\end{figure}

\begin{figure}[!ht]
  \begin{center}
    \subfigure[Cell density values.]{\includegraphics[height=5.5cm]{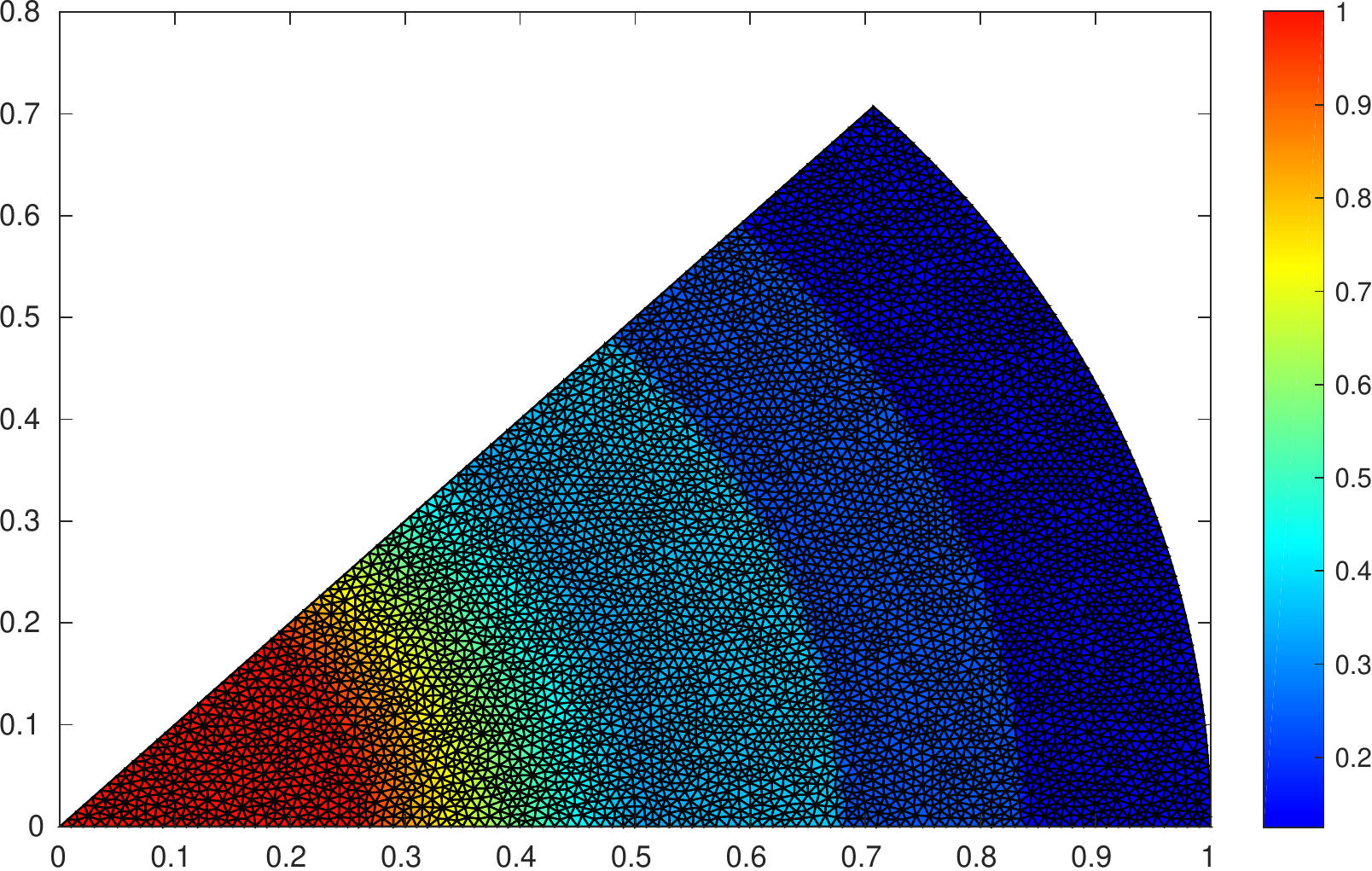}\label{fig_sod_3rd-order_den}}\hspace*{2mm}
    \subfigure[Submean values versus radii.]{\includegraphics[height=5.5cm]{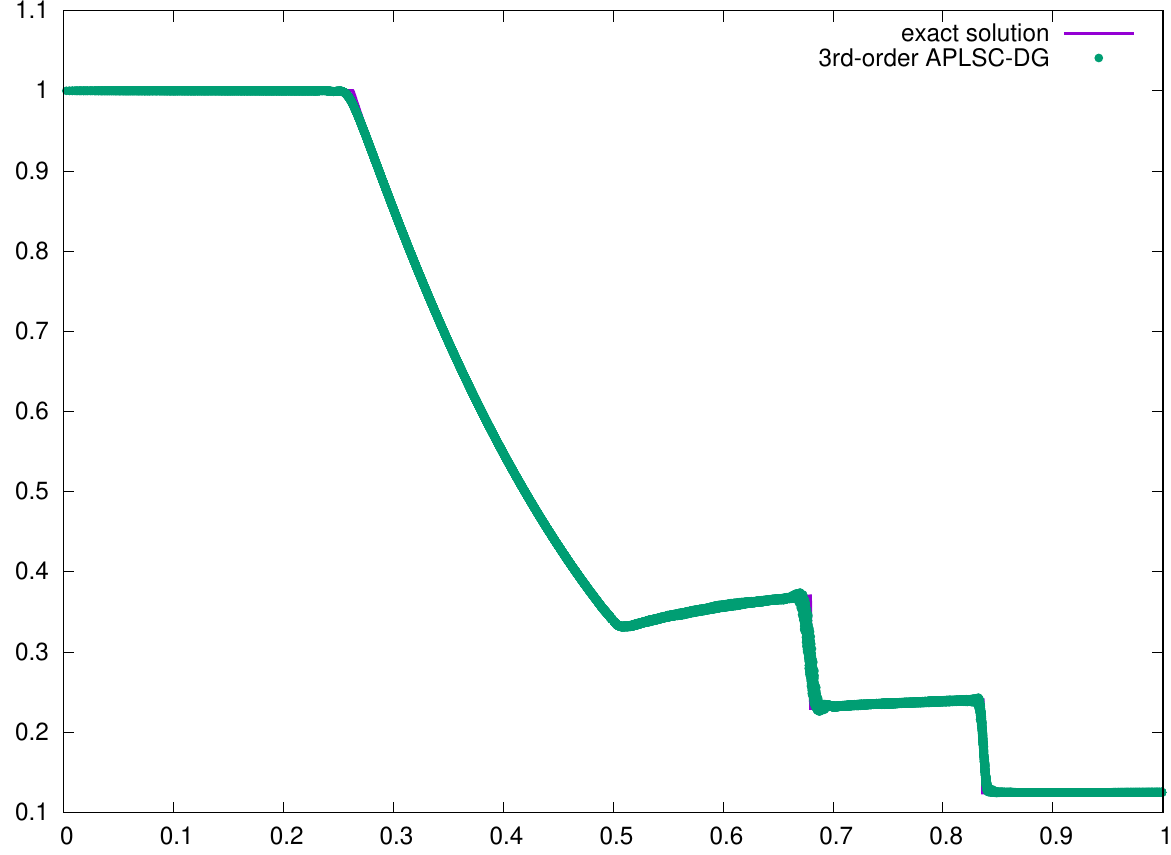}\label{fig_sod_3rd-order_profile}}
    \caption{3rd-order APLSC-DG solution for the polar Sod shock tube problem on 10571 cells.}
  \label{fig_sod_3rd-order}
  \end{center}
\end{figure}

In Figure~\ref{fig_sod_6th-order}, the 6th-order APLSC-DG scheme has been used on a very coarse anisotropic mesh made of 230 triangular cells. In the light of Figure~\ref{fig_sod_den}, one can see how the radial wave structure has been accurately capture, even in this coarse mesh case. Figure~\ref{fig_sod_profile}, where the subcell mean values versus the subcell centroid radii $\sqrt{x^2+y^2}$ confirms this statement as the different points for a given radius do coincide. In Figure~\ref{fig_sod_detect}, subcells corrected during the different Runge-Kutta stages of the last time iteration are colored accordingly the amount of first-order FV correction used. It illustrates how this \textit{a posteriori} correction procedure has been activated on zones corresponding to the solution loss of smoothness, meaning the left and right ends of the expansion fan, the contact discontinuity and the shock. One can also observe how the correction does operate locally inside the cell at a subcell scale, allowing the preservation of DG subcell high accurate resolution.
Let us emphasize that this \textit{a posteriori} correction procedure is not limited to the case of very high-order of accuracy on coarse grids. It also preforms very well at second or third order. See for instance Figure~\ref{fig_sod_3rd-order} where the third-order APLSC-DG numerical solution on 10571 triangular cells is depicted.

Figure~\ref{fig_sod_3rd-order} depicts how the numerical solution produced by our third-order APLSC-DG scheme respects the radial symmetry of the flow, and tends to converge to the exact solution solution in this quite refined mesh context.

\subsubsection{Sedov point blast problem.}
\label{subsubsect_Euler_sedov}

We consider the Sedov problem for a point-blast in a uniform medium. An exact solution based on self-similarity arguments is available, see for instance \cite{Kamm07}. The initial conditions are characterized by $(\rho_0,p_0,\bs{v}_0)=(1,10^{-14},\bs{0})$, and the polytropic index is equal to $\frac{7}{5}$. We set an initial delta-function energy source at the origin prescribing the pressure in a control volume, yet to be defined, containing the origin as follows, $p_{or}=(\gamma-1)\, \frac{\varepsilon_0}{v_{or}}$, where $v_{or}$ denotes the measure of the chosen control volume and $\varepsilon_0$ the total amount of release energy. By choosing $\varepsilon_0=0.244816$, as suggested in \cite{Kamm07}, the solution consists of a diverging infinite strength shock wave whose front is located at radius $r=1$ at $t=1$, with a peak density reaching 6. The computational domain is defined in polar coordinates by $(r,\theta)\in [0,\, 1.2]\times [0,\frac{\pi}{4}]$. Similarly to the polar Sod shock tube problem, we prescribe symmetry boundary conditions at the boundaries $\theta=0$ and $\theta=\frac{\pi}{4}$, and an outflow condition at $r=1.2$.\\

Regarding the control volume in which the delta-function energy will be dropped off, generally the cell containing the origin is considered. Here, to make this test case even more challenging, we choose to restrict the energy source only to the one subcell containing the origin. This means that initially, in one grid element the pressure in one subcell will be set to $p_{or}$, while in the remainder of the cell the pressure will be $10^{-14}$. Let us further emphasize that generally in this test case, because one cannot simulate vacuum, the initial pressure is set to $10^{-6}$ over the domain, except at the origin. Here, to make it once again more challenging, we set the initial pressure to $10^{-14}$.\\

\begin{figure}[!ht]
  \begin{center}
    \subfigure[$t=0$.]{\includegraphics[height=5.cm]{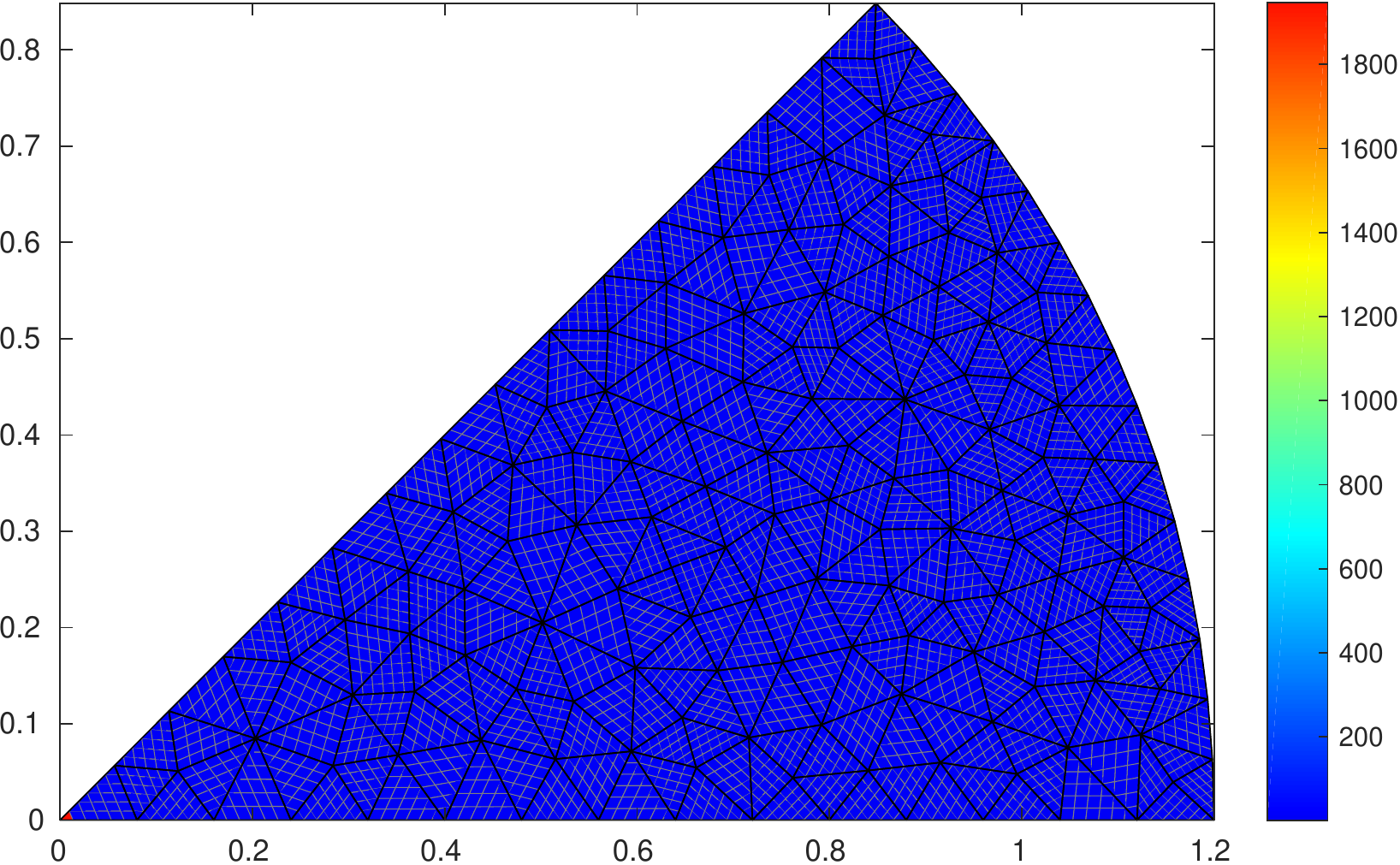}\label{fig_sedov_6th_ene_ini}}\hspace*{2mm}
    \subfigure[$t=1$.]{\includegraphics[height=5.cm]{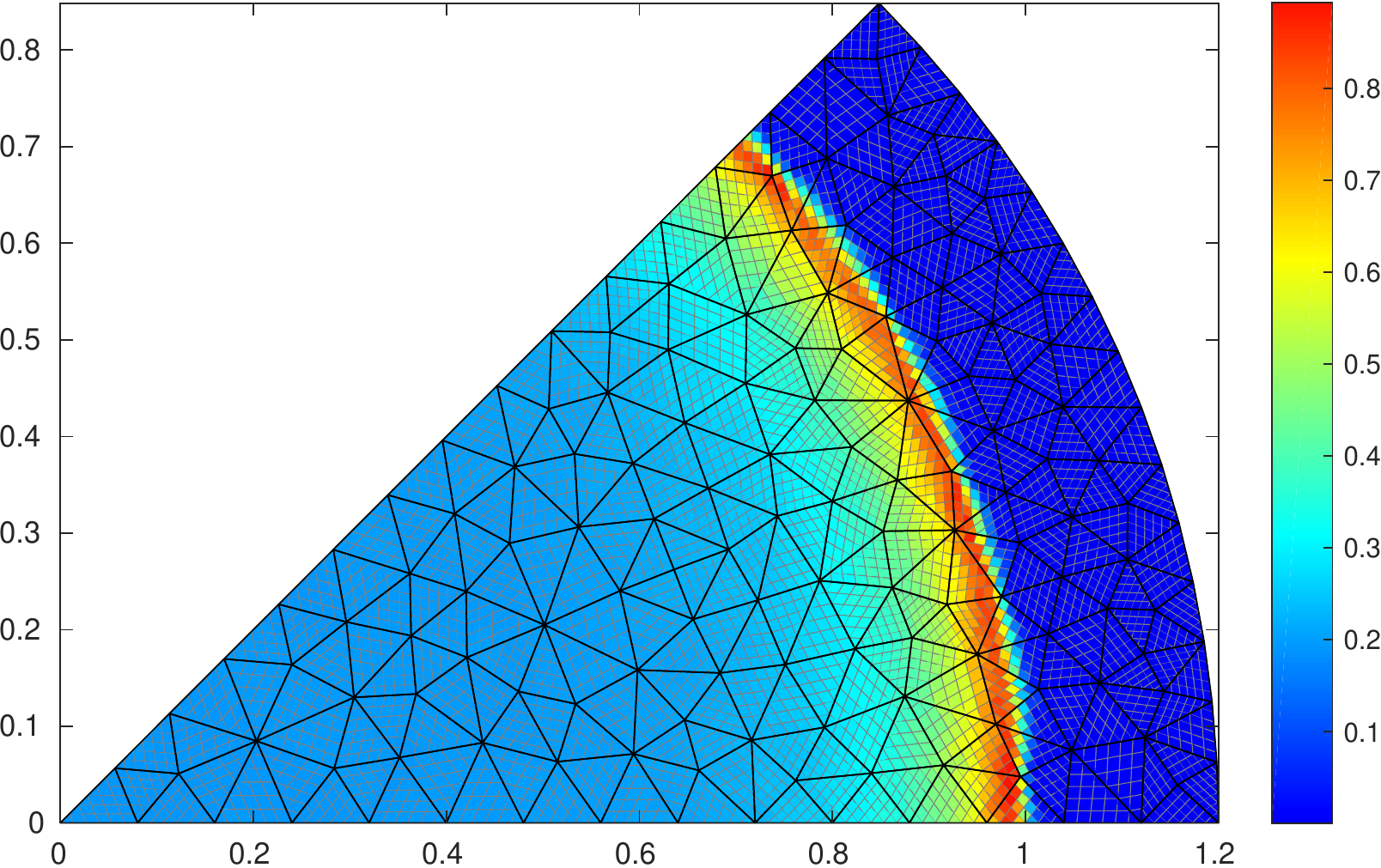}\label{fig_sedov_6th_ene}}
    \caption{6th-order APLSC-DG solution for the Sedov problem on 271 cells: subcell mean total energy values.}
  \label{fig_sedov_6th-order}
  \end{center}
\end{figure}

We run this modified Sedov point blast problem with the sixth-order APLSC-DG scheme on a very coarse grid made of 271 triangular cells. In this particular case, the amount of total energy contained in the subcell located at the origin reaches $1947.5$, while in the rest of the cell as well as in the remainder of the domain the total energy is set to 2.5E-14, see Figure~\ref{fig_sedov_6th-order}. Any scheme lacking positivity-preserving property or a rigorous stabilization technique would fail solving this test problem. 
\begin{figure}[!ht]
  \begin{center}
    \includegraphics[height=6.5cm]{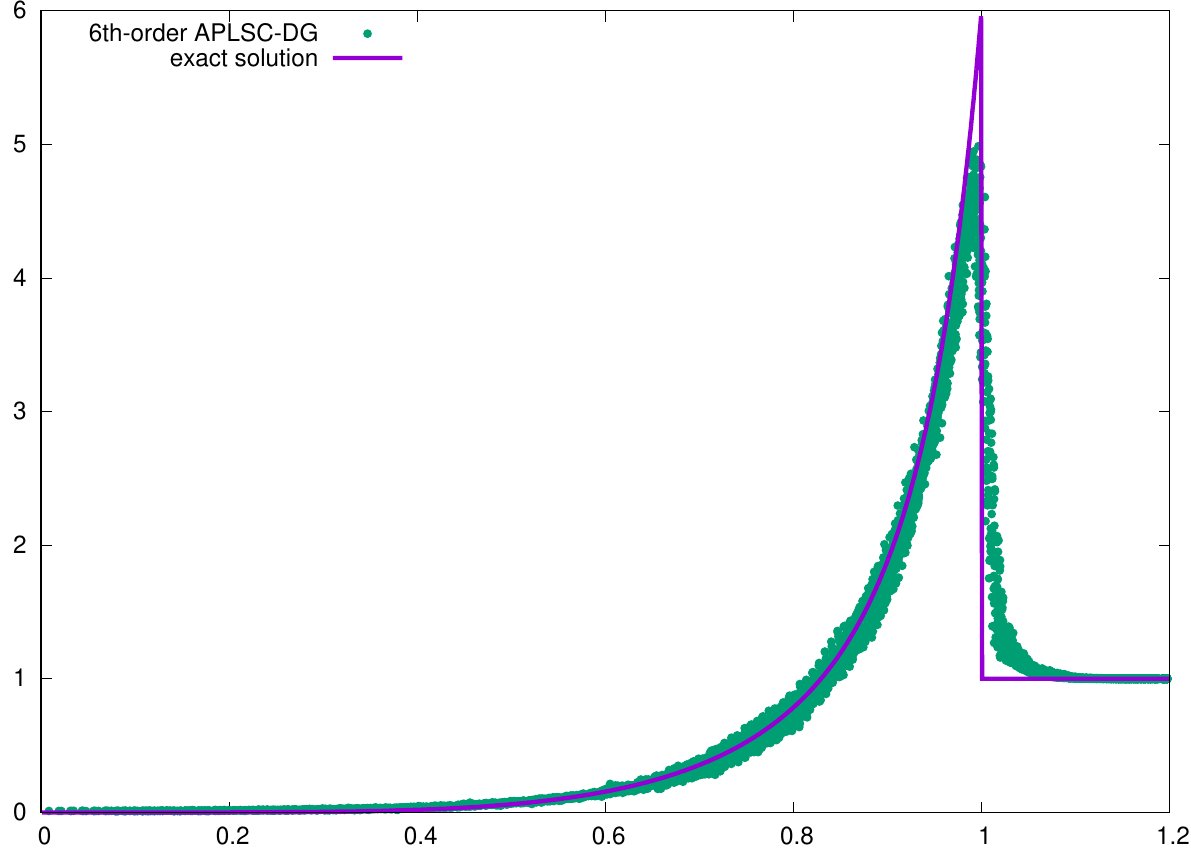}
    \caption{6th-order APLSC-DG solution for the Sedov problem on 271 cells: submean density versus radius.}
  \label{fig_sedov_6th-order_profile}
  \end{center}
\end{figure}
In Figure~\ref{fig_sedov_6th_ene}, one can see how the circular aspect of the shock has been accurately captured by the scheme, and the shock wave front is correctly located. This latter further goes different cells, enlightening the very precise subcell resolution of the APLSC-DG method. The numerical solution produced remains quite close to the one-dimensional self-similar exact solution, see Figure~\ref{fig_sedov_6th-order_profile}.
\begin{figure}[!ht]
  \begin{center}
    \subfigure[Cell mean total energy map.]{\includegraphics[height=5.5cm]{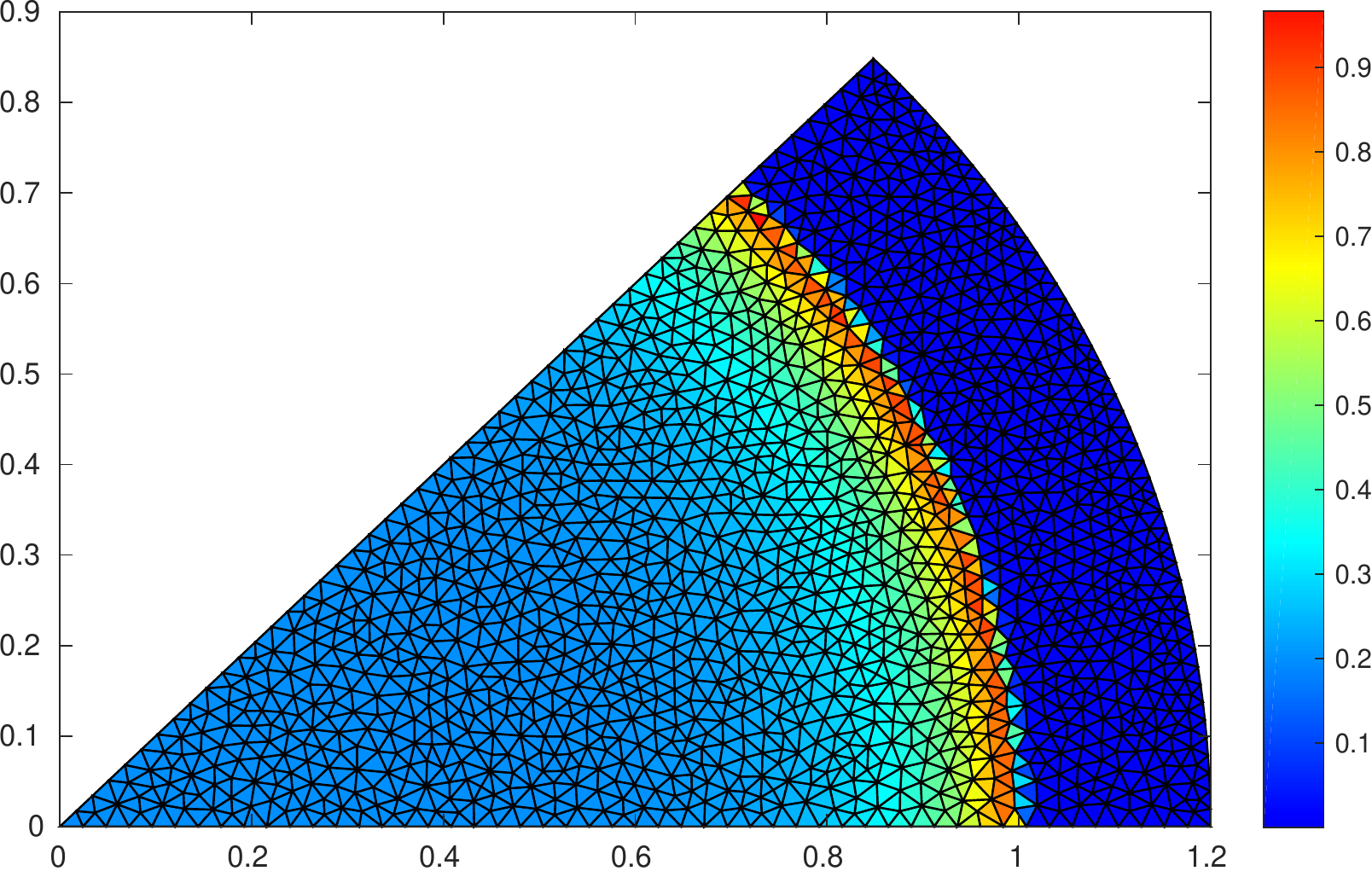}\label{fig_sedov_3rd_ene}}\hspace*{2mm}
    \subfigure[Subcell mean density versus radii.]{\includegraphics[height=5.5cm]{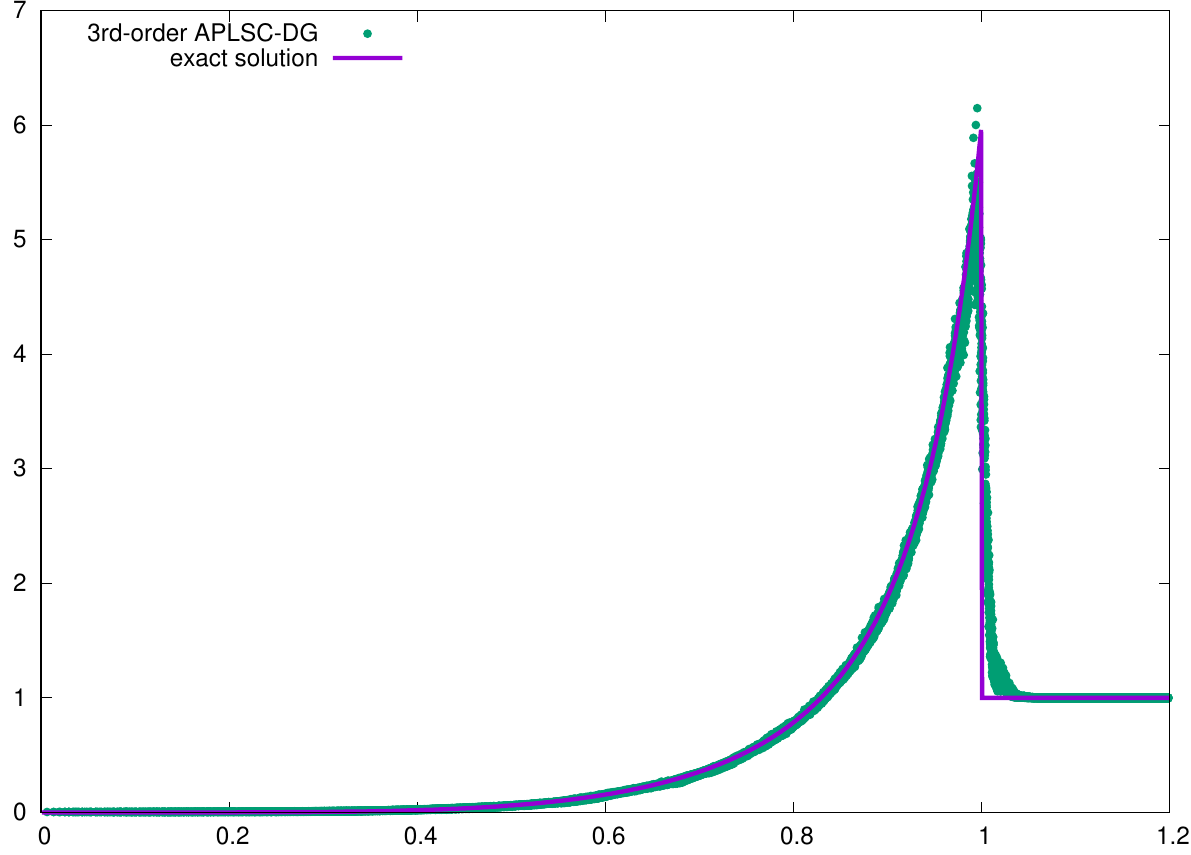}\label{fig_sedov_3rd_den_profile}}
    \caption{3rd-order APLSC-DG solution for the Sedov point blast problem at $t=1$ on 2894 cells.}
  \label{fig_sedov_3rd-order}
  \end{center}
\end{figure}
Once again, this \textit{a posteriori} local subcell correction procedure designed for DG schemes performs also very well for lower order methods, as depicted in Figure~\ref{fig_sedov_3rd-order} where a finer grid made of 2894 cells has been used. Indeed, the numerical solution is very close to the one-dimensional analytical solution. In Figure~\ref{fig_sedov_3rd_ene}, only the cell total energy means values are represented, and not the submean values as we generally do, for a better readability of the results in this finer grid context.

\vspace*{5mm}
\subsubsection{The forward-facing step problem.}
\label{subsubsect_Euler_step}

\begin{figure}[!ht]
  \begin{center}
    \includegraphics[width=16.5cm]{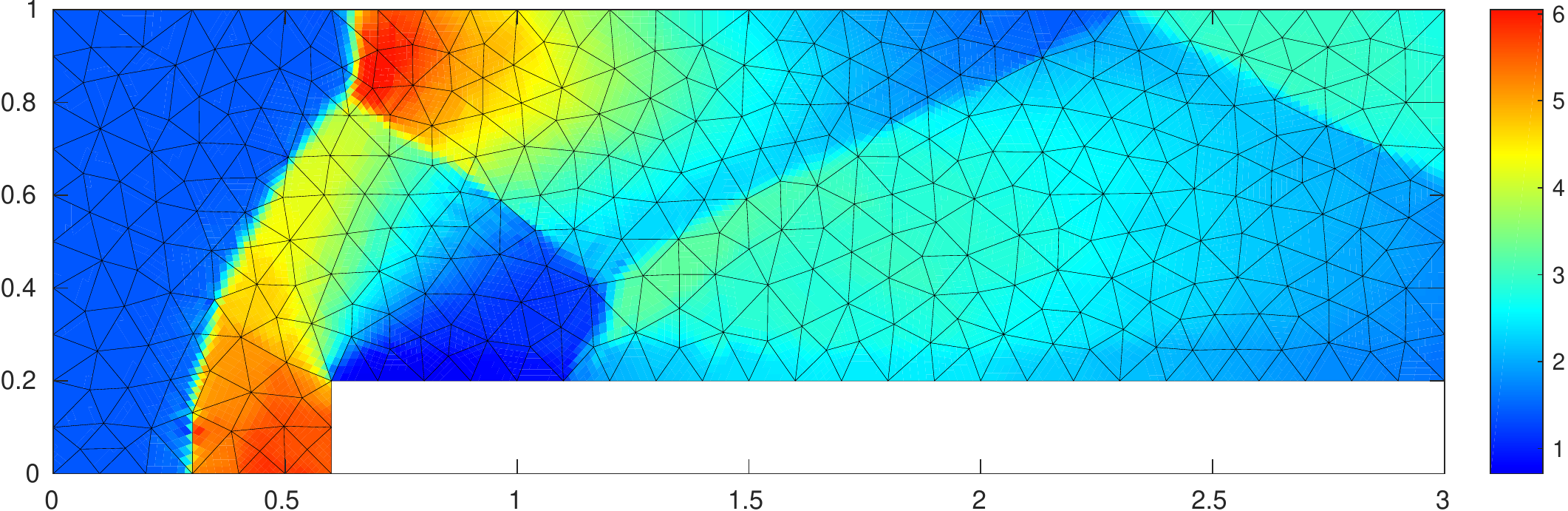}
    \caption{6th-order APLSC-DG solution for the facing step problem on 680 cells at $t=4$: submean density map.}
  \label{fig_step_6th-order}
  \end{center}
\end{figure}

We now consider the forward facing step problem, which has been initially introduced by A. Emery in \cite{emery_1968}, and further studied by P. Woodward and P. Colella in \cite{woodward_1984}. This challenging test case consists in a Mach 3 flow in a 3 units in length and 1 unit in width wind tunnel. Initially, the tunnel is filled with a gamma gas law with $\gamma=\frac{7}{5}$, which everywhere has density $\rho_0=1.4$, pressure $p_0=1$ and velocity $\bs{v}_0=\(3,\, 0\)\tra$. The 0.2 high step being located at $x=0.6$, the computational domain is then $\([0,\,3]\times[0,\,1]\)\setminus\([0.6,\,3]\times[0.2,\,1]\)$. Gas with this density, pressure and velocity is continually fed in from the left-hand boundary. Let us emphasize that unlike as it is generally done, we did not refine the mesh near the corner, see Figure~\ref{fig_step_6th-order}, nor modify in any way our APLSC-DG scheme.
In Figure~\ref{fig_step_6th-order}, the numerical solution obtained by means of 6th-order APLSC-DG scheme on an unstructured grid made of 680 cells is displayed. Let us note that despite the coarseness of the mesh used, the shocks and the rarefaction fan created around the corner are quite well resolved, while ensuring an low oscillatory behavior. This result demonstrates once again the high capability of the presented \textit{a posteriori} local subcell correction combined with high-order discontinuous Galerkin methods.

\section{Conclusion}
\label{conclusion}

The paper aims at presenting the two-dimensional unstructured extension of the new correction technique of DG schemes introduced in \cite{vilar_aplsc_1D}. This \textit{a posteriori} procedure relies on the expression of DG methods as a FV-like scheme on a subgrid. By means of this theoretical part, we modify at the subcell level the so-called reconstructed fluxes only where the uncorrected DG scheme has failed. Consequently, only very few subcells require this particular treatment. In this paper, a new version of the correction procedure is also introduced, where a convex blending of high-order DG reconstructed fluxes and first-order FV fluxes is applied in the vicinity of troubled zones. For the remaining subcells, the submean values obtained through the uncorrected DG method are kept, as they have been detected as admissible by troubled zone criteria. This correction procedure allows us to retain the very precise subcell resolution of DG schemes, along with addressing the issues of spurious oscillations or non-entropic behavior. A wide number of test cases on different problems have been used to depict the good performance and robustness of the presented correction technique. Different types of cell subdivision have been compared.\\

In the future, we intend to extend this \textit{a posteriori} correction technique to moving grid configurations, both ALE and Lagrangian formalisms, as well as the case of curvilinear meshes. We also plan to adapt this local subcell reconstructed flux correction framework to the \textit{a priori} paradigm, by means of the Flux-Corrected Transport (FCT) methodology, in order to obtain an automatic very high-order and property preserving scheme.

\section*{Acknowledgment}
R.A. was partially funded  by SNF project "Structure preserving and fast methods for hyperbolic systems of conservation laws" number 200020$\_$204917.



\bibliographystyle{plain}
\bibliography{APLSC_2D}

\end{document}